%% file: publication.tex
\documentclass{article}[11pt]
\usepackage[papersize={8.5in,11in},top=4.0cm,bottom=2.0cm,left=2.0cm,right=2.0cm]{geometry} % right margin was on 2.0cm, width was on 8.5in, changed for review by Gottfried
\usepackage{subcaption}
\usepackage[square,numbers]{natbib}
\usepackage{amsmath,amsfonts,mathrsfs,amsfonts,amsthm}
\usepackage{graphicx}
\usepackage{epstopdf}
\usepackage{dsfont}
\usepackage{color}
\usepackage{mathtools}
\usepackage{setspace}
\usepackage{float}
\usepackage{floatflt}
\usepackage{booktabs}
\usepackage{siunitx}
\newcommand{\trp}{^{\textnormal{\scriptsize T}} }

\usepackage{tikz}
\usepackage{pgf}
\usepackage{pgfplots}
\usepackage[utf8]{inputenc} \DeclareUnicodeCharacter{2212}{-}

\pgfplotsset{compat=newest}
\usetikzlibrary{arrows.meta}
\usetikzlibrary{pgfplots.groupplots}
\usetikzlibrary{matrix}
\usetikzlibrary{cd}

\theoremstyle{plain}
\newtheorem{thm}{Theorem}[section]
\newtheorem{prop}[thm]{Proposition}

\newtheorem{cor}[thm]{Corollary}
\newtheorem{lem}[thm]{Lemma}
\theoremstyle{remark}
\newtheorem{rem}[thm]{Remark}
\newtheorem{exm}[thm]{Example}

\usepackage{caption}
\usepackage{nicefrac}
\usepackage{xargs}
\usepackage{xcolor}

\DeclareMathOperator{\rank}{\mathrm{rank}}
\DeclareMathOperator{\id}{\mathds{1}}
\usepackage{fancyhdr}
\usepackage{color}
\definecolor{light}{gray}{0.50}
\definecolor{heavy}{gray}{0.35}
\definecolor{black}{gray}{0.0}
\definecolor{dgreen}{rgb}{0.0,0.7,0}
\definecolor{dred}{rgb}{0.9959,0,0}
\definecolor{green}{rgb}{0.0,0.99599,0.0}
\definecolor{purple}{rgb}{0.6,0.0,0.4}

\newcounter{kleincommentno}
\def\imgPath{image/}
\graphicspath{{image/}}
\begin{document}

\title{Balanced data assimilation for highly oscillatory mechanical systems}

\author{Gottfried Hastermann\thanks{Freie Universit\"at Berlin, Institut f\"ur Mathematik, Arnimallee 6, D-14195 Berlin, Germany} \and
  Maria Reinhardt\thanks{Universit\"at Potsdam, Institut f\"ur Mathematik, Karl-Liebknecht-Str. 24/25, D-14476 Potsdam, Germany} \and
  Rupert Klein\thanks{Freie Universit\"at Berlin, Institut f\"ur Mathematik, Arnimallee 6, D-14195 Berlin, Germany} \and
  Sebastian Reich\thanks{Universit\"at Potsdam, Institut f\"ur Mathematik, Karl-Liebknecht-Str. 24/25, D-14476 Potsdam, Germany}
}

\maketitle

\begin{abstract}
  Data assimilation algorithms are used to estimate the states of a dynamical
  system using partial and noisy observations. The ensemble Kalman filter has
  become a popular data assimilation scheme due to its 
  simplicity and robustness for a wide range of application areas. Nevertheless, 
  this filter also has limitations due to its inherent assumptions of 
  Gaussianity and linearity, which can manifest themselves in the form of dynamically 
  inconsistent state estimates. This issue is investigated here for balanced, slowly 
  evolving solutions to highly oscillatory Hamiltonian systems which are prototypical 
  for applications in numerical weather prediction. It is demonstrated that the 
  standard ensemble Kalman filter can lead to state estimates that do not satisfy the 
  pertinent balance relations and ultimately lead to filter divergence. Two remedies 
  are proposed, one in terms of blended asymptotically consistent time-stepping schemes, 
  and one in terms of minimization-based post-processing methods.
  The effects of these modifications to the standard ensemble Kalman filter are 
  discussed and demonstrated numerically for balanced motions of two 
  prototypical Hamiltonian reference systems.
\end{abstract}
\vspace{1cm}

\vfill
\newpage

\doublespacing{}

% ==============================================================================
% ==============================================================================
% ==============================================================================

\section{Introduction}
A problem dating back as far as the advent of numerical weather prediction is
the incorporation of physical
observations into a dynamical model with more than one time scale.
The famous first forecast of L.\ F.\ Richardson~\cite{Lynch2014} 
failed due to the choice of an unbalanced
initial condition gained from observations.
In essence the observational data did not satisfy certain discrete energy
balances and this triggered artificial oscillations in the pressure, ultimately 
leading to erroneous results. In the context of data assimilation, several 
solutions to the related problem of finding balanced initial data were proposed 
over the last decades.
Lynch~\cite{Lynch1992} suggested to apply a digital filter after every
assimilation step to eliminate spurious fast oscillations, and this 
technique was adopted in the weather prediction context with some success. 
Strategies that incorporate the observational
data in the model evolution in a gradual and smooth way instead of using all
the information about the observation at one single point in time have been
suggested for example in~\cite{Bloom1996} and~\cite{Kay2010}. Kepert~\cite{kepert2009}
proposed a method to overcome the issue of artificial balances triggered by localized 
Bayesian data assimilation. Here ``localization'' refers to approaches designed to 
avoid spurious long-range correlations by allowing the model state in a given grid 
point of a flow simulation to be influenced only by data found within a given 
maximum distance from it. Kepert suggested to localize the data filter in the 
streamfunction and velocity potential fields rather than in the velocity or momentum 
variables. With a similar goal, Gottwald~\cite{gottwald2014} incorporated additional 
climatological information in the assimilation process in the sense of 
variance-limited Kalman filters, \cite{GottwaldEtAl2011}, so as to drive the 
model's level of imbalance towards its climatological mean.

In the context of variational data assimilation for slow-fast 
Hamiltonian systems the issue was addressed, e.g., by Cotter~\cite{Cotter2013}. 
Variational methods, in contrast to the ``filtering techniques'' referred to
in the last paragraph, aim to optimize the match between simulation and 
observation over an entire time window in the past, i.e., not only just at the
time when the observations arrive. Cotter's approach differs 
from the ones cited above in that he explicitly uses 
an analytical fast-slow transformation of variables which he assumes can be 
derived from the structure of the system's Hamiltonian and which defines its 
relevant slow manifold. In this situation, he defines artifical Hamiltonian 
dynamics that drives a system state from arbitrary states towards close-by 
states on the slow manifold relatively quickly, but still on the slow time scale. 
He then suggests to use this artificial dynamics to constrain a variational 
data assimilation method (4DVAR) so as to produce a balanced state as the 
initial condition for the next forecast at the end of the data assimilation 
time window.

Here we propose two alternative approaches to addressing the balancing
problem for data assimilation based on filtering techniques. The first 
approach relies on the ensemble-based Bayesian sequential data assimilation 
paradigm and is designed as a post-processing step in the filtering procedure 
that penalizes imbalances and is structurally similar to the 3DVAR method, 
\cite{kalnay_atmospheric_2002}. Two alternatives for the computational 
implementation of this post-processing step are briefly discussed, comprising a 
Gauss-Newton minimization and a pseudo-time evolution.

The second approach proposed in this paper is, in contrast, incorporated in the 
forward simulation and relies on the ability of our asymptotically consistent
numerical method to seamlessly represent balanced and unbalanced dynamics.
In some aspects, this approach is similar to that of Cotter~\cite{Cotter2013} 
in that we use analytical knowledge regarding the full oscillatory and a nearby 
reduced slow dynamics. Yet, our ansatz is not tied closely to Hamiltonian
structure, and our balancing strategy is a direct part of the forward 
simulation rather than being incorporated in the data assimilation procedures.
The key idea, first formulated in \cite{BenacchioEtAl2014}, is to first apply 
some known filtering technique for data assimilation that is not specifically
designed to maintain physical balances, but to start the subsequent forward
simulation by several time steps with a ``blended model'' that interpolates
in a judicious way between the reduced slow and the full oscillatory 
dynamics. The discretization of the interpolating model family is dissipative
with respect to the fast modes, while properly advancing the balanced modes.
As a consequence, when the full dynamics becomes active after the blending
time window, fast oscillations have been removed and nearly balanced 
conditions prevail throughout the remaining majority of the forward 
simulation time steps.

The rest of this section introduces the class of highly-oscillatory 
nonlinear finite dimensional test problems utilized in this paper, and discusses 
the failure of a Bayesian data assimilation procedure which has motivated our 
work. Section~\ref{sec:ProposedMethods} describes the two balanced data 
assimilation approaches proposed in this paper. Section~\ref{sec:results} 
compares the performance of several data assimilation techniques for the 
oscillatory test problem. Section~\ref{sec:conclusions} provides a summary
and an outlook to future work.

% ==============================================================================
% ==============================================================================
% ==============================================================================

\subsection{Model problem}\label{sec:model}

With atmospheric models in mind as a motivation, we propose numerical techniques 
that allow the user to follow the slow evolution of a system with multiple time
scales starting from balanced initial data. In doing so, we 
restrict to the finite-dimensional setting in this paper, in line with 
Lorenz' seminal investigations of oscillatory systems and predictability in
\cite{Lorenz1963,Lorenz2006} and, more specifically, with studies into 
the existence and properties of slow manifolds for multiple time scale
Hamiltonian systems in \cite{Bokhove1996,Lynch2002,Camassa1995}. 

In particular, we discuss sequential data assimilation for 
highly-oscillatory systems with Hamiltonian energy functional
\begin{equation}\label{eq:Hamiltonian}
  H^\varepsilon (q,p)	= \frac{1}{2} p\trp p + \frac{1}{2\varepsilon^2}
  g(q)\trp K g(q) + V(q) ,
\end{equation}
with momenta and coordinates \(p,q\in \mathbb{R}^N\). Here \(V:
\mathbb{R}^{N} \rightarrow \mathbb{R}\) is a potential energy, \(g :\mathbb{R}^N \to
\mathbb{R}^L, \, L\leq N\) gives rise to rapid oscillations with a diagonal
matrix of force constants \(K = \text{diag}(k_1, ..., k_L), k_i > 0 \), and
\(\varepsilon\) is a stiffness parameter satisfying \(0 < \varepsilon \ll 1\). 
The associated Hamiltonian
equations of motion are then given by
\begin{equation}\label{eq:model1}
  \begin{aligned}
    \dot{q} & = p
    \\
    \dot{p} & = - \varepsilon^{-2} G(q)\trp K g(q) - \nabla V(q),
  \end{aligned}
\end{equation}
where \(G(q) := Dg(q) \in \mathbb{R}^{L\times N}\) denotes the Jacobian matrix of
\(g\) at \(q\).
These equations pose challenges in their numerical treatment as well as for
sequential data assimilation techniques in the limit \(\varepsilon \to 0\).
We observe that solutions of~\eqref{eq:model1} preserve the Hamiltonian
energy functional~\eqref{eq:Hamiltonian}, and that bounded energy, i.e., \(
H^\varepsilon (q,p) = {\cal O}(1)\) as \(\varepsilon \to 0\), implies \(g(q)
= {\cal O}(\varepsilon)\). In other
words, for \(\varepsilon \ll 1\) solutions \(q\) of bounded energy have to
stay close to the constraint manifold
\begin{equation}
  {\cal M} = \{ q\in \mathbb{R}^N: \, \|g(q)\| = 0 \}.
\end{equation}
From here on we will assume that \(\text{rank}(G(q)) = L < N\) within the
domain of interest,  so that an
explicit local decomposition of \(q\) into fast and slow modes is possible 
according to the following
\begin{rem}:
  Let \(\Omega\subseteq \mathbb{R}^N\) be open and bounded and let
  \(G(q) \in \mathbb{R}^{L\times N}\) with
  \(\rank G(q) = L < N\) for all \(q \in \Omega \), then the linear map
  \(\mathcal{P}_q: \mathbb{R}^N\rightarrow\mathbb{R}^N\) given by
  \begin{align}\label{eq:ProjectionOperator}
    \mathcal{P}_q \coloneqq G\trp(q) {(G(q)G\trp(q) )}^{-1}G(q)
  \end{align}
  is an orthogonal projection.
\end{rem}
\begin{rem}:
  \(\mathcal{P}_q^{\perp}\) denotes the orthogonal projection onto the
  orthogonal
  complement of the image of \(\mathcal{P}_q\).
  For every \(q \in \mathcal{M}\) its image is included in the corresponding
  tangent space to \(\mathcal{M}\),
  % i.e. \(\mathcal{P}_q^\perp \mathbb{R}^N \rightarrow T_q\mathcal{M}\).
  in fact \(\mathcal{P}_q^\perp \mathbb{R}^N = T_q\mathcal{M}\). 
\end{rem}
\begin{rem}:
  As presented in~\cite{BenettinEtAl1987} we can decompose the Hamiltonian
  energy functional in~\eqref{eq:Hamiltonian} into fast, slow and coupling
  energies after introducing a local coordinate transform into slow and fast
  variables (c.f.~Lemma \eqref{lem:split} below). The ``fast'' 
  part of the Hamiltonian describes rapid oscillations in the fast 
  variables orthogonal to the slow manifold~\(\mathcal{M}\) and depends on 
  the slow variables, i.e., on the current nearest point on~\(\mathcal{M}\),
  only parametrically. In the original coordinates this fast or 
  ``oscillatory'' part of the Hamiltonian reads
  \begin{equation}\label{eq:oscHamiltonian}
    H_{\rm osc}^\varepsilon (q,p) = \frac{1}{2} p\trp \mathcal{P}_q p +
    \frac{1}{2\varepsilon^2}
    g(q)\trp K g(q).
  \end{equation}
  %
%  %
%  \begin{equation}\label{eq:HamiltonianSplit}
%    H^\varepsilon (q,p) = \frac{1}{2} p\trp \mathcal{P}_q p +
%    \frac{1}{2\varepsilon^2}
%    g(q)\trp K g(q) + \frac{1}{2} p\trp (1-\mathcal{P}_q) p + V(q).
%  \end{equation}
%  %
  The above-mentioned coupling terms vanish for \(\varepsilon \rightarrow
  0\) in the present scenario of small amplitude oscillations, for diagonal 
  positive definite $K$, and for eigenfrequencies of the individual
  components of the fast oscillations that are independent of the 
  position on the manifold, see~\citep{RubinUngar1957,Bornemann1997}.
  As a consequence, 
  the trajectories given by the full Hamiltonian in~\eqref{eq:Hamiltonian} follow
  the evolution determined by~\eqref{eq:oscHamiltonian} closely for short times. 
  %From a different point of view we observe that in the correct local coordinate
  %system, the motion orthogonal to the constraint manifold is almost harmonic.
  %
\end{rem}
Henceforth we will assume \(g \) to be locally smooth and \(G(q)\) to have
full rank \(L\) for all \(q\in \mathbb{R}^N\) satisfying \(\|g(q)\| \leq C\)
for sufficiently large constant \(C>0\). To state the setting more
rigorously, we consider solutions \(q^\varepsilon,p^\varepsilon \in
C^1([0,T],\mathbb{R}^N)\) to~\eqref{eq:model1} given the initial conditions
\begin{equation} \label{eq:initcond}
  \begin{aligned}
    q^\varepsilon(0) & = q^0_0 + \varepsilon \bar{q}, \quad \bar{q} \in
    \mathbb{R}^{N}                                                         \\
    p^\varepsilon(0) & = p^0_0 + \varepsilon \bar{p}, \quad  \bar{p} \in
    \mathbb{R}^{N}
  \end{aligned},
  \quad
  (q_0^0,p_0^0) \in \mathcal{TM},
\end{equation} 
where \(\mathcal{TM}\) denotes the tangential bundle of \(\mathcal{M}\) which
we will interpret as a manifold in phase space, i.e.,
\begin{equation}
  \mathcal{TM} = \{(q,p) \in \mathbb{R}^{2N}: q\in \mathcal{M} \land p \in
  T_q\mathcal{M} \}\,.
\end{equation}
Therefore the initial data is, up to a
perturbation of order \(\varepsilon\),
\textit{tangential}~\cite{RubinUngar1957}, and we note in passing that the oscillatory energy \eqref{eq:oscHamiltonian} is of order ${\cal O}(1)$ in this case.
Rubin and Ungar~\cite{RubinUngar1957} proved existence of a unique
solution \(q^0,p^0\) to the differential algebraic equation system
\begin{equation}\label{eq:limiteq}
  \begin{alignedat}{3}
    \dot{q}^0 & = p^0					 & \qquad q^0(0) & =
    q^0_0 \in \mathcal{M}
    \\
    \dot{p}^0 & = - G\trp(q^0) K \lambda - \nabla V(q^0) & \qquad p^0(0) & =
    p^0_0 \in T_{q^0_0}\mathcal{M}
    \\
    g(q^0)   & = 0\,
  \end{alignedat}
\end{equation} 
and convergence \(q^\varepsilon(t) \overset{\epsilon \rightarrow
  0 }{\longrightarrow} q^0(t)\),
\(p^\varepsilon(t) \overset{\epsilon \rightarrow 0 }{\longrightarrow} p^0(t)\)
uniformly for \(t\in [0,T]\) in the ``tangential case'' $(\bar q = 0, \bar p = 0)$.
They obtained a similar result for the non-tangential case when the eigenfrequencies of
the fast oscillations normal to the constraint manifold are independent of position. 
This is the scenario we pursue in our computational examples below. See 
\cite{Bornemann1997} for an alternative elegant rigorous analysis in the case of 
co-dimension~1. The Lagrange multiplier \(\lambda \in C^1([0,T],\mathbb{R}^L)\) can
the be determined algebraically for every \(t \in [0,T]\) by 
\begin{equation} \label{equ:lagrangian}
  0 = \ddot{g}(q^0) = -G(q^0) \left[\nabla V(q^0) + G\trp(q^0) K \lambda \right]
  + \sum_{i=1}^N \sum_{j=1}^N
  \left({(p^0)}\trp e_{i,j} p^0\right) \,  \frac{\partial^2 g(q^0)}{\partial q_i
    \partial q_j},
\end{equation}
where \(e_i\) is the \(i\)-th cartesian unit vector of \(\mathbb{R}^N\) and
\(e_{i,j} = e_i \otimes e_j \in \mathbb{R}^{N\times N} \). Therefore this
differential algebraic system is of index \(3\) and we obtain the additional
hidden constraint \( G(q^0)p^0=0 \) by differentiation of \(g(q^0)=0\) with
regard to the parameter. Under additional assumptions~\cite{BenettinEtAl1989}
could prove that solutions \((q^\varepsilon,p^\varepsilon)\), initially
\(\varepsilon\)-close to \((q^0,p^0)\), stay \(\varepsilon\)-close for
exponentially long times.

\begin{exm}\label{ex:doublependulum}
  Consider a chain of \(L\) mass points with positions \(r_i\in \mathbb{R}^D\)
  and momenta \(v_i\in \mathbb{R}^D\) for \(i = 1, ..., L\).
  The first point (denoted by subscript $0$) is assumed to be fixed, all points have equal mass and are
  under the influence of a constant unidirectional
  force \(a_0 e_D\) with $a_0 \in \mathbb{R}$ and $e_D$ the unit vector in the last cartesian direction in $\mathbb{R}^D$. (Here we think of the gravitational force.)
  All points are pairwise connected by \(L\) (linear) elastic bonds,
  characterized by
  their force coefficients \(K = \mathrm{diag}(k_1,\dots k_L)\) and their
  equilibrium lengths \(l_i>0\) for \(i= 1\dots L\).
  Observing \(N =D L\), we can describe the evolution of this mechanical system
  by~\eqref{eq:model1} if we collect the components of all positions and momenta in \(q\) and \(p\), respectively.
  By means of classical mechanics we then conclude
  \begin{equation}
    \label{eq:pendulumchain}
    \begin{aligned}
      g(q)  =
              & \begin{pmatrix}\|r_1 \| - l_1 & \dots & \| r_{i}  - r_{i-1}\| -
        l_i            & \cdots & \| r_{L}	- r_{L-1}\| - l_L
      \end{pmatrix}\trp \\
      V(q)  = & a_0 \sum_{i=1 }^L	e_D\trp \, r_i,
    \end{aligned}
  \end{equation}
  %
  % where in this case \(e_D\in \mathbb{R}^D\) is the unit vector in the last component.
  If not stated differently, we assume \(D=2\).
  \begin{figure}[htt]
    \centering
    \input{\imgPath doublependulum.eps_tex}
    \caption{A graphical depiction of~\eqref{eq:model1}
      using~\eqref{eq:pendulumchain} for the case of two mass points \(L=2\)
      moving
      in the plane \(D=2\).}
  \end{figure}
  Note that this is a genuinely nonlinear model with at least two time scales
  even for the simplest case of the elastic pendulum with \(L=1\).
  Due to the work of~\cite{BenettinEtAl1987}, we expect solutions
  \(q^\epsilon \) for this case to stay close to the
  solutions of the classical pendulum \(q^0\) for long times.
  Since our goal is to observe and predict chaotic slow dynamics, we wish
  for a system with solutions \(q^0 \), which already exhibit chaotic behavior.
  For this reason we have to consider a slightly more complex model and choose
  \(L=2\).

  As the constraint manifold is of codimension \(L=2\) the results of~\cite{RubinUngar1957} are
  not applicable. Nevertheless one can, under additional assumptions,
  conclude convergence of the solution \(\lim_{ \varepsilon \rightarrow 0}
  q^\epsilon =q^0\) by the means of~\cite{Takens1980,Bornemann1997}. For
  \(L>1\) and non-resonant configurations~\cite{BenettinEtAl1989} proved
  solutions \(q^\varepsilon\) to stay close to \(q^0\) for long times, if
  initially so.
\end{exm}
\begin{rem}
  As demonstrated in~\cite{Reich1995}
  \begin{equation} \label{eq:soft_constraint}
    g(q^{\varepsilon}) = \varepsilon^2
    \lambda(q^{\varepsilon},p^{\varepsilon}),
  \end{equation}
  with \(\lambda (q,p) \) determined by~\eqref{equ:lagrangian}, is a
  better approximation to the slow dynamics of~\eqref{eq:model1} than the
  zeroth order balance relation \((q,p) \in T{\cal M}\). Replacing the
  constraint in~\eqref{eq:limiteq} by~\eqref{eq:soft_constraint} leads to the
  concept of soft or flexible constraints as introduced
  in~\cite{Reich1995,Zhou2000}.
\end{rem}

\begin{rem}
  It should be noted that initial conditions with unconstrained momentum of
  the form
  \begin{equation}
    \begin{alignedat}{4}
      q^\varepsilon_0 & =q^0_0+\varepsilon \tilde{q}_0 &  & \qquad q^0_0\in
      \mathcal{M}, \tilde{q}_0 \in \mathbb{R}^N
      \\
      p^\varepsilon_0 & =\tilde{p}_0		       &  & \qquad \tilde{p}_0
      \in \mathbb{R}^N\,.
    \end{alignedat}
  \end{equation}
  do not generally follow this theory. In this case, an additional
  force term can appear in the limiting equations~\eqref{eq:limiteq}.
  See~\cite{RubinUngar1957,Takens1980,Bornemann1997} for more details.
\end{rem}

% ===================================================================================
% ===================================================================================

\subsection{Bayesian data assimilation}

When describing physical processes by models there are several sources of
uncertainties, such as model errors or an uncertainty about the initial
conditions. Ensemble-based data assimilation combines model outputs with, 
possibly also somewhat erroneous, observational data to estimate a probability 
distribution over model states conditioned on the observations. Owing to the
usually high dimension of the model state space, probability densities are 
generally approximated by the empirical probability densities represented by 
ensembles of individual realizations of model runs.

Variational data assimilation seeks estimates of the trajectory of model 
states over an entire observation time interval by solving a related error 
minimization problem. A well known candidate of this kind is \(4\)D-VAR, as 
explained e.g.\ in~\cite[p.~186]{reichbook}.

An alternative is sequential Bayesian data assimilation. In a forecast step 
this method evolves the empirical probability distribution by forward 
simulation of the ensemble members until a new observation \(y_{\rm obs}\) 
becomes available. The resulting distribution is called the prior or forecast 
distribution, \(\pi^{\rm f}\). In the analysis step, the observations are 
assimilated by applying Bayes' theorem to generate an improved posterior or 
analysis distribution, \(\pi^{\rm a}\), that accounts for the observational 
evidence. This Bayesian step reads
\begin{align}
  \pi^{\rm a}(z\vert y_{\rm obs}) 
  = \frac{\pi(y_{\rm obs}\vert z)}{\pi^f(y_{\rm obs})}\pi^{\rm f}(z),
\end{align} 
where \(z = (q\trp,p\trp)\trp \in \mathbb{R}^{2N}\) denotes model states.

In this work we focus on the second approach, and on how to apply sequential
Bayesian data assimilation to models of the form~\eqref{eq:model1}.
To this end we consider the deterministic evolution under these model
equations given normally distributed initial data
\(
z(0) \sim  \mathcal{N}(z_0,Q)
\)
where \(z_0 \in T\mathcal{M}\). Furthermore we assume linear observations
\begin{align}\label{eq:model_obs}
  y_{\rm obs}(t_k)=H_{\rm{obs}}z(t_k)+\zeta
\end{align}
where \(H_{\rm obs} \in \mathbb{R}^{I \times 2N}\) is the matrix representing
the linear observation map, and $\zeta \sim \mathcal{N}(0,R)$ is the
measurement error with Gaussian statistics. Hereby \(I\in \mathbb{N}\) is the
dimension of the observation space.

For linear models, Gaussian measurement error, and Gaussian initial data, the
\textit{Kalman filter} solves the problem of optimally matching the analysis 
distributions to the observations~\cite{Kalman1960}. Since the Gaussian structure
of probabilities is exactly preserved in this case, the prior and posterior 
densities are completely characterized by their means 
\(\bar{z}^{\rm f}(t_k),\bar{z}^{\rm a}(t_k)\) and covariances 
\(P^f(t_k),P^a(t_k)\) at time \(t_k\).

When the model equations are nonlinear and therefore the forecast
distribution is not Gaussian anymore, we still can recover the main idea of
the Kalman filter and approximate \(\bar{z}^{\rm f}(t_k)\) and \(P^{\rm
f}(t_k)\) by their empirical counterpart and use the \textit{ensemble Kalman
filter} (EnKF)~\cite{Evensen2003} to obtain the posterior mean
\(\bar{z}^a(t_k)\) and covariance \(P^a(t_k)\). To be more specific we take
the existing ensemble 
\begin{align}
Z^a(t_{k-1}) = \Bigl(z_i^{\rm a}(t_{k-1})\Bigr)_{i = 1}^{M} \in \mathbb{R}^{2N\times M}
\end{align} 
which represents the distribution \(\pi^a(z, t_{k-1})\) and evolve each member 
according to the model equations in time until~\(t_k\).
Now the resulting ensemble 
\(Z^f(t_k) = {\bigl(z_i^{\rm f}(t_{k})\bigr)}_{i = 1}^{M}\) 
samples the prior density \(\pi^f(z,t_k)\) and we use \(\bar{z}^{\rm f}(t_k)\approx
\frac{1}{M}\sum_{i=1}^M z_i^{\rm f}(t_k)\eqqcolon \bar{z}^{\rm f}\) and \(P^{\rm
f}(t_k) \approx \frac{1}{M-1} \sum_{i=1}^M (z_i^{\rm f}(t_k) - \bar{z}^{\rm
f}) (z_i^{\rm f}(t_k) - \bar{z}^{\rm f})\trp\) to estimate the first and
second moments of \(\pi^f(z, t_{k})\). To finally transform the prior samples
to samples of the posterior, we assume a linear
transformation
\begin{align}
  z_j^{\rm a}(t_k)=\sum_{i=1}^M z_i^{\rm f}(t_k) \sigma_{ij}(t_k)\qquad j =
  1,\ldots,M \label{eq:ensembleDA},
\end{align}
but are still left with a choice of the transformation matrix coefficients 
\(\sigma_{ij}\)~\cite{reichbook}. Our choice in this study will be the 
ensemble square root filter (ESRF) as described e.g.~in~\cite[p.211-212]{reichbook}.
The corresponding transfer matrix reads
\begin{align}
  \sigma = \text{diag}(w - \frac{1}{M}) + S,
\end{align} 
where 
\begin{subequations}
\begin{align}
   & S = \left( I + \frac{1}{M-1}(H_{obs}A)^{\rm
    T}R^{-1}H_{obs}A\right)^{-\frac{1}{2}} \in \mathbb{R}^{M\times M}, \\
   & w = \frac{1}{M}\sum_{i=1}^M e_i-
  \frac{1}{M-1}S^2A^{\rm T}H_{obs}^{\rm T} R^{-1}(H_{obs}\bar{z}^{\rm f}-y_{\rm
  obs}) \in \mathbb{R}^{M},
\end{align}
\end{subequations}
with $e_i \in \mathbb{R}^M, \ (e_i)_j = \delta_{ij}$, and where 
the ``ensemble anomalies'' are 
\begin{align}
   & A=\left[ z_1^{\rm f}(t_k)-\bar{z}^{\rm f}(t_k) \quad \dots \quad
  z_M^{\rm f}(t_k)-\bar{z}^{\rm f}(t_k) \right] \in \mathbb{R}^{2N\times M}.
\end{align}
Using the ensemble square root filter we avoid the perturbation of the 
observations as necessary for non deterministic versions of the 
EnKF~\cite{Whitaker2002}. Nevertheless our statements do not depend on 
the specific choice made here.

% ==============================================================================

\subsubsection{Failure of the plain ensemble square root filter}

%Our state variable is \(z = (q\trp ,p\trp )\trp \in \mathbb{R}^{2N}\).
Although the Hamiltonian (\ref{eq:Hamiltonian}) is
conserved under the model dynamics~\eqref{eq:model1}, i.e.,
\begin{equation}
  H^\varepsilon \left({z_i^\varepsilon}^{\rm f}(t_{k+1})\right) 
  = H^\varepsilon\left({z_i^\varepsilon}^{\rm a}(t_k)\rule{0pt}{10pt}\right)\,,
\end{equation}
it is not conserved under transformation~\eqref{eq:ensembleDA} which implements
the data assimilation step.
In particular, one often observes a severe increase in the oscillatory energy
(\ref{eq:oscHamiltonian}), i.e.
\begin{equation}
  H_{\rm osc}^\varepsilon\left( {z_i^\varepsilon}^{\rm a} (t_{k})\rule{0pt}{10pt}\right) 
  \gg H_{\rm osc}^\varepsilon\left({z_i^\varepsilon}^{\rm f}(t_k)\right)\,,
\end{equation}
which, in practice, can lead to a destabilization of the simulation after a few
data assimilation cycles. The reader is referred to 
\cite{HarlimMajda2010,KellyEtAl2015} for rigorous analyses of such possible 
catastrophic filter divergences. An explicit example of this effect is also 
provided below, see Fig.~\ref{fig:esrffailure}.
For linear scalar balance relations, the situation can be controlled however, as
summarized in the following 
\begin{rem}\label{lem:invariant}
  Let \(\sigma \in \mathbb{R}^{M\times M}\) be the transformation matrix of a
  linear ensemble transform filter. Let furthermore \(g(q) = G \, q\) with 
  $G \in \mathbb{R}^{L\times N}$ a constant positive semi-definite matrix, 
  so that $g$ is linear. Then for every
  ensemble of prior samples \(q^{\rm f}_j\in \mathbb{R}^{N}\) and posterior
  samples \(q^{\rm a}_j\in \mathbb{R}^{N}\) with \(j\in \{1\dots M\}\)
  \begin{align}
    \label{equ:invariant}
    g\left(q^{\rm a}_j(t_k) \right) \leq  C \max_{i = 1\dots m} g\left( q^{\rm f}_i(t_k)
    \right)
  \end{align} 
  at every time point \(t_k\) with \(C=\sum_{i=1}^m |\sigma_{ij}|\).
  In fact, due to the linearity of \(g\) we can immediately conclude
  \begin{equation}
    g\left(q^{\rm a}_j\right) 
    = g\left(\sum_{i=1}^m q^{\rm f}_i \sigma_{ij} \right) 
    = \sum_{i=1}^m  g\left(q^{\rm f}_i \sigma_{ij} \right)
    = \sum_{i=1}^m  \sigma_{ij} g\left(q^{\rm f}_i \right) 
    \leq \sum_{i=1}^m  |\sigma_{ij}| \max_{l = 1\dots m} g\left(q^{\rm f}_l \right)\, .
  \end{equation}
\end{rem}
\begin{rem}\label{cor:invariant}
  If the ensemble of prior samples is exactly balanced, i.e., satisfies
  \(g(q^{\rm f}_i)=0\) for every \(i=1\dots M\) then the
  ensemble of posterior samples will satisfy
  \(g(q^{\rm f}_i)=0\) for every \(i=1\dots M\), too.
\end{rem}
Note that for nonlinear \(g\) such control is not available, since 
neither~\eqref{equ:invariant} nor Remark~\ref{cor:invariant} remain
valid in general.

Although not of immediate importance for the current assimilation cycle, the
assimilation reduces the mean distance of the ensemble to the observations as
expected, and yet the subsequent forecast can be drastically wrong nevertheless. 
In the case of rather small \(\varepsilon\) this can ultimately lead to filter
divergence. An example of this situation is illustrated in
Figure~\ref{fig:esrffailure}, which also shows results with the improved
balanced DA procedure to be described shortly.

\begin{figure}[H]
  \centering
  \input{\imgPath esrf_failure.pgf}
  \caption{
    The divergence of an ensemble square root filter (ESRF) applied to the
    elastic stiff double pendulum as described in
    Example~\ref{ex:doublependulum}. The initial state was chosen as
    \(q_0={(1,0,2,0)}\trp\), \( p_0 = {(0,0,0,0)}\trp \) and the model
    parameters according to Table~\ref{tab:double_parameters}. The left
    figure depicts the ensemble trajectories (ens.) of the second mass point
    and the associated ensemble means (est.). The plain ESRF diverges after a
    couple of assimilation cycles, whereas the balanced version (bESRF)
    performs qualitatively well. The right figure shows the residual
    \(\|g(q)\|\). The residual of the reference solution is non-zero as
    expected, but stays small. The residual of the balanced ESRF is
    drastically lower than that of the plain ESRF. We obtain these results
    using the method proposed in~\eqref{eq:bmethod} below and the setup as
    described in Chapter~\ref{sec:results}, c.f.~Table~\ref{tab:double_parameters}. 
    In contrast to the situation there, we increase the initial spread and 
    choose \(\rho_0 = 0.1\).
  }\label{fig:esrffailure}
\end{figure}
%

% ==============================================================================
% ==============================================================================
% ==============================================================================

\section{Proposed methods}
\label{sec:ProposedMethods}

To overcome the above-mentioned issue, we propose two different methods. The
first, subsequently called ``penalty method'', observes and corrects the
balance residuals after the assimilation algorithm. For this purpose we solve
a minimization problem structurally similar to the 3DVar method (see
e.g.~\cite{kalnay_atmospheric_2002}). The second, subsequently called
``blended time stepping method'', is an extension of ideas first formulated
in~\cite{BenacchioEtAl2014}. This approach does not modify the assimilated
states but leverages an intermediate model as part of the forecast step that
drives the evolution towards balanced states.

% ==============================================================================
% ==============================================================================

\subsection{Ensemble based penalty method}

Let \(\hat{z}_i = \left({\hat{q}_i}^T, {\hat{p}_i}^T\right)^T\) 
denote the coordinates and momenta of an ensemble \(\hat Z\) provided by 
applying a linear ensemble transform filter to the forecasts 
\(Z^f = \left(z^f_i\right)_{i=1}^{M}\).
The $\hat{z}_i$ would represent the analysis ensemble if we were to ignore 
further balancing requirements. To improve the balance of these states, and 
to obtain the final analysis ensemble \(Z^a\), we first generate an ensemble of 
updated imbalances by applying the ensemble Kalman filter 
transformation~\eqref{eq:ensembleDA} to the forecasted values of \(g\), i.e.,
\begin{equation} \label{transform_g}
  \hat{g}_j \coloneqq \sum_i^M g(q_i^{\rm f}) \,\sigma_{ij}\quad j \in \{1,\dots,
  M\}\,,
\end{equation}
and then minimize \(L: \mathbb{R}^{N\times M} \rightarrow
\mathbb{R}\), which acts on the ensemble of positions $Q = \left(q_i\right)_{i=1}^M$ 
only, and is defined through
\begin{subequations}\label{eq:MinimizationProblem}
\begin{align}
  \label{eq:costfunctional}
  L(Q)     
    & = \frac{1}{2}\sum_{i=1}^M 
      \Bigl((q_i-\hat{q}_i)\trp B (q_i - \hat{q}_i) +  S_i(q_i) \Bigr)
  \\
  S_i(q_i) & = (g(q_i) - \gamma \hat{g}_i)\trp \Lambda (g(q_i) - \gamma \hat{g}_i) 
\end{align}
\end{subequations}
after each assimilation procedure in a post processing step. Note that
ensemble members $i,j$ are not coupled in this minimization problem, so that
it is equivalent to $M$ independent minimizations of smaller size. 

The post-processed balanced posterior ensemble $Z^{\rm a}$ is now given as minimizers, i.e.,
\begin{equation}
  Z^{\rm a} = \left(\arg \min_{Q} L(Q), \hat P \right),
\end{equation}
where $\hat P$ denotes the ensemble of momenta from $\hat Z$.

In \eqref{eq:MinimizationProblem}, \(B\in \mathbb{R}^{N\times N}\) is a symmetric 
positive definite matrix and \(\Lambda\in \mathbb{R}^{L \times L}\) is a positive 
definite diagonal matrix which weighs the importance of the proposed analysis
ensemble against that of improved balance. The matrix $B$ is typically given
by the inverse of the empirical covariance matrix of the ensemble $\hat Q$, i.e.,
the ensemble of coordinates from $\hat Z$. The matrix $\Lambda$ is chosen to be diagonal
with a common scaling factor $\ell>0$ which controls the impact of the balancing
terms $S_i$ in (\ref{eq:costfunctional}). Furthermore, the parameter \(\gamma \in
[0,1]\) controls the quality of balance to be achieved relative to the assimilated
balance residuals~\(\hat{g}_i\). Whereas the $\hat g_i$ may already be reduced in amplitude 
relative to the forecast values $g(q_i^f)$ depending on the structure of the weights
$\sigma_{ij}$, a stronger enforcement of balance in a single assimilation step
may be advantageous. This is achieved by choosing values $\gamma < 1$. The extreme
case of $\gamma = 0$, together with $\|B\| \ll \|\Lambda\|$ would amount to a 
projection of the forecasted states onto the balanced manifold, which - in turn - 
might not be desirable either as this may suppresses physically meaningful smaller
imbalances, too.  At this stage, the concrete choice of $\gamma$ remains subject to
the particular application context. Additional physical considerations will provide further
guidance towards a best practice in the choice of parameters.

%\klein{[Frage: Sollten wir nicht auch noch etwas zur Wahl von $B$ und $\Lambda$ sagen?]}

\begin{rem}
  We deliberately refrain from combining the assimilation algorithm and this
  post processing step into a single assimilation step for the sake of transparency and to advertise the
  flexibility of this post processing approach.
\end{rem}

% ==============================================================================
\subsubsection{Gauß-Newton minimization}

The gradient of~\eqref{eq:costfunctional} evaluated at any minimizer thereof
vanishes i.e.
\begin{equation}
  \label{eq:criticalpoint}
  0 = \frac{\partial L}{\partial q_i} (Q) = B (q_i - \hat{q}_i) 
  + \nabla S_i(q_i) \qquad \forall i\in \{1,\dots,M\}.
\end{equation}
For the approximate numerical solution of this system of coupled nonlinear 
equations we first summarize here the Gauss-Newton algorithm, which we modify 
slightly in the next section to obtain the method actually used in sample 
calculations below. 

Let $Q^k$ denote the $k$th iterate of the ensemble of positions. Then one 
linearizes the cost functional~\eqref{eq:costfunctional} by
\begin{align}
  g(q_i)     & \approx g(q_i^k) + G(q_i^k)(q_i - q_i^k)
\end{align}
and thus obtains
\begin{equation}
  \label{eq:criticalpointlinear}
  \begin{split}
    \left(B + 	G\trp(q_i^k) \Lambda G(q_i^k)\rule{0pt}{10pt}\right)(q_i - q_i^k) 
    = - B (q_i^k- \hat{q}_i) -  \nabla S_i(q^k_i)
  \end{split} \qquad \forall i\in \{1,\dots M\}
\end{equation}
as the linearized critical point condition.  Solving this system of linear 
equations for the ensemble of increments \( \Delta Q^k \coloneqq Q^{k+1} -  Q^k \) 
is equivalent to the inversion of 
\begin{equation}
  \label{eq:metric}
  \tilde{B}_i^k  \coloneqq {\left(B +  G(q_i^k)\trp \Lambda G(q_i^k)
    \right)}. 
\end{equation}
Each of these matrices is indeed invertible, since \(B\) is symmetric positive 
definite and $\Lambda$ is strictly positive diagonal. This allows one to formulate 
the update increment for each ensemble member independently, i.e., 
\begin{equation}
  \label{eq:increment}
    \Delta q_i^k = - \bigl({\tilde{B}}^k_i\bigr)^{-1} \left(B (q_i^k- \hat{q}_i) 
    + \nabla S_i(q^k_i)   \right).
\end{equation}
\begin{rem}\label{rem:SMWformula}
  For \(L\ll N\) one may reduce the computational costs for each iteration
  step by expressing \(\bigl({\tilde{B}}^k_i\bigr)^{-1}=(B +  G(q_i^k)\trp \Lambda
  G(q_i^k))^{-1}\) by the Sherman-Morrison-Woodbury formula~\cite[p.
  51]{golub1996}.
\end{rem}
%
%Subsequently we follow the Quasi-Newton method again and descend along the
%gradient using
%%
%\begin{equation}\label{eq:gaussnewtonupdate}
%  q_i^{k+1} - q_i^{k} = h \Delta q^k_i.
%\end{equation}
%%
%Hereby \(h>0\) is a step size determined by a line search for a local minimum
%of the balance functional \(L\), for details see e.g.~\cite{quarteroni2007}.

% ==============================================================================

\subsubsection{Continuous formulation and a modified search direction}

Instead of using the gradient descent proposed by the Gauss-Newton method,
we aim here to solve~\eqref{eq:criticalpoint} using a different direction of 
descent. To this end we introduce a modified search direction,
obtained by replacing \(\Lambda\) in~\eqref{eq:metric} by \(\Lambda h\) where
\(h>0\). In analogy with~\eqref{eq:increment} we denote the ensemble increment 
obtained by this method by \(\Delta \tilde{Q}^k\), and we note that 
all fixed points of \( Q^{k} \mapsto Q^k + h \Delta\tilde{Q}^k \) are
solutions of~\eqref{eq:criticalpoint} and therefore again candidates for
minimizers.

A single iteration by this method turns out to be one step with pseudo-time 
increment $\Delta s = h$ of a stable numerical integrator for the auxiliary 
differential equation system 
\begin{equation}
    \label{eq:3DVar}
    \frac{d}{ds} q_i = - (q_i-\hat{q}_i)  
    - B^{-1}G\trp(q_i)\Lambda(g(q_i) - \gamma \hat{g}_i) 
    \qquad \forall i\in \{1,\dots, M\}\,.
\end{equation}
This integrator, augmented with an adaptive step size optimized for rapid
reduction of the residuum, and applied to the initial value problem with initial 
data \(q_i(0) = \hat{q}_i\) constitutes the proposed minimizer for the 
balancing cost function in \eqref{eq:MinimizationProblem}. This is also our
preferred method used in the test cases below.
\begin{rem}
  The evolution governed by~\eqref{eq:3DVar} is a gradient flow driven by
  \(L\) and the geometry of \( \mathrm{diag}(B^{-1},\dots, B^{-1}) \in
  \mathbb{R}^{MN \times MN}\). We therefore expect the solution
  of~\eqref{eq:3DVar} to converge to an equilibrium solution \(q_\infty =
  \lim_{t\rightarrow \infty} q(t)\) satisfying~\eqref{eq:criticalpoint}.
\end{rem}
\begin{prop}
  The numerical method governed by~\eqref{eq:increment} with $\Lambda$ replaced
  by $h\Lambda$ for a pseudo-time increment $s^{k+1}-s^k$, i.e.,
  \begin{equation}
    \label{eq:3DVarContDiscrete}
    q_i^{k+1} = q_i^{k} - h {\left(B + h  G\trp(q_i^k)\Lambda G(q_i^k)\right)}^{-1}
    \left(B (q_i^k- \hat{q}_i)  
    +  \nabla S_i(q^k_i)   \right)
  \end{equation}
  is consistent with~\eqref{eq:3DVar}.
  For \(h\in (0,1)\) there exists \(\delta>0\) such that every 
  sequence~\({(q_i^{k})}_{k\in \mathbb{N}}\) determined 
  by~\eqref{eq:3DVarContDiscrete} and starting in the open ball 
  \(B_\delta(q_i^\infty)\) around an equilibrium solution \(q_i^{\infty}\)
  converges to that solution.
  \begin{figure}[!h]
    \begin{tikzcd}[row sep=large,column sep=large]
      q_i^k \arrow[d,"h \rightarrow 0"] \arrow[dr,"k \rightarrow \infty"]    \\
      q_i(t) \arrow[r,"t \rightarrow \infty"]  & q_i^\infty
    \end{tikzcd}
    \centering
    \caption{The commuting diagram shows the stability property of discretization~\eqref{eq:3DVarContDiscrete}.}
  \end{figure}
\end{prop}
\begin{proof}
  We recall that \(B\in \mathbb{R}^{N\times N}\) is invertible and bounded as
  it is finite dimensional. The expansion given by the Neumann series gives
  \begin{equation}
  \begin{aligned}
    {\left( B + h  G\trp(q_i^k)\Lambda G(q_i^k) \right)}^{-1} 
      & = {\left( \id + h B^{-1}  G\trp(q_i^k)\Lambda G(q_i^k) \right)}^{-1} B^{-1} 
        = \sum_{k=0}^\infty {\left(-B^{-1} G\trp(q_i^k) \Lambda G(q_i^k) \right)}^k B^{-1} 
        \\
      & = B^{-1} + \mathcal{O}(h).
  \end{aligned}
  \end{equation}
  This allows us to conclude consistency by a standard Taylor argument, expanding the 
  solution \(q_i\) at \(t^k\) and assuming \(q_i^n=q_i(t^k)\). For this purpose first 
  rewrite 
  \begin{subequations}
  \begin{align}
    q_i^{k+1} 
      & = q_i^{k} - h {\left(B + h  G\trp(q_i^k)\Lambda G(q_i^k)\right)}^{-1} 
          \frac{\partial L}{\partial q_i}(Q^k) 
          \\ 
      &= q_i^{k} - h {B}^{-1} \frac{\partial L}{\partial q_i}(Q^k) + \mathcal{O}(h^2)
  \end{align}
  \end{subequations}
  and subsequently conclude 
  \begin{align}
  \left\| q_i(t^{k+1})-q_i^{k + 1}\right\| 
  = \left\| q_i(t^k) - h {B}^{-1} \frac{\partial L}{\partial q_i}(Q(t^k)) 
    - q_i^k  + h {B}^{-1} \frac{\partial L}{\partial q_i}(Q^k) + \mathcal{O}(h^2) 
    \right\| = \left\|O(h^2) \right\|.
  \end{align}
  This implies global first order consistency and the first part of the 
  statement. For the second part let \( h \in (0,1)\). 
  Subtracting \(q_i^{\infty}\) on both sides and furthermore using 
  \begin{equation}
    0 = B(q_i^\infty - \hat{q}_i) +  G(q^k_i)\Lambda (g(q^k_i)-\bar{\hat{g}})
  \end{equation}
  allows us to conclude equivalence of~\eqref{eq:3DVarContDiscrete} and the
  following identity,
  \begin{equation}
    \begin{aligned}
    q_i^{k+1} - q^{\infty}_i 
      & = q_i^{k} - q^{\infty}_i 
          - h {\left(B + h  G\trp(q_i^k)\Lambda G(q_i^k)\right)}^{-1} 
            \left(B (q_i^k- q_i^{\infty}) +  G\trp(q^k_i)\Lambda(g(q^k_i) - g(q^\infty_i) ) \right) 
        \\
      & = \left( \id - h {\left(B + h  G\trp(q_i^k)\Lambda G(q_i^k)\right)}^{-1} 
          \left(B  +  G\trp(q^k_i)\Lambda G(q^k_i) \right) \right)(q^k_i - q^\infty_i ) 
          + h \, r(q_i^k,q_i^k-q_i^\infty)\,.
    \end{aligned}
  \end{equation}
  The last equality is valid as long as \(q^{k}\in B_\rho(q^\infty)\) for sufficiently small \(\rho>0\). 
  In this case we can apply Taylor expansion which also gives us 
  \begin{equation}
   R(w_i) \coloneqq \sup_{v \in B_\rho(q^\infty)} \|r(v_i,w_i)\| \in \mathcal{O}(\|w_i\|^2).
  \end{equation}   
  The fact that \(B\) and \(G(q_i^k)\trp \Lambda G(q_i^k)\) are both
  symmetric positive definite allows us to conclude the following estimate
  \begin{subequations}
  \begin{align}
    \|q_i^{k+1} - q^{\infty}_i\| 
      & \leq  \|  \id - h {\left(B + h  G\trp(q_i^k)\Lambda G(q_i^k)\right)}^{-1} 
        \left(B +  G\trp(q^k_i)\Lambda G(q^k_i)\right)\| \| q_i^{k} - q^{\infty}_i\| 
        + h R(q_i^k-q_i^\infty)
        \\ 
      & \leq  (1-h) \| {\left(B + h  G\trp(q_i^k)\Lambda G(q_i^k)\right)}^{-1} B  
        \| \| q_i^{k} - q^{\infty}_i\| + h R(q_i^k-q_i^\infty)
        \\
      & \leq  \frac{(1-h) \|B\|}{\| {\left(B + h  G\trp(q_i^k)\Lambda G(q_i^k)\right)}\|}   
        \| q_i^{k} - q^{\infty}_i\| + h R(q_i^k-q_i^\infty)
        \\ 
      & \leq  (1-h) \| q_i^{k} - q^{\infty}_i\| + h C_1 \| q_i^{k} - q^{\infty}_i\|^2
        \\
      & \leq  (1 - h  + h C_1 \| q_i^{k} - q^{\infty}_i\|) \| q_i^{k} - q^{\infty}_i\| 
        \\
      & \leq  C_2 \| q_i^{k} - q^{\infty}_i\|.
  \end{align}
  \end{subequations}
  Here the constant satisfies \(C_1<1\) as long as \(q_i^k\) is already
  close enough to \(q_i^{\infty}\). If so, then we immediately obtain \(
  \left\| q_i^{k+1}-q_i^{\infty}\right\| C_1 <1\) and therefore \(C_2<1\)
  too. An inductive argument finally implies convergence to the equilibrium
  solution \(q_i^{\infty}\) for sufficiently close initial value.
\end{proof}

\subsection{Blended time-stepping}
Motivated by the results of~\cite{BenacchioEtAl2014}, we introduce a numerical
time stepping scheme that extends a classical projection approach by subsequent
blending steps. These steps begin with a classical projection step and then 
access successively less constrained intermediate models along a one-parameter
model family that continuously bridges between the unconstrained original
model~\eqref{eq:model1} and the reduced and its fully constrained 
limit~\eqref{eq:limiteq}. This approach was originally developed in the context 
of incompressible fluid dynamics where the singular perturbation arises by the 
vanishing Mach number limit \(\textrm{Ma}\rightarrow 0\). In addition to the 
classical projection schemes introduced in~\cite{Chorin1967}, much effort was 
spent on developing asymptotic preserving low Mach number numerical schemes 
and variants thereof, \citep[see, e.g.,][]{KleinEtAl2001,CordierEtAl2012,Jin2012}.
The essential point is their ability to blend between the (weakly)
compressible and the incompressible dynamics without additional stability
constraints. As observed in~\cite{BenacchioEtAl2014}, solving the
incompressible model immediately after the assimilation for one or two time
steps and subsequently blending smoothly back to the compressible model over 
another few time steps can further reduce artificial imbalances caused by data 
assimilation relative to an approach that simply projects the system state 
onto the incompressible manifold in one step and then proceeds with the 
compressible model.

To adapt this strategy to our situation, we introduce the following family
of blended models, controlled by the blending parameter $\alpha$, 
\begin{subequations}\label{eq:blendedmodel}
\begin{equation}
  \begin{alignedat}{3}
    \dot{q} 
      & = p\, ,
        & \qquad q(0) 
          & = q_0 \in \mathbb{R}^N
            \\
    \dot{p} 
      & = - G\trp(q) K 
            \Bigl(\frac{\alpha}{\varepsilon^2} \, g(q) + (1 - \alpha) \, \lambda^0 \Bigr)- \nabla V(q)\,,
        & \qquad p(0) 
          & = p_0 \in \mathbb{R}^N
%           \\
%    G(q) p 
%      & = 0\,.
  \end{alignedat}
\end{equation} 
where $\lambda^0$ is the lagrangian multiplier as calculated from the limit problem 
\begin{equation}\label{eq:limitmodel2}
  \begin{alignedat}{3}
    \dot{q} 
      & = p\, ,
        & \qquad q(0) 
          & = q_0 \in \mathbb{R}^N
            \\
    \dot{p} 
      & = - G\trp(q) \, K\, \lambda^0 - \nabla V(q)\,,
        & \qquad p(0) 
          & = p_0 \in \mathbb{R}^N
            \\
    G(q) \, p 
      & = 0\,.
  \end{alignedat}
\end{equation} 
\end{subequations}
Obviously, for $\alpha = 1$ we recover the unconstrained dynamics from 
\eqref{eq:model1}, whereas for $\alpha = 0$ and provided $(q_0, p_0) \in T{\cal M}$, 
we recover the constrained dynamics from \eqref{eq:limiteq}. Note, however, 
that in this blended model we consciously use the hidden limit constraint for the 
momenta, i.e., $G(q)\, p = 0$, rather than the original constraint $g(q) = 0$. 
As we will demonstrate below, this ensures a desired dissipative behavior towards 
balanced solutions when $\varepsilon$ is small but non-zero and 
$0 < \alpha < 1$. For the limiting cases $\alpha \in \{0,1\}$, however, the 
non-dissipative symplectic integrator is maintained.

Let us denote by 
\begin{align}\label{eq:blending}
  z^{n+1}=\psi^\alpha_h(z^n) \qquad \alpha
  \in [0,1]
\end{align}
a numerical discretization of \eqref{eq:blendedmodel}, so that the operator
$\psi^\alpha_h$ advances a given solution $z^n$ at time $t^n$ by a time step
$h = t^{n+1}-t^n$ to the next time level.
For \(\alpha = 0\) we employ a projection method that keeps the momenta tangential 
to the manifold $g(q) = g(q_0)$, whereas for \(\alpha = 1\) we use the nearly 
energy-preserving Stömer-Verlet scheme to resolve the unconstrained 
model~\eqref{eq:model1}. The time discretization for \(\alpha\in (0,1) \) will 
have dissipative character as shown in sections~\ref{sssec:DissipativeModel}, \ref{sssec:DissipativeDiscretization} below, and it is designed to efficiently 
remove any artificially introduced oscillatory energy from the system.

\begin{figure}[h!]
  \centering
  \begin{subfigure}[b]{0.49\textwidth}
    \includegraphics[width=\textwidth]{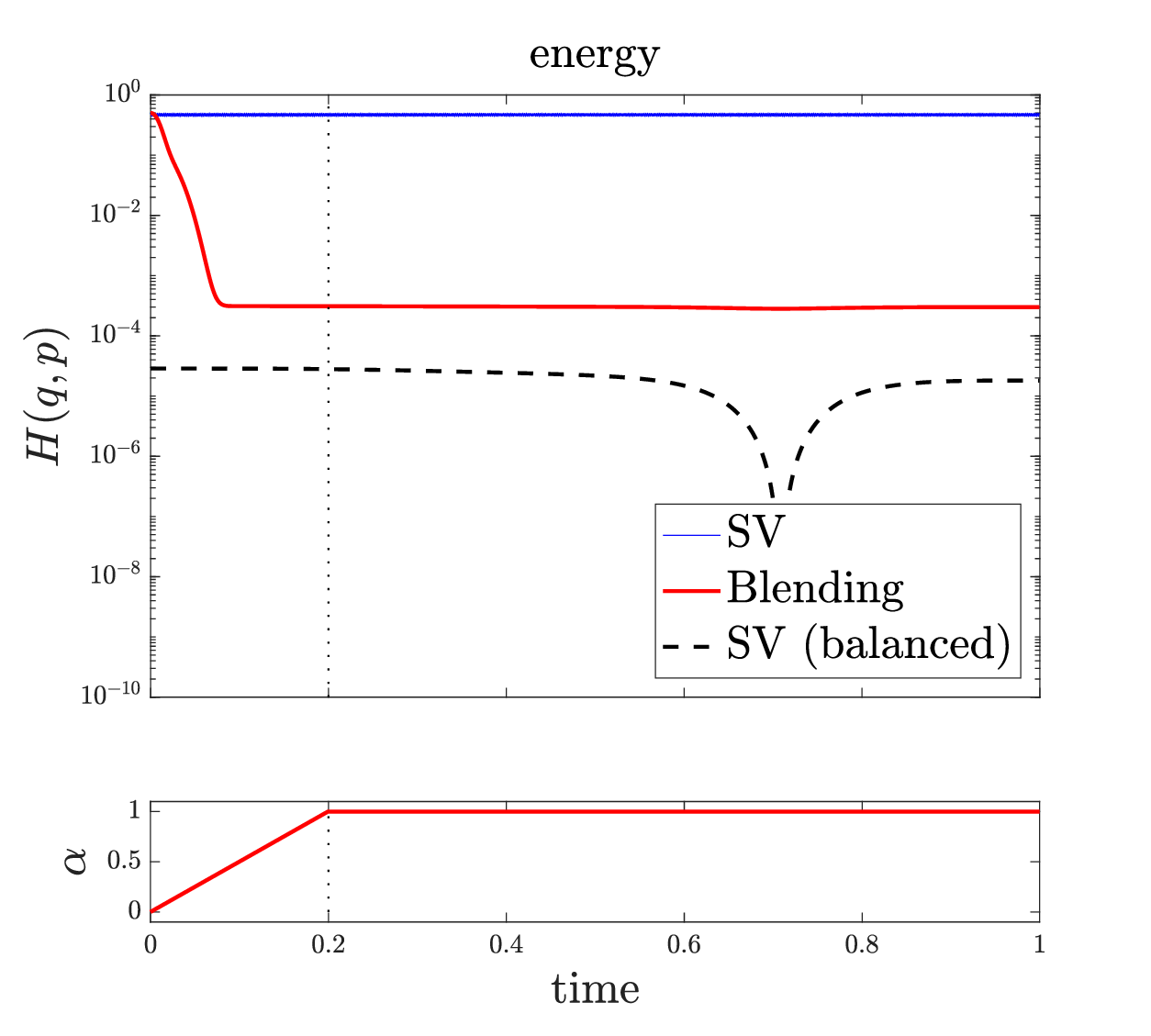}
  \end{subfigure}
  \begin{subfigure}[b]{0.49\textwidth}
    \includegraphics[width=\textwidth]{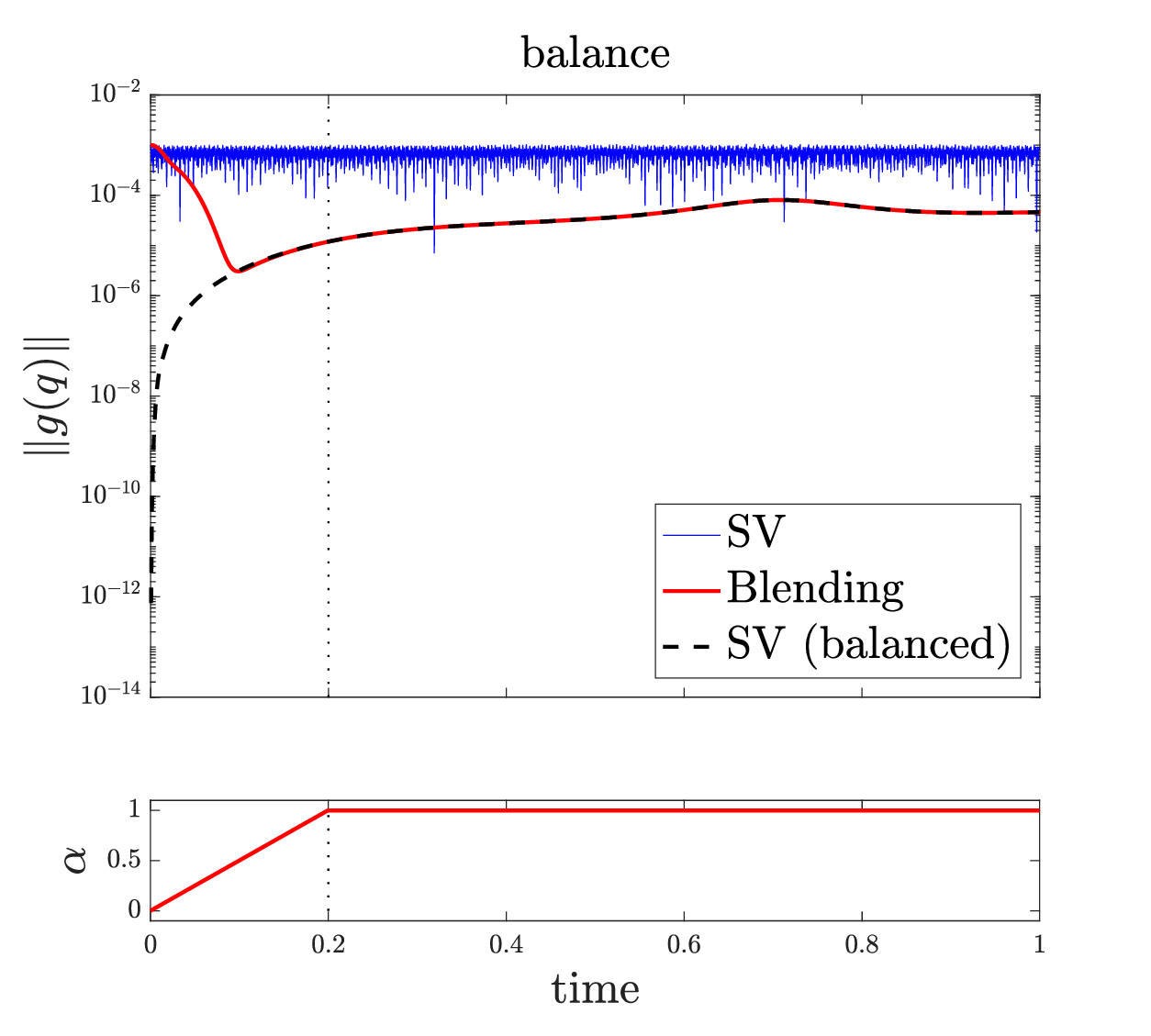}
  \end{subfigure}
  \caption{
    Energy and residual of balance relation for the stiff elastic
    double pendulum (c.f.~Example~\ref{ex:doublependulum} and 
    Fig.~\ref{fig:esrffailure}).  Initially
    unbalanced numerical solutions based on the St\"ormer-Verlet~(SV) (blue)
    and the blending method (red), respectively. The corresponding parameter
    \(\alpha \) for the blended time stepping is shown and the start of the
    pure forecasting region is marked by a vertical dotted line. As the
    reference we display the data for an initially balanced solution,
    computed again by the St\"ormer-Verlet method (dashed).
  }\label{fig:numericalmethod}
\end{figure}

In this approach we accept a non-vanishing consistency
error with respect to the fast model when evolving the system with \(\alpha
\in [0,1) \). But, as discussed in the beginning, we can assume the
solutions of the unconstrained system to stay \(\varepsilon \)-close to
the solutions of the constrained one. This enables us to locally decompose
the consistency error into two parts, one in \(\mathcal{M} \), caused by the
nonlinearity of \(V\) and another one orthogonal to \( \mathcal{M}\). The slow
first part is assumed to be captured by the data assimilation, whereas the
second fast part is small of order \(\mathcal{O}(\varepsilon)\) as discussed 
before. Yet, when artificial imbalances are introduced through a data assimilation 
procedure, the latter assumption ceases to be valid. In this case, the discrete evolution 
of the blended method should rapidly dampen the fast oscillations orthogonal to 
\(\mathcal{M}\) as long as \(\alpha \in (0,1)\) and until they attain the 
correct magnitude of \(\mathcal{O}(\varepsilon)\). To this end we propose to 
use the blending method~\eqref{eq:blending} as follows. Let us denote the 
\textit{blending window} by \(k\geq 1\) and start our forecast at time \(t_n\). 
Let \(\eta \) be the number of forecast time integration steps. Then the 
following two steps are repeated in every forecast cycle 
(c.f.~Figure~\ref{fig:blending}).
\begin{enumerate}
  \item \textbf{Blending:}
        Let \( \alpha \in \mathbb{R}^k\) such that \(0=\alpha_1 < \alpha_2 \leq
        \dots \leq \alpha_{k-1} <
        \alpha_k=1\). Integrate until \(t_{n+k}\) using
        \begin{equation} q_{n+k} = \left(\psi_h^{\alpha_k} \circ
          \psi_h^{\alpha_{k-1}} \circ \dots	 \circ
          \psi_h^{\alpha_{1}}\right)\left(q_n\right)
        \end{equation}
  \item  \textbf{Forecast:}
        Obtain forecast at \(t_{n+\eta}\) by evolving \(q_{n+k}\) along
        \(\psi^1_h\) for \(\eta-k\)~--~time steps.
\end{enumerate}
\begin{figure}[h]
  \centering
  \input{\imgPath blending.tex}
  \caption{
    Blended time stepping applied after analysis with blending window \(k\).
    The numerical flow \(\psi^\alpha_h\) with step width \(h\) is given
    by~\eqref{eq:blending}.
  }\label{fig:blending}
\end{figure}
Figure~\ref{fig:numericalmethod} illustrates the qualitative behaviour of the
blended time stepping for the stiff elastic double pendulum, introduced in
Example~\ref{ex:doublependulum}, with slightly unbalanced initial
coordinates. We observe that with respect to the balance residual the blended
time stepping improves the situation drastically. After short time the
residuals of the initially unbalanced and initially balanced solution match.
Since we dissipate energy in the fast variables (c.f.
Lemma~\ref{lem:dissipative}) as long as \(\alpha\notin \{0,1\} \), the
overall energy of the system decreases as expected and the slow rotational
motion of the stiff double pendulum is therefore resolved reasonably well
with regard to balance and energy. Nevertheless, due to the lack of a priori
knowledge so far, we chose \(\alpha \in [0,1] \) and as a linear function of 
time, but as we
already can guess from the form of the decay, this is a brute force and
suboptimal choice in the sense that we can find a smaller range of \(\alpha\) 
within which the solution relaxes to the slow motion more quickly. We
leave the development of an optimized control of the blending sequence for
future work.
%

% ==============================================================================

\subsubsection{A dissipative model family}
\label{sssec:DissipativeModel}

The aim of the present and the next sections is to analyse the 
proposed balanced data assimilation strategy, and in particular its behavior 
within the blending window while $\alpha \in (0,1)$, in the spirit of
a ``modified equation'' analysis~\citep{WarmingHyett1974}. The modified
equation system that effectively represents our discretization of the
blended model reads
\begin{equation}\label{eq:model3}
  \begin{aligned}
    \dot{\tilde{q}} 
      & = \tilde{p}
        \\
    \dot{\tilde{p}}  
      & = - \frac{1}{\varepsilon^{2}} G(\tilde{q})\trp K g(\tilde{q}) 
        - d\, \mathcal{P}_{\tilde{q}} \,\tilde{p} 
        - \nabla V(\tilde{q}),
  \end{aligned}
\end{equation}
where \(d\) is a scalar damping coefficient. A key characteristic of this 
model is that, owing to the projection $\mathcal{P}_{\tilde{q}}$ in the 
damping term (see \eqref{eq:ProjectionOperator}), the dissipation term 
acts predominantly on the fast oscillatory components of the solution, 
thereby pulling initially unbalanced states towards balanced conditions 
as time evolves. In 
section~\ref{sssec:DissipativeDiscretization} we will show that a particular
time discretization of the blended model in \eqref{eq:blendedmodel} is 
consistent with solutions of this dissipative surrogate model when the
relationship $\alpha = (1-h\,d)$ between the blending parameter $\alpha$ from 
\eqref{eq:blendedmodel}, the damping parameter $d$ from \eqref{eq:model3}, and 
the time step size $h$ is observed. This analogy, in the sense of a modified 
equation analysis of the discretization, will largely explain the success of 
the proposed blending strategy in controlling unwanted imbalances in a data 
assimilation procedure. See also Fig.~\ref{fig:commutingdiagram} below. 
Note that a similar damping term was first proposed in the context of numerical 
stabilization techniques in~\cite{Gear1986}. 

Since fast and slow energy 
parts of the Hamiltonian~\eqref{eq:Hamiltonian} can be separated only by an
asymptotic argument and are coupled nonlinearly, we do not expect the
surrogate system in \eqref{eq:model3} to completely dissipate the 
fast energy of~\eqref{eq:Hamiltonian}. Nevertheless, in reasonably well 
separated cases the impacts of finite $\varepsilon$ and nonlinearity will 
be negligible, especially in the context of data assimilation,
where the correct slow energy itself is anyway known only with limited
accuracy. 

To develop an intuition regarding the behavior of solutions
to~\eqref{eq:model3} we discuss the arguably simplest model in the class of
such problems with multiple scales, the uncoupled harmonic oscillator.
\begin{exm}[Damped harmonic oscillator]
  Let \(q=(\mu,\nu)\trp\), \(p=(\eta,\zeta )\trp\),
  \(K=\mathrm{diag}(\varepsilon^2,1)\) and \(V(q) =0 \)
  \begin{subequations}
  \begin{alignat}{3}
    \label{eq:harmonicslow}
    \dot{\mu} & = \eta	 \qquad   &             & \dot{\eta} &
              & = - \mu                                          \\
    \label{eq:harmonicfast}
    \dot{\nu} & = \zeta  \qquad
              &                 & \dot{\zeta} &            & = -
    \frac{1}{\varepsilon^2} \nu - 4 d \zeta
  \end{alignat}
  \end{subequations}
  The well known analytical solutions for the damped harmonic
  oscillator are given by \( \mu = \mu_0 \cos(t) + \nu_0 \sin(t) \) and
  \begin{equation}
    \nu = \begin{cases}
      e^{- 2 d t} \left( \nu_0 + t (2 d \nu_0 + \zeta_0)  \right) 
      & 2 d = \frac{1}{\varepsilon} \\
      e^{-2 d t} \left( (\nu_0+\frac{\zeta_0}{\omega_d}) e^{i \omega_d t} +
      (\nu_0-\frac{\zeta_0}{\omega_d}) e^{- i \omega_d t}  \right) 
      & 2 d \neq
      \frac{1}{\varepsilon}
    \end{cases}
  \end{equation}
  where the frequency for the fast damped component is given by \(\omega_d =
  \sqrt{\frac{1}{\varepsilon^2} - 4 d^2}\).
  We immediately realize that \(d =0\) gives us the solution for the highly
  oscillatory system~\eqref{eq:model1} and furthermore relaxes to the
  constraint (in this case also slow) manifold exponentially. In the
  general nonlinear and coupled case we present the corresponding result in
  Lemma~\ref{lem:dissipative} below.

  For the overdamped limit, i.e., for \(d \varepsilon \rightarrow \infty \) as
  \(\varepsilon \rightarrow 0\) and \(d \rightarrow \infty \), we conclude
  uniform convergence to the same solution \((\mu,0)\) as for the constrained
  system as long as \(\nu_0 \in o(1/d)\) and \(\zeta_0 \in o(1)\).
  Again this result can be stated in more general form and is presented in
  Lemma~\ref{lem:singpert}.
\end{exm}

The following lemma summarizes the well known (c.f.~\cite{Reich1999})
split of variables into a slow tangential and a fast normal part. It
will enable us to identify slow and fast variables with respect to the
different asymptotic limits.
\begin{lem}\label{lem:split}
  Any solution \((q, p)\) of System~\eqref{eq:model3} can be split into components 
  \(\mu,\eta \) and \(\nu,\zeta \) which satisfy
  \begin{equation}
    \label{eq:split}
    \begin{aligned}
      q & = E\trp \mu  + G\trp{(GG\trp)}^{-\frac{1}{2}} \nu     \\
      p & = E\trp \eta + G\trp{(GG\trp)}^{-\frac{1}{2}} \zeta\,.
    \end{aligned}
  \end{equation}
  In the new coordinates ((\(\mu \), \(\eta \)), (\(\nu \), \(\zeta \))),
  system~\eqref{eq:model3} is equivalent to
  \begin{subequations}
  \begin{alignat}{3}
    \label{eq:slow}
    \dot{\mu} & = \dot{E} q  + \eta \qquad
              & \dot{\eta}
              & = \dot{E} p - E\nabla V(q)                                        \\
    \label{eq:fast}
    \dot{\nu} & = \left(\frac{d}{dt} {(GG\trp)}^{-\frac{1}{2}}G \right) p + \zeta
    \qquad    & \dot{\zeta}
              & =\left(\frac{d}{dt} {(GG\trp)}^{-\frac{1}{2}}G \right) p -
    \varepsilon^{-2} {(GG\trp)}^{\frac{1}{2}} K g(q) - d \zeta -
    {(GG\trp)}^{-\frac{1}{2}} G \nabla V(q).
  \end{alignat}
  \end{subequations}
  Furthermore the fast energy \(H_{osc}\) satisfies
  \begin{equation}
    \label{eq:oscSlowFast}
    H^\varepsilon_{osc}((\nu,\zeta),(\mu,\eta)) = \frac{1}{2}\zeta\trp\zeta +
    \frac{1}{2\varepsilon^2}g(q(\mu,\nu))\trp Kg(q(\mu,\nu))
  \end{equation}
\end{lem}
\begin{proof}
  For the sake of readability we omitted the argument \(q\) for \(G,E\)
  in the statement of the Lemma and will do so throughout the proof. 
  We split momenta tangential and orthogonal to
  \(\mathcal{M}\) denoted by \(\eta \) and \(\zeta \) as well as the
  coordinates denoted by \(\mu \) and \(\nu \), respectively. More concretely
  we choose
  \begin{equation}
    \label{eq:choice}
    \begin{alignedat}{3}
      \mu & = E q \qquad & \eta  & = E p \\
      \nu & = {(GG\trp)}^{-\frac{1}{2}}G q \qquad & \zeta & =
      {(GG\trp)}^{-\frac{1}{2}}G p
    \end{alignedat}
  \end{equation}
  where the columns of \(E \in \mathbb{R}^{(N-L) \times N}\) are an
  orthonormal basis of \( {\left(\mathcal{P}_q \mathbb{R}^N\right)}^\perp \).
  It is easy to check that \(E\trp E\) is a orthogonal projection onto
  \({\left(\mathcal{P}_q \mathbb{R}^N\right)}^\perp \) and since orthogonal
  projections onto a fixed subspace are unique (for every \(q\)), we already
  know \(E\trp E = \mathcal{P}^\perp_q = \id - \mathcal{P}_q\).
  Substituting~\eqref{eq:choice} into the right hand side of~\eqref{eq:split}
  we get
  \begin{equation}
    E\trp E q + G\trp(GG\trp)^{-1}G q 
    = \mathcal{P}^\perp_q q + \mathcal{P}_q q 
    = q.
  \end{equation} 
  Since we used the same geometry to split the momenta this already
  implies~\eqref{eq:split}.

  We differentiate~\eqref{eq:choice} and use system~\eqref{eq:model3} for
  \(\dot{q}\) and \(\dot{p}\). Since by construction \(G\trp E=0=E\trp G\)
  most of the terms drop and we conclude~\eqref{eq:slow} and~\eqref{eq:fast}
  after some straightforward algebraic manipulation. Recalling \(GG\trp \) is
  symmetric positive definite, the last statement~\eqref{eq:oscSlowFast}
  finally follows from
  \begin{subequations}
    \begin{align}
      H^\varepsilon_{osc}(q,p) & = \frac{1}{2} p\trp \mathcal{P}_q p +
      \frac{1}{2\varepsilon^2} g(q)\trp K g(q)
      \\
                               & = \frac{1}{2} p\trp
      G{\trp{(GG\trp)}^{-\frac{1}{2}}}\trp {(GG\trp)}^{-\frac{1}{2}}G p +
      \frac{1}{2\varepsilon^2} g(q)\trp K g(q)                         \\
                               & = \frac{1}{2} \zeta\trp \zeta +
      \frac{1}{2\varepsilon^2} g(q)\trp K g(q).
    \end{align}
  \end{subequations}
\end{proof}
\begin{lem}\label{lem:dissipative}
  Let \(d = \frac{c}{\varepsilon}\), with \(c>0\) fixed. Then solutions
  to~\eqref{eq:model3} which initially satisfy \(H^\varepsilon_{osc}(q,p) \in
  \mathcal{O}(\varepsilon^{-2})\), dissipate fast energy down to some
  residual of order \(H^\varepsilon_{osc} \in \mathcal{O}(1)\), if only
  \(\varepsilon\) is sufficiently small.
\end{lem}
\begin{proof}
  Again we will omit the arguments of \(E\) and \(G\) for notational
  convenience. Additionally we introduce \(\Gamma \coloneqq
  (GG\trp)^{-\frac{1}{2}}G\).
  We will prove the statement by arguments from geometric singular perturbation
  theory~\cite{Fenichel1979}.
  For this purpose we split system~\eqref{eq:model3} into slow and fast
  parts by the means of Lemma~\ref{lem:split}.
  Subsequently we multiply by \(\varepsilon \) and rescale \(\hat{\zeta} =
  \varepsilon \zeta \) which results in
  \begin{subequations}
  \begin{alignat}{3}
    \dot{\mu}                     & = \dot{E} q  + \eta \qquad
                                  &
    \dot{\eta}                    & = \dot{E} p - E\nabla V(q)                 \\
    \varepsilon \dot{\nu}         & = \varepsilon \dot{\Gamma} q + \hat{\zeta}
    \qquad                        &
    \varepsilon \dot{\hat{\zeta}} & = \varepsilon^2 \dot{\Gamma} p -
    {(GG\trp)}^{\frac{1}{2}} K g(q) - \varepsilon d \hat{\zeta} - \varepsilon^2
    {(GG\trp)}^{-\frac{1}{2}} G \nabla V(q).
  \end{alignat}
  \end{subequations}
  We denote the right hand side of the fast variables by
  \begin{align}
    F((\nu,\hat{\zeta}),(\mu,\eta),\varepsilon(d)) =
    \begin{pmatrix}
      \varepsilon \dot{\Gamma} + \hat{\zeta} \\
      - \varepsilon^2 \dot{ \Gamma} p - {(GG\trp)}^{\frac{1}{2}} K g(q) -
      \varepsilon d \hat{\zeta} - \varepsilon^2 {(GG\trp)}^{-\frac{1}{2}} G
      \nabla
      V(q)
    \end{pmatrix}
  \end{align}
  In the limit \(d \varepsilon = {\rm const.} \), \(\varepsilon \rightarrow
  0\) we identify the critical manifold as
  \begin{equation}
    \hat{\mathcal{TM}} \coloneqq \left\{(\eta ,\hat{\zeta},\nu, \mu)\in
    \mathbb{R}^{2N}: g(q(\mu,\nu))=0 \land \hat{\zeta} =0 \right\}  .
  \end{equation}
  Next we prove normal hyperbolicity of the critical manifold, i.e., we show that
  there are no eigenvalues of \(\frac{\partial F}{\partial
    (\nu,\hat{\zeta})}\) with vanishing real part.
  The gradient evaluated on the manifold \(\hat{\mathcal{TM}}\) and for
  \(\varepsilon=0\) is given by the block matrix
  \begin{equation}
    DF\coloneqq\frac{\partial}{\partial (\nu,\hat{\zeta})}
    F((\nu,\hat{\zeta}),(\mu,\eta),\varepsilon)\vert_{(\nu,\hat{\zeta}) \in
      \hat{\mathcal{TM}},\varepsilon = 0}=
    \begin{pmatrix}
      0                                                 & \id  \\
      - {(GG\trp)}^{\frac{1}{2}}K{(GG\trp)}^\frac{1}{2} & -\id
    \end{pmatrix}.
  \end{equation}
  To compute the eigenvalues of this non symmetric matrix, we first
  recall that \({(GG\trp)}^\frac{1}{2}\) is symmetric positive definite and
  since \(K\) is a strictly positive diagonal matrix,
  \(\tilde{K}\coloneqq{(GG\trp)}^\frac{1}{2}K
  {(GG\trp)}^\frac{1}{2}\) is symmetric and positive definite, i.e., it
  has \(L\) positive eigenvalues \(\omega_{\tilde{K},j}>0\).
  Therefore we conclude zero is no eigenvalue of \(DF\) by \(\det DF =
  \det(-\id)\det(-\tilde{K}) \).
  Using the Schur complement again we argue for some eigenvalue \(\omega \neq 0
  \) of \(DF\)
  \begin{equation}
    \det (DF - \omega \id) = \det (-\omega \id) \det(-(1+\omega)\id -
    \frac{1}{\omega}\tilde{K}).
  \end{equation}
  The determinant vanishes if and only if there is \(j \in \{1,\dots, L\} \)
  such that
  \begin{equation}
    -\omega (\omega + 1) = \omega_{\tilde{K},j}.
  \end{equation}
  Solving this quadratic equation already provides us with all possible
  eigenvalues by
  \begin{align}
    \omega_{\pm,j}=\frac{-1\pm \sqrt{1-4 w_{\tilde{K},j}}}{2}.
  \end{align}
  We directly observe \(\mathrm{Re}\, \omega_{\pm,j}<0\) for all \(j \in
  \{1,\dots L\}\) and therefore notice that \(\hat{\mathcal{TM}}\) is
  normally hyperbolic. By finally applying Fenichel's theorem we obtain
  existence of slow manifolds \(S^\varepsilon\) (c.f.~\cite{Kuehn2015})
  \(\varepsilon\)-close to a compact submanifold \(S\subset
  \hat{\mathcal{TM}}\) of the critical one as long as \(\varepsilon\) is
  sufficiently small. More specifically we conclude for any
  \(((\eta^\varepsilon, \zeta^\varepsilon),
  (\mu^\varepsilon,\nu^\varepsilon)) \in S^\varepsilon \)
  \begin{subequations}
  \begin{align}
    {g(q(\mu^\varepsilon,\nu^\varepsilon))}\trp K
    g(q(\mu^\varepsilon,\nu^\varepsilon))                          & \leq c_1
    \varepsilon^2
    \\
    {\hat{\zeta}^{\varepsilon \mathrm{T}} }\hat{\zeta}^\varepsilon & \leq  c_2
    \varepsilon^2
  \end{align}
  \end{subequations}
  and therefore
  \begin{equation}
    \max_{ (\eta^\varepsilon, \zeta^\varepsilon),
      (\mu^\varepsilon,\nu^\varepsilon)) \in S }
    H^\varepsilon_{osc}((\eta^\varepsilon, \zeta^\varepsilon),
    (\mu^\varepsilon,\nu^\varepsilon)) \in \mathcal{O}(1).
  \end{equation}
  Another consequence of Fenichel's theorem is that the dynamical behaviour
  of the linearization \(DF\) of the fast subsystem on the critical manifold
  already determines the dynamical behaviour of solutions starting off a slow
  manifold \(S^\varepsilon\). Since all eigenvalues of \(DF\) have negative
  real part, we conclude \(S\) as well \(S^\varepsilon\) is attracting.
  Therefore any solution starting nearby will approach some \(S^\varepsilon\)
  which finally implies the energy dissipation as stated.
\end{proof}
Subsequently we will use~\eqref{eq:model3} to establish a model hierarchy
which resembles the analytical counterparts discretized by the blended
numerical method~\eqref{eq:blending}. The following two lemmata concern the
behaviour of the limit cases \(d\rightarrow 0 \) and \(d\rightarrow \infty
\). The first one is based on the classical result of continuous dependency
on initial data and parameters for ordinary differential equations with
continuously differentiable right hand side. In both cases we fix
\(\varepsilon>0\) and omit this standard proof.
\begin{lem}\label{lem:contdep}
  Let \(\varepsilon>0\) be fixed. Solutions \((\tilde{q},\tilde{p})\) of the
  dissipative system~\eqref{eq:model3} approach solutions of the purely
  Hamiltonian system~\eqref{eq:model1} as \(d\rightarrow 0\)\,.
\end{lem}
For the other part \(d \rightarrow \infty \) we use again geometric
singular perturbation theory and we can conclude a slightly different type of 
statement in terms of invariant manifolds.
\begin{lem}\label{lem:singpert}
%  For sufficiently large \(d\) and \(\varepsilon^{2} \in o(1/d)\) there
%  exists a Manifold \(\mathcal{M}_{1/d}\) which lies within
%  \(\mathcal{O}(1/d)\) of any compact subset of
%  \(\mathcal{M}_\infty\coloneqq\{(q,p) \in \mathbb{R}^{2N}:\zeta(q,p) = 0\}\)
%  (c.f. Lemma~\ref{lem:split}) and is invariant under the evolution
%  of~\eqref{eq:model3}. Furthermore every solution starting off, but
%  sufficiently close to this subset will approach \(\mathcal{M}_{1/d}\).
Let \(d\) be sufficiently large and \(\varepsilon^{2} \in o(1/d)\). For every compact subset of \(\mathcal{M}_\infty\coloneqq\{(q,p) \in \mathbb{R}^{2N}:\zeta(q,p) = 0\}\)
  (c.f. Lemma~\ref{lem:split}) there exists a manifold \(\mathcal{M}_{1/d}\) which lies within \(\mathcal{O}(1/d)\) of this subset
  and is locally invariant under the evolution
  of~\eqref{eq:model3}. Furthermore every solution starting
  sufficiently close to \(\mathcal{M}_{1/d}\) will approach~\(\mathcal{M}_{1/d}\).
\end{lem}
\begin{proof}
  As pointed out we aim to apply geometric singular perturbation theory again.
  Therefore we start as before by splitting slow and fast momenta explicitly
  utilizing Lemma~\ref{lem:split}.
  Contrary to the situation in Lemma~\ref{lem:dissipative} the coordinates then
  are both slow variables.
  By dividing the momentum equation in~\eqref{eq:fast} by \(d\) and passing to the limit
  \(d\rightarrow \infty\) we obtain the critical manifold as
  \(\mathcal{M_\infty}=\{(q,p) \in \mathbb{R}^{2N}:\zeta(q,p) = 0\}\).
  We denote the right hand side of the momentum equation in~\eqref{eq:fast} by
  \(F((\eta,\zeta),q,1/d)\) and linearize on \(\mathcal{M}_\infty\).
  \begin{align}
    DF\coloneqq\frac{\partial}{\partial \zeta}
    F((\eta,\zeta),q,1/d)\vert_{(\nu,\hat{\zeta}) \in
      \mathcal{M_\infty},\varepsilon = 0}= -\left(GG\trp\right)
  \end{align}
  Since \(GG\trp\) is positive definite and by assumption \(\rank (GG\trp) =
  L\) we conclude that \(-\left(GG\trp\right)^{1/2}\) has exactly \(L\)
  negative Eigenvalues.
  \(M_\infty\) is therefore normally hyperbolic and we now infer by Fenichel's
  theorem~\cite{Kuehn2015} existence of an invariant (with respect
  to~\eqref{eq:model3}) manifold \(\mathcal{M}_{1/d}\),
  \(1/d\)~--~close to a compact subset of our choice of
  \(\mathcal{M}_\infty\),
  exactly as stated.
  Since we additionally have only a stable subspace on \(\mathcal{M}_\infty \)
  we gain the attractive behavior of \(\mathcal{M}_{1/d}\) by the same theorem.
\end{proof}
\begin{rem}\label{rem:modelEquiv}
  The same statement is true if we take a compact submanifold of
  \(\mathcal{M}_\infty\) and therefore also cover the case where the
  evolution starts on the constraint manifold \(\mathcal{M} \subseteq
  \mathcal{M}_\infty\).
\end{rem}
%

%\begin{cor}
%  Let \(d\) be sufficiently large and \(\varepsilon^{2} \in o(1/d)\), then  
%  every compact subset of \(\mathcal{M}_\infty\coloneqq\{(q,p) \in \mathbb{R}^{2N}:\zeta(q,p) = 0\}\)
%  (c.f. Lemma~\ref{lem:split}) is \(\mathcal{O}(1/d)\) close \klein{[It remains to be clarified which distance is meant here.]} to an invariant 
%  manifold $\mathcal{M}_{1/d}$ of the evolution described by~\eqref{eq:model3}. 
%  Furthermore, every solution starting sufficiently close to \(\mathcal{M}_{1/d}\) 
%  will approach \(\mathcal{M}_{1/d}\) asymptotically for long times.
%\end{cor}
%

\begin{cor}
Let \(d\) be sufficiently large and \(\varepsilon^{2} \in o(1/d)\). For every compact subset of the constraint manifold \(\mathcal{M}\) there exists a manifold \(\mathcal{M}_{1/d}\) which lies within \(\mathcal{O}(1/d)\) of this subset
  and is locally invariant under the evolution
  of~\eqref{eq:model3}. Furthermore every solution starting
  sufficiently close to \(\mathcal{M}_{1/d}\) will approach \(\mathcal{M}_{1/d}\).
\end{cor}

%  Hausdorff distance
%$$
%d(X,Y) = \max(\sup_{x\in X} \inf_{y\in Y} |x-y|, \sup_{y\in Y} \inf_{x\in X} |x-y|)
%$$

\begin{rem}
  Although the preceding corollary tells us there is at least one slow
  manifold for large \(d\) that satisfies the constraint, this does not imply
  we approach one of this kind, when starting slightly off the original
  constraint manifold \(\mathcal{M}\).
\end{rem}

So far we have only considered the analytical properties of the dissipative
system~\eqref{eq:model3}. Building upon the insights gained, we now propose a
related numerical method.

% ==============================================================================

\subsubsection{Dissipative discretization of the intermediate blended models}
\label{sssec:DissipativeDiscretization}

Here we will establish a relation between discrete solutions to the blended 
model from \eqref{eq:blendedmodel} and solutions of the dissipative surrogate model 
from \eqref{eq:model3} in the spirit of a modified equation 
analysis~\citep{WarmingHyett1974}. In other words, we will argue that trajectories 
produced by the blended method will locally relax to the constraint manifold by 
similar means as in the context of the numerical stabilization of solvers
for differential algebraic equations~\cite{ascher1994}. It will be useful to 
consider the following generalization of the constrained system~\eqref{eq:limiteq}, 
which we will call ``relaxed constraint system'' below, 
\begin{equation} \label{eq:modeldinfty}
  \begin{alignedat}{3}
    \dot{q} & = p				   & \qquad q(0) & = q_0 \in
    \mathbb{R}^{N}
    \\
    \dot{p} & = - G(q)\trp K \lambda - \nabla V(q) & \qquad p(0) & = p_0, \quad
    G(q_0)p_0 = 0
    \\
    g(q)    & =g(q_0)\,.
  \end{alignedat}
\end{equation}
\begin{rem}\label{rem:tangentialequivconstrained}
  After differentiating the constraint, yielding \(G(q)p = 0\), \(\lambda \) is given as
  before by~\eqref{equ:lagrangian}.
  For \(g(q_0)=0\) this system is equivalent to the constrained
  system~\eqref{eq:limiteq} in the sense that \(\mathcal{M}\) is invariant under
  the evolution in time following~\eqref{eq:modeldinfty}. Due to continuous dependency on initial data we
  furthermore conclude that solutions to~\eqref{eq:modeldinfty} approach
  solutions of~\eqref{eq:limiteq} as \(g(q_0) \rightarrow 0\).
\end{rem}

As pointed out previously, our method is supposed to be consistent
with the unconstrained system~\eqref{eq:model1} and system~\eqref{eq:limiteq}
for \(\alpha = 1\) and \(\alpha=0\), respectively. For the first case we
furthermore require high fidelity in our approximation of energy conservation 
for discrete solutions as well and choose the method to be 
symplectic~\cite{HairerEtAl2010}, i.e., such that the gradient of the discrete 
flow \(D\psi^1_h\) satisfies
\begin{equation}
  \label{eq:symplectic}
  (D\psi^1_h)\trp J (D\psi^1_h) = J \coloneqq
  \begin{pmatrix}0 & \id \\ -\id & 0\end{pmatrix}.
\end{equation}
This property is shared with the analytical flow and responsible for exact
conservation of the energy functional \(H\) for analytical solutions as well
as preservation of a modified, close by, energy functional for discrete 
solutions. For an extensive presentation and discussion on this topic see
e.g.~\cite{HairerEtAl2010} or~\cite{leimkuhler_reich_2005}.

We build our method on the symplectic (c.f.~\eqref{eq:symplectic}) St\"ormer
Verlet method,
\begin{equation}\label{eq:SV}
  \begin{aligned}
    q^{n+\frac{1}{2}} & = q^n + h \, p^{n}
    \\
    p^{n+1}           & = p^{n} 
    - h \nabla V(q^{n+\frac{1}{2}}))                           
    - \frac{h}{\varepsilon^2} {G\trp(q^{n+\frac{1}{2}})} K g(q^{n+\frac{1}{2}}) 
    \\
    q^{n+1}           & = q^{n+\frac{1}{2}} + h \, p^{n+1}\,.
  \end{aligned}
  %\tag{SV}
\end{equation}
Despite its simplicity this method performs exceptionally well and is 
extensively discussed in detail, e.g., in~\cite{HairerEtAl2010}.

In contrast to~\eqref{eq:model1}, the constrained model equations~\eqref{eq:limiteq}
are a system of differential algebraic equations; their differentiation index is~\(3\)~\cite{Ascher1998}.
Solving these equations numerically leaves the choice of either fulfilling the constraint
exactly or of accepting a numerical approximation error for \(g(q^n) = 0\).
In the Hamiltonian context, the first choice suggests, e.g., the 
SHAKE and RATTLE schemes (c.f.~\cite{leimkuhler_reich_2005}).
Given initial states \(q^n\) and tangential momenta \(p^n\), both algorithms
use a projection to ensure \(q^{n+1} \in \mathcal{M}\) and
\(p^{n+1} \in T_{q^{n+1}}\mathcal{M} \).

The alternative of accepting approximation errors for the constraint itself
relies on index reduction of the analytical system and subsequent
discretization (c.f.~e.g.~\cite{Ascher1998}).
In this context, a common task is to design stabilized methods~\cite{ascher1994}
which allow for a discrete evolution close to the constraint manifold 
such that the error on the constraint stays small for long times.

Since we do not aim to run the constrained model \(\alpha=0\) for more than 
a few time steps in the blended method, we will employ an index reduction
approach but ignore the issue of stabilization at this point.
Motivated by the St\"ormer-Verlet method we propose a projection method
for~\eqref{eq:limiteq} which satisfies the hidden constraint
\begin{align}
  G(q)p=0
\end{align}
up to a given tolerance and the constraint \(g(q)=0\) in~\eqref{eq:limiteq}
up to a global error of order \(\mathcal{O}(h^2)\).
The proposed discretization of the blended model \eqref{eq:blendedmodel} reads  
\begin{equation}\label{eq:bmethod}
  \begin{aligned}
    q^{n+\frac{1}{2}}  
      & = q^n + \frac{h}{2} p_n                   
        \\
    p^{n+1,\alpha}            
      & = p^n - h \nabla V(q^{n+\frac{1}{2}}) 
          - h G\trp(q^{n+\frac{1}{2}}) 
            K \left(\frac{\alpha}{\varepsilon^2} g(q^{n+\frac{1}{2}}) 
                  + (1 - \alpha)\lambda^{n+\frac{1}{2}} 
            \right)
        \\
    q^{n+1,\alpha}            
      & = q^{n+\frac{1}{2}} +\frac{h}{2} p^{n+1,\alpha}\,,  
%        \\
%    G(q^{n+1,\red{\alpha}}) p^{n+1,\red{\alpha}} 
%      & = 0.
  \end{aligned}
\end{equation}
where \(\lambda^{n+\frac{1}{2}}\) satisfies the weakly nonlinear equations
\begin{subequations}\label{eq:it}
\begin{align} 
p^{n+1,0}
  & = p^n 
    - h \nabla V(q^{n+\frac{1}{2}}) 
    - h G\trp(q^{n+\frac{1}{2}}) K \lambda^{n+\frac{1}{2}}                                                 
    \label{eq:it1}\\
h G\left(q^{n+\frac{1}{2}} 
+ \frac{h}{2}p^{n+1,0}\right) G\trp(q^{n+\frac{1}{2}}) \,K  \lambda^{n+\frac{1}{2}}
  & = G\left(q^{n+\frac{1}{2}} + \frac{h}{2} p^{n+1,0}\right) 
      \left(p^n - h \nabla V(q^{n+\frac{1}{2}})\right)\,.
    \label{eq:it2}
  \end{align}
  \end{subequations}
as derived from the relaxed constraint of tangential motion, i.e., $G(q)\, p = 0$. 

For this method, we now state and prove several consistency results:

\begin{lem}\label{lem:consistentSlow}
  Let \(\alpha = 0\) and \(\kappa \in \mathbb{N}\). Let \((q^\kappa,
  p^\kappa)\) be the numerical solution given by applying
  method~\eqref{eq:bmethod} \(\kappa \)~times to initial data \((q^0, p^0)
  \in \mathbb{R}^{2N}\), which satisfy \(G(q^0)p^0 = 0\). Then \((q^\kappa,
  p^\kappa)\) is consistent with the analytical solution \((q,p)\)
  of the constrained system~\eqref{eq:modeldinfty} at time \(T = h\kappa \) for initial condition
  \(q^0, p^0\).
  More specifically, 
  \begin{subequations}
  \begin{align}
    \|p^\kappa-p(T)\|      & \leq c h 	           \\
    \|q^\kappa-q(T)\|      & \leq \tilde{c} h^2  \\
    \|G(q^\kappa)p^\kappa \| & =0	                  \\
    \label{eq:errconstr}
    \|g(q(T))-g(q^\kappa)\|  & \leq \hat{c} h
  \end{align} 
  \end{subequations}
  where the constants \(\tilde{c},\hat{c},c\) are independent of \(h\) and \(\kappa \).
\end{lem}
\begin{proof}
  The proof is following~\cite{leimkuhler1994}. While first order consistency
  is essentially proven by a classical Taylor expansion argument, one still
  needs to address the algebraic constraint. At the continuous level this is
  readily achieved by reference to~\eqref{equ:lagrangian} which becomes
  \begin{align}
    \label{equ:analyticalexplicit}
    K \lambda  = \left(G(q)G(q)\trp\right)^{-1} 
    \left(- G(q) \nabla V(q)
          + \sum_{i=1}^L \sum_{j=1}^L  \left({p}\trp e_{i,j} p\right) \,
            \frac{\partial^2 g(q)}{\partial q_i \partial q_j}
    \right).
  \end{align}
  The momentum equation in~\eqref{eq:limiteq} is then equivalent to
  \begin{align}\label{equ:analyticclosed}
    \dot{p} = - \mathcal{P}_q^\perp \nabla V(q)  +
    G(q)\trp(G(q)G(q)\trp)^{-1} \sum_{i=1}^L \sum_{j=1}^L
    \left({p}\trp e_{i,j} p\right) \,  \frac{\partial^2 g(q)}{\partial q_i
      \partial q_j}.
  \end{align}
  For the discrete case we observe
  \begin{subequations}
  \begin{align}
    G(q^{n+1,0})p^n & = G(q^n)p^n + \sum_{i=0}^{N}\sum_{j=0}^{N} (q^{n+1,0}-q^n)
    e_{i,j} p^n \frac{\partial^2
      g(q^{n})}{\partial q_i \partial q_j } + \mathcal{O}(h^2)
    \\
                  & =	h \sum_{i=0}^{N}\sum_{j=0}^{N} \frac{p^{n+1,0}+p^n}{2}
    e_{i,j} p^n \frac{\partial^2
      g(q^{n+\frac12})}{\partial q_i \partial q_j } + \mathcal{O}(h^2),
  \end{align}
  \end{subequations}
  by Taylor expansion, where \(G(q^{n})p^n = 0 \) holds due to the
  tangential update of the previous time step or if applicable, by
  the initial condition in~\eqref{eq:bmethod}.
  Using this identity and the second and fourth update rules
  in~\eqref{eq:bmethod}
  we can express \(h K\lambda \) explicitly by
  \begin{subequations}
  \begin{align}
    hK\lambda = & \left( G\left(q^{n+1,0}\right)
    G\trp(q^{n+\frac{1}{2}})\right)^{-1}	G\left(q^{n+1,0}\right) \left(p_n -
    h \nabla V(q^{n+\frac{1}{2}}) \right)
    \\
    = & \ h {\left( G\left(q^{n+1,0}\right) G\trp(q^{n+\frac{1}{2}})\right)}^{-1}	
    \left(- G\left(q^{n+1,0}\right)  \nabla V(q^{n+\frac{1}{2}})
      + \left(\sum_{i=0}^{N}\sum_{j=0}^{N} \frac{p^{n+1,0}+p^n}{2} e_{i,j} p^n
      \frac{\partial^2 g(q^{n+\frac12})}{\partial q_i \partial q_j}
      \right) 
    \right)
    \\ \nonumber
    & 
    + \mathcal{O}(h^2).
%    \\
%    = & - h {\left( G\left(q^{n+1,0}\right)
%    G\trp(q^{n+\frac{1}{2}})\right)}^{-1}	G\left(q^{n+1,0}\right)  \nabla
%    V(q^{n+\frac{1}{2}})
%    \\ \nonumber
%    & + h {\left(G\left(q^{n+1,0}\right)
%        G\trp(q^{n+\frac{1}{2}})\right)}^{-1}
%    \left(\sum_{i=0}^{N}\sum_{j=0}^{N} \frac{p^{n+1,0}+p^n}{2} e_{i,j} p^n
%    \frac{\partial^2 g(q^{n+\frac12})}{\partial q_i \partial q_j
%    }
%    \right) + \mathcal{O}(h^2).
  \end{align}
  \end{subequations}
  Expanding by Taylor again, this identity now enables us to rewrite the
  momentum update to
  \begin{equation}
    \label{eq:explicittmomentupdate}
    \begin{aligned}
    p^{n+1,0} = p^n &- h \mathcal{P}_{q^{n+\frac{1}{2}}}^\perp  \nabla
    V(q^{n+\frac{1}{2}}) \\ &- h
    G\trp(q^{n+\frac{1}{2}})	\left( G(q^{n+\frac{1}{2}})
    G\trp(q^{n+\frac{1}{2}})\right)^{-1}
    \sum_{i=0}^{N}\sum_{j=0}^{N}
    \frac{p^{n+1,0}+p^n}{2} e_{i,j} p^n \frac{\partial^2 g(q^{n+\frac12})}{\partial q_i \partial
      q_j } + \mathcal{O}(h^2).
    \end{aligned}
  \end{equation}
  Comparing the momentum in closed form as given in~\eqref{equ:analyticclosed}
  with the discretization~\eqref{eq:explicittmomentupdate}, we verify second
  order local consistency and therefore first order global consistency
  \begin{align}
    \|p^{\kappa} - p(T)\| \leq	c h\,.
  \end{align}
  To bound the consistency error in the coordinates one combines the update
  rules in~\eqref{eq:bmethod} to
  \begin{align}
    q^{n+1,0} = q^n + h \frac{p^n + p^{n+1,0}}{2}.
  \end{align}
  Since \(p_{n+1}\) this a locally second order consistent approximation we obtain
  \begin{align}
    q^{n+1,0} = q^n + h \frac{p(t^n) + p(t^{n+1})}{2} + \mathcal{O}(h^3).
  \end{align}
  The implicit midpoint rule is of local third consistency order and we obtain
  \begin{align}
    \|q^{n+1,0}-q(t^{n+1})\| =  \mathcal{O}(h^3)
  \end{align} and therefore
  \begin{align}
    \|q^\kappa-q(T)\| \leq	\tilde{c} h^2.
  \end{align}
  The bound for the constraint~\eqref{eq:errconstr} follows then directly by
  expanding \(g(q(T))\) around \(g(q^n)\).
\end{proof}
\begin{cor}
  Let additionally \(g(q^0)=0 \), then \((q^\kappa,p^\kappa) \) is consistent
  with the constraint system~\eqref{eq:limiteq}.
\end{cor}

So far we consider only initial data with momenta satisfying tangency to the 
constraint manifold, i.e., \(G(q)p=0\). This is of course necessary in the 
context of consistency, since the underlying model is not well posed otherwise. 
Nevertheless the proposed usage in data assimilation procedures introduces 
exactly such initial data. The subsequent two statements will clarify what 
to expect if we apply method~\eqref{eq:bmethod} to general initial data 
while \(\alpha = 0 \).
\begin{lem}
  For \(\alpha = 0\) the method in~\eqref{eq:bmethod} approximates the
  projection of momentum \(\mathcal{P}^\perp_q p\) in the following sense.
  \begin{align}
    \|\mathcal{P}^\perp_{q^n} {p^n} - p^{n+1,0} \| \leq c h
  \end{align}
\end{lem}
\begin{proof}
  Use of the expression for \(K\lambda \) as stated in~\eqref{eq:it2} and
  subsequent Taylor expansion yield
  \begin{subequations}
  \begin{align}
    p^{n+1,0} & = p^n - G\trp(q^{n+\frac{1}{2}})  {\left(G\left(q^{n+1,0}\right)
        G\trp(q^{n+\frac{1}{2}})\right)}^{-1}
    G\left(q^{n+1}\right) p^n  + \mathcal{O}(h)
    \\
            & = p^n - G\trp(q^{n})  {\left(G\left(q^{n}\right)
        G\trp(q^{n})\right)}^{-1}     G\left(q^{n}\right)p^n +
    \mathcal{O}(h)                                                         \\
            & = \mathcal{P}^\perp_{q^n}p^n + \mathcal{O}(h).
  \end{align}
  \end{subequations}
\end{proof}
\begin{cor}
  Let \(\alpha = 0\) and \(q^0,p^0 \in \mathbb{R}^N\), then method~\eqref{eq:bmethod} is globally
  first order consistent to the solution given by the constraint system~\eqref{eq:modeldinfty} and
  balanced initial data \(q^0\) and \(\mathcal{P}^\perp_{q^0}p^0\).
\end{cor}

%\begin{exm}[Harmonic oscillator]
%  Let \(q=(\mu,\nu)^T\), \(p=(\eta,\zeta)^T\),
%  \(K=\mathrm{diag}(\varepsilon^2,1)\) and \(V(q) =0 \)
%  \begin{alignat}{3}
%    \label{eq:harmonicslow}
%    \dot{\mu} & = \eta	 \qquad   &             & \dot{\eta} &
%              & = - \mu                                          \\
%    \label{eq:harmonicfast}
%    \dot{\nu} & = \zeta  \qquad
%              &                 & \dot{\zeta} &            & = -
%    \frac{1}{\varepsilon^2} \nu
%  \end{alignat}
%  This system has no coupling between the fast and slow components and
%  furthermore solutions will stay on the constraint manifold \(\mathcal{M} =
%  \{\nu = g(q) = 0\}\) if initially so and additionally fulfill \(\zeta = 0\).
%  For any initial values \(q_0=(\mu_0,\nu_0)^T\), \(p_0=(\eta_0,\zeta_0)^T\)
%  the numerical method introduced in \eqref{eq:tmethod} in this example is a
%  consistent approximation to
%  \begin{alignat}{3}
%    \dot{\mu} & = \eta	 \qquad &       & \dot{\eta} &        & = - \mu \\
%    \dot{\nu} & = 0  \qquad
%              &               & \zeta &            & = 0  .
%  \end{alignat}
%\end{exm}

Next we establish a consistency result for the blending method for $\alpha \not= 0$,
which will provide a connection between the discrete evolution of the blended model 
from~\eqref{eq:blendedmodel} on the one hand, and the continuous evolution of the damped
surrogate model \eqref{eq:model3} on the other hand.
\begin{lem}\label{lem:consistentBlending}
  Let \(\kappa \in \mathbb{N}\) and let \(\alpha = \max(0,1-dh)\).
  Furthermore let \((q^\kappa,p^\kappa)\) be the numerical solution given by
  applying~\eqref{eq:bmethod} \(\kappa\) times to initial data
  \((q^0,p^0)\in\mathbb{R}^{2N}\).
  Then \((q^\kappa,p^\kappa)\) is consistent with the solution \((q,p)\) of the
  dissipative
  system~\eqref{eq:model3} at time \(T=h\kappa\).
  More specifically, % there holds
  \begin{subequations}
  \begin{align}
    \|q^\kappa-q(T)\| & \leq c h^2       \\
    \|p^\kappa-p(T)\| & \leq \tilde{c} h
  \end{align} 
  \end{subequations}
  where \(c,\tilde{c}\) is a constant independent of \(h\) and
  \(\kappa \).
\end{lem}
\begin{proof}
  We start by expressing \(\lambda \) explicitly and rewriting the momentum
  update as
  \begin{subequations}
  \begin{align}
    p^{n+1,\alpha} 
    		& = p_n - h \nabla V(q^{n+\frac{1}{2}}) 
             - \alpha \, \frac{h}{\varepsilon^{2}} {G\trp(q^{n+\frac{1}{2}})} K g(q_{n+\frac{1}{2}}) 
    		 - (1-\alpha) \, h	{G\trp(q^{n+\frac{1}{2}})} K \lambda
    		\\
            & = p^n 
              - h \nabla V(q^{n+\frac{1}{2}}) 
              - \alpha \, \frac{h}{\varepsilon^{2}}  {G\trp(q^{n+\frac{1}{2}})} K g(q^{n+\frac{1}{2}})
    		\\
            & \nonumber \phantom{= p^n }	
              - (1-\alpha) \, {G\trp(q^{n+\frac{1}{2}})} \left( G\left(q^{n+1}\right)
    			G\trp(q^{n+\frac{1}{2}})\right)^{-1} G\left(q^{n+1,\alpha}\right) 
			   \left(p_n - h \nabla V(q^{n+\frac{1}{2}}) \right)
    \\
          	& =  p^n 
	          - h (\tilde{\mathcal{P}}^\perp_{q^{n+\frac{1}{2}}} 
                   + \alpha \tilde{\mathcal{P}}_{q^{n+\frac{1}{2}}}) \nabla V(q^{n+\frac{1}{2}}) 
              - \frac{h}{\varepsilon^{2}} \alpha G\trp(q^{n+\frac{1}{2}}) K g(q^{n+\frac{1}{2}}) 
              - (1-\alpha)\tilde{\mathcal{P}}_{q^{n+\frac{1}{2}}} p^n
  \end{align} 
  \end{subequations}
  where 
  \(\tilde{\mathcal{P}}_{q^{n+\frac{1}{2}}} =
  G\trp(q^{n+\frac{1}{2}}) {\left(G\left(q^{n+1,\alpha}\right)
  G\trp(q^{n+\frac{1}{2}})\right)}^{-1}  G\left(q^{n+1,\alpha}\right)\) satisfies by
  Taylor expansion
  \(\tilde{\mathcal{P}}_{q_{n+\frac{1}{2}}} = \mathcal{P}_{q_{n+\frac{1}{2}}} +
  \mathcal{O}(h)\).
  For sufficiently small \(h\) we now can substitute \(\alpha = \max(0,1-hd)\)
  by \(1-hd\) since \(d\) is independent of $h$. Expanding by Taylor we find first order  
  consistency with the dissipative model \eqref{eq:model3} in the sense of modified equation 
  analysis,
  \begin{subequations}
  \begin{align}
    p^{n+1,\alpha} 
      & = p^n - h (\tilde{\mathcal{P}}^\perp_{q^{n+\frac{1}{2}}} 
        + (1 - h d) \tilde{\mathcal{P}}_{q^{n+\frac{1}{2}}}) \nabla V(q^{n+\frac{1}{2}}) 
        - \frac{h}{\varepsilon^{2}} (1 - h d) G\trp(q^{n+\frac{1}{2}}) K g(q^{n+\frac{1}{2}}) 
        - h d \tilde{\mathcal{P}}_{q^{n+\frac{1}{2}}} p_n
        \\
      & = p^n 
        - h  
          \left( 
            \varepsilon^{-2} G\trp(q^{n+\frac{1}{2}}) K g(q^{n+\frac{1}{2}}) 
          + d\, \mathcal{P}_{q^{n+\frac{1}{2}}} p^n 
          + \nabla V(q^{n+\frac{1}{2}}) 
          \right)
        + \mathcal{O}(h^2).
  \end{align}
  \end{subequations}
  The update of the coordinates can be rewritten as
  \begin{equation}
    q^{n+1} = q^n + h\frac{p^n + p^{n+1}}{2}
  \end{equation} 
  and therefore the same proof as in
  Lemma~\ref{lem:consistentSlow} leads to the statement.
\end{proof}

\begin{rem}
  Although it may seem odd at a first glance that we have to change 
  \(\alpha\) depending on \(h\) to achieve convergence to solutions
  of the dissipative surrogate model, this
  does not contradict the fact that given an a priori choice of \(\alpha \),
  discrete solutions of the blended method are consistent approximations to
  solutions of system~\eqref{eq:model3} for a certain \(d\). Thus, the
  discretized blended model inherits the tendency to approach the balanced
  manifold from the dissipative model as desired. The relationship between
  $\alpha, d$, and $h$ is also what motivates our reference to the modified
  equation approach. 
\end{rem}
To complete and illustrate the overall picture we collect most of the
preceding results and references for the proof of the following commuting
diagram.
\begin{prop}\label{prop:diagram}
  Let \(\alpha=\max(0,1-hd)\). Let \(\phi \) and \(\psi \) be the analytical
  and
  numerical flows respectively, with regard to the models as mentioned below,
  then the diagram in Figure~\ref{fig:commutingdiagram} commutes.
  \begin{figure}[!h]
    \begin{tikzcd}[row sep=large,column sep=large]
      \psi_h^{1} \arrow[r,leftarrow,"d \rightarrow 0"] \arrow[d,"h
        \rightarrow 0"] & \psi_h^{\alpha} \arrow[r,"d
        \rightarrow \infty"] \arrow[d,"h \rightarrow 0"]  & \psi_h^{0}
      \arrow[d,"h \rightarrow 0"]&  \\
      \phi^{\varepsilon} \arrow[r,leftarrow,"d \rightarrow 0 "] \arrow[rrr,
        bend right,"\varepsilon \rightarrow 0"]
      & \tilde{\phi}^{d} \arrow[r,"d \rightarrow \infty"] &
      \tilde{\phi}^{\infty}  \arrow[r,Leftrightarrow, "\text{on }
        \mathcal{M}"] & \phi^0 \\
    \end{tikzcd}
    \centering
    \caption{The commuting diagram shows the connections between the analytical
      model hierarchy given by the flow \(\phi^\epsilon \) of the Hamiltonian
      system~\eqref{eq:model1}, the flow \(\tilde{\phi}^d\)
      of the dissipative system~\eqref{eq:model3} and the
      flow \(\tilde{\phi}^\infty \) of the relaxed constrained
      model~\eqref{eq:modeldinfty}.
      Furthermore the diagram depicts the consistency results, of the
      St\"ormer-Verlet method~\eqref{eq:SV} and the blended
      method~\eqref{eq:blending} denoted by \(\psi_h^\alpha \).
      The flow for the constrained model~\eqref{eq:limiteq} is
      denoted by \(\psi^0\). }\label{fig:commutingdiagram}
  \end{figure}
\end{prop}
\begin{proof}
  The consistency of the Störmer-Verlet method is stated in
  e.g.~\cite{HairerEtAl2010} and for an overview of all the other connections
  in
  the commuting diagram in Figure~\ref{fig:commutingdiagram} we refer to
  Figure~\ref{fig:proofcommutingdiagram}.
  \begin{figure}[!h]
    \begin{tikzcd}[row sep=large,column sep=large]
      \psi_h^{1} \arrow[r,leftarrow,"\text{\eqref{eq:bmethod}}"]
      \arrow[d,"\text{\cite{HairerEtAl2010}}"] & \psi_h^{\alpha}
      \arrow[r,"\text{\eqref{eq:bmethod}}"]
      \arrow[d,"\text{Lem.~\ref{lem:consistentBlending}}"]  & \psi_h^{0}
      \arrow[d,"\text{Lem.~\ref{lem:consistentSlow}}"]&
      \\
      \phi^{\varepsilon} \arrow[r,leftarrow,"\text{Lem.~\ref{lem:contdep}}"]
      \arrow[rrr, bend
        right,"\text{\cite{RubinUngar1957}}"] & \tilde{\phi}^{d}
      \arrow[r,"\text{Lem.~\ref{lem:singpert}}"] &
      \tilde{\phi}^{\infty}  \arrow[r,Leftrightarrow, "\text{Rem.~\ref{rem:modelEquiv}}"] & \phi^0 \\
    \end{tikzcd}
    \centering
    \caption{The diagram depicts the same situation as in
      Figure~\ref{fig:commutingdiagram}, but refers to the previously
      established
      results and relevant literature, instead of the limits.
    }\label{fig:proofcommutingdiagram}
  \end{figure}
\end{proof}

% ==============================================================================
% ==============================================================================
% ==============================================================================

\section{Numerical Results}
\label{sec:results}

For experiments in the context of data assimilation one immediate obstacle
arises from potential model errors.
We avoid this question by considering an initially balanced reference 
solution of~\eqref{eq:model1} which is approximated by the St\"{o}rmer-Verlet method~\eqref{eq:SV}.
Henceforth this solution will be denoted by \(z^{\text{ref}}\).
The observations then are given by \(y_{\rm{obs}}(t_k) = H_{\rm{obs}}z^{\text{ref}}(t_k) + \zeta_k\).
Hereby \(H_{\rm{obs}} z(t_k) = q(t_k)\) and \(\zeta_k\) is the realization of the normally distributed measurement error at
some time \(t_k = k \Delta t_{\text{obs}}\), when the observation becomes available.
We assume the measurement error to have zero mean and covariance \(R=\rho I\).
The resulting evolution of observations is assimilated by the proposed data
assimilation scheme.
The advantage of this setup is the straightforward assessment of the quality of
the data assimilation method by comparing the reference solution to, e.g., the
ensemble members or the point estimate of their mean.

According to (c.f.~\cite{reichbook}), due to finite ensemble sizes the true
covariances of the posterior distributions are underestimated in ensemble
based data assimilation methods. One technique to address this issue is
ensemble inflation which amounts to an artificial increase of the spread of
the ensemble after each assimilation step by
\begin{align}
  z_i^{\rm new}:=\bar{z}+\sigma_{\text{infl}} (z_i-\bar{z}).
\end{align}
In the experiments we apply the ensemble inflation as the last step of the
assimilation procedure. For the comparison of the presented methods we choose
again the stiff elastic double pendulum from Example~\ref{ex:doublependulum}
as the dynamical model. The initial ensemble is constructed from copies of
the initial state of \(z^{\rm ref}\) except for the components tangential to
the constraint manifold \(\mathcal{M}\). For those we perturb the reference
state normally distributed with covariance \(\rho_0\). Subsequently we
eliminate the resulting normal component i.e. balance the initial data by
minimizing~\eqref{eq:costfunctional} for \(B=0\) and \(\gamma=0\). This is a
natural modification of the penalty method for the first step and ensures the
samples are spread along the constraint manifold only. In comparing the
different balancing methods, the initial data are always generated in the way
just described.

For the blended time stepping method we choose a linear ramp for \(\alpha \)
as depicted in Fig.~\ref{fig:numericalmethod} where \(\alpha=0\) initially and
\(\alpha=1\) at the end of the blending window. The analysis
of the damped model equation in \eqref{eq:model3} is based on linearization 
and suggests that, similarly to the situation with the harmonic oscillator, 
one can find values of the blending and damping parameters \(\alpha \) and 
\(d\), respectively, that imply dynamics close to the aperiodic case.
As we do not further investigate the question for optimal \(\alpha\) we
choose a linear ramp to step through different values of the damping
coefficient as brute force approach.

For the numerical values of the parameters of the experiments we refer to
Table~\ref{tab:double_parameters}.
\begin{table}[h]
  \centering
  \begin{tabular}{cccccccccccccc}
    \toprule
    \(B\) & 
    \(L\) & 
    \(l\) & 
    {\(\varepsilon \)} & 
    \(K\) & 
    {\(a_0\)}  & 
    {\(\Delta t\)} & 
    \(M\) & 
    \(H_{\rm{obs}} z\) & 
    {\(\Delta t_{\rm obs}\)} & 
    \(\rho\) & 
    \(\rho_0\)  & 
    \(\sigma_{\text{infl}}\) & 
    T \\ \midrule
    \(\id\) & 
    \(2\) & 
    \((1,1)\) & 
    \(0.001\) & 
    \( {\rm diag}(1,0.04)\) & 
    \(9.81\) & 
    \(0.001\) & 
    \(20\) & 
    \(q\) & 
    \(0.1 \) & 
    \(0.05\) & 
    \(0.05\) & 
    \(1.05\) & 
    \(500\) \\ \bottomrule
  \end{tabular}
  \caption{
    Parameters for the numerical experiments using the double pendulum
    model.%
    The model parameters are given by the equilibrium lengths \(l \in
    \mathbb{R}^L\),
    the scale separation parameter \(\varepsilon \), the stiffness matrix \(K\)
    and the gravity \(a_0\). 
    The model is discretized by the St\"ormer-Verlet method~\eqref{eq:SV} with
    step width \(\Delta t\).
    \(M\) denotes the ensemble size, \(\Delta t_{\text{obs}}\) the interval
    between two observations and \(Hz\) the observed variable.
    We choose \(\rho \) as covariance of the measurement error, \(\rho_0 \)
    as the initial uncertainty i.e.~the covariance of the initial ensemble and
    \(\rho_\text{infl}\) as the inflation factor.
    Finally \(T \) denotes the duration of the experiment.
  }\label{tab:double_parameters}
\end{table}
As with regard to the implementation details,
we minimize the functional~\eqref{eq:costfunctional} using either the 
Broyden-Fletcher–Goldfarb-Shanno (BFGS)~\cite{broyden_convergence_1970,fletcher_new_1970,goldfarb_family_1970,shanno_conditioning_1970} method as implemented in scipy~\cite{2020SciPy-NMeth} 
or the proposed algorithm of~\eqref{eq:3DVarContDiscrete}.
For the first we require a tolerance of \(10^{-8}\) and as initial values we
choose the results of the plain EnKF. In the second case we choose a fixed step size of \(h=10^{-3} \) and iterate as long as
the maximal absolute value of the increment~\eqref{eq:increment} exceeds
\(10^{-8}\). 
Additionally we need to solve a nonlinear system for the implicit part of the
blended time stepping method~\eqref{eq:bmethod} when \(\alpha=0\). This system 
is solved using the scipy~\cite{2020SciPy-NMeth}
wrapper for the modified Powell method from the MINPACK~\cite{More:126569} subroutine
\textit{hybrd}. The initial value is the zero vector of dimension \(L\) and
the tolerance for the nonlinear problem is set to double precision i.e.
\(10^{-16}\).
To quantify the error of the methods we use the time averaged root mean square
error as given in~\cite{reichbook}
\begin{equation}\label{eq:rmse}
  \mathrm{TRMSE}(Z)=\sqrt{\frac{1}{N_T}\sum_i^{N_T}{\|\hat{Z}(t_i)-Z(t_i)\|}^2}.
\end{equation}
Hereby \(\hat{Z}\) denotes the estimate for the quantity \(Z\) and both are
evaluated at \(N_T\) time points \(t_k\in [0,T]\).

This score is shown below as a function of the tuning
parameters of the respective method. For comparison we furthermore show the
results for the unmodified ensemble Kalman filter. As seen in
Figures~\ref{fig:res_penalty_coordinates}~--~\ref{fig:res_blending}, the
forecast quality for the coordinates and the momenta improve drastically when
choosing appropriate tuning parameters for the respective methods.

For the penalty method we realize from
Figures~\ref{fig:res_penalty_coordinates} and~\ref{fig:res_penalty_momenta}
that we obtain the best results when forcing the analysis balance residual
of each ensemble member to be close to the respective one inferred from the
forecast. We can enforce this by the penalty method when setting \(\gamma =
1\). We also find that increased weights do add to the forecast quality only
up to certain extent.

The blending method only allows for one tuning parameter, the blending window
size. Comparing several choices in Figure~\ref{fig:res_blending}, we 
obtain the best results when choosing a window large enough to capture a
full period of the less stiff spring in the blending window. This happens
approximately around \(\frac{2 \pi \varepsilon}{\sqrt{k_2}}\approx 0.3 \Delta
t_{\mathrm{obs}}\).
\begin{figure}
  \input{\imgPath results_penalty.pgf}
  \caption{The left panel depicts the time averaged root mean square error
  (TMRSE) of the coordinates obtained by the penalty method 
  minimizing the functional~\eqref{eq:costfunctional} and 
  using the previously mentioned BFGS
  solver. The right panel shows the results for the same experiment,
  but using the penalty method solved by the descent with modified search
  direction from~\eqref{eq:3DVarContDiscrete}.
  In orange we depict the results obtained by the unmodified ensemble Kalman filter. }\label{fig:res_penalty_coordinates}
\end{figure}
\begin{figure}
  \input{\imgPath results_vel_penalty.pgf}
  \caption{The left and right panels show the time averaged root mean square
  error in the tangential component of the unobserved momenta, for the
  penalty method solved by the BFGS and again~\eqref{eq:3DVarContDiscrete}
  respectively. In orange we depict the results obtained by the unmodified
  ensemble Kalman filter.}\label{fig:res_penalty_momenta}
\end{figure}
\begin{figure}
  \input{\imgPath results_blending.pgf}
  \input{\imgPath results_vel_blending.pgf}
  \caption{
  The left panel shows the time averaged root mean square error for the
  coordinates obtained by the blending method. The right panel displays the
  same for the tangential component of the unobserved momenta. For comparison
  the results obtained by the unmodified ensemble Kalman filter are shown in
  orange. The abscissa describes the ratio of the length of the blending
  window and the observation interval.
  }\label{fig:res_blending}
\end{figure}

% ==============================================================================
% ==============================================================================
% ==============================================================================

\section{Conclusions}
\label{sec:conclusions}

Aiming to improve data assimilation for slow solutions of highly 
oscillatory systems, this paper suggests two principally different extensions 
of ensemble-based data assimilation algorithms. The first approach modifies 
the data assimilation scheme itself and consists of a rather generic post 
processing step involving the minimization of a cost functional that quantifies 
the oscillatory solution content.

The second approach utilizes the ability of asymptotically consistent 
numerical schemes which provide seamless access to both the highly oscillatory 
systems of interest and to the reduced differential-algebraic counterparts 
describing motions on the associated slow manifolds. Following ideas 
first formulated in~\cite{BenacchioEtAl2014}, this method filters 
oscillatory, off-manifold, components arising in the course of an 
assimilation step by starting the subsequent forward simulation with a few 
time steps of the reduced dynamics and blending the solver back to the 
full system over another couple of time steps. Whereas the full and reduced
dynamics are represented by (nearly) energy preserving integrators, the 
intermediate systems accessed during the blending phase are designed here to 
selectively dissipate the oscillatory solution components. This prohibits the
re-introduction of oscillations in the course of the blending procedure. A 
rigorous justification of the blended time-stepping method by asymptotic 
analysis is provided. The optimal parameter choice within the blending time 
window remains a topic for further investigation.

It is demonstrated that both methods perform well in terms of forecast quality 
and allow accurate state estimation in situations where the standard ensemble
Kalman Filter fails to do so. The dependency of the forecast skill on 
the respective tuning parameters behaves as expected in our prototypical test
case of the elastic double pendulum. Both approaches leave room for further
improvement and extension, however. Thus, e.g., seamless incorporation of 
the balancing step in a Bayesian filter may help optimizing both in terms of 
accuracy, balancing quality, and efficiency.  Also, as one referee of this 
paper pointed out, the implicit particle filter technique of A.~Chorin and
co-workers, \cite{ChorinTu2009,ChorinEtAl2010}, may provide a means of 
effectively steering an entire ensemble to predominantly sample balanced 
states. This would be achieved by penalizing the probability of unbalanced 
states within the space of probability distributions which the method has 
access to. For the blending approach, aside from optimized sequences 
of the blending parameter, additional gains are conceivable when the fast part
of the dynamics is known to be linear, as is the case in atmospheric flow 
applications.

The broader application area for the two proposed stabilization techniques is
ensemble-based data assimilation for geophysical processes. This application
area shares the situation of small oscillatory energy and conservative motion
along a slow manifold. Ongoing research therefore investigates the efficacy of 
the proposed methods for data assimilation into multidimensional geophysical 
flow models. Depending on the pertinent spatio-temporal scales, several 
different dominant balances emerge in those models \citep{Klein2010}, with 
geostrophic balance as a prominent example. In contrast to the present work, 
these balances are often essentially linear so that Lemma~\ref{lem:invariant} 
applies and only weak generation of imbalances by the standard EnKF is expected. 
Considerable imbalances are introduced, however, by spatial localization in the 
assimilation algorithm, a measure that is used to avoid artificial global scale 
correlations. Localization destroys the linearity of the filter transformation 
and therefore gives rise to stronger imbalances again. The methodologies proposed 
here both directly translate to this context, since neither of the algorithms 
leverages the linearity of the filter. Ongoing studies investigate these issues 
for the rotational shallow water equations as well as for a vertical slice model 
of the atmosphere.

% ==============================================================================
% ==============================================================================
% ==============================================================================

\section*{Acknowledgments}	
This research has been partially funded by Deutsche Forschungsgemeinschaft
(DFG) through grant CRC 1114 ``Scaling Cascades in Complex Systems'', Project
Number 235221301, Project A02 ``Multiscale data and asymptotic model
assimilation for atmospheric flows''.
\bibliographystyle{abbrvnat}
\bibliography{references}
\end{document}

%% file: image/doublependulum.eps_tex
%% Creator: Inkscape 1.0beta1 (32d4812, 2019-09-19), www.inkscape.org
%% Accompanies image file 'doublependulum.eps' (pdf, eps, ps)
\begingroup%
  \makeatletter%
  \providecommand\color[2][]{%
    \errmessage{(Inkscape) Color is used for the text in Inkscape, but the package 'color.sty' is not loaded}%
    \renewcommand\color[2][]{}%
  }%
  \providecommand\transparent[1]{%
    \errmessage{(Inkscape) Transparency is used (non-zero) for the text in Inkscape, but the package 'transparent.sty' is not loaded}%
    \renewcommand\transparent[1]{}%
  }%
  \providecommand\rotatebox[2]{#2}%
  \newcommand*\fsize{\dimexpr\f@size pt\relax}%
  \newcommand*\lineheight[1]{\fontsize{\fsize}{#1\fsize}\selectfont}%
  \ifx\svgwidth\undefined%
    \setlength{\unitlength}{66.20337604bp}%
    \ifx\svgscale\undefined%
      \relax%
    \else%
      \setlength{\unitlength}{\unitlength * \real{\svgscale}}%
    \fi%
  \else%
    \setlength{\unitlength}{\svgwidth}%
  \fi%
  \global\let\svgwidth\undefined%
  \global\let\svgscale\undefined%
  \makeatother%
  \begin{picture}(1,1.66077704)%
    \lineheight{1}%
    \setlength\tabcolsep{0pt}%
    \put(0,0){\includegraphics[width=\unitlength]{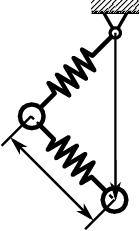}}%
    \put(0.24576073,0.19783253){\color[rgb]{0,0,0}\makebox(0,0)[lt]{\lineheight{1.25}\smash{\begin{tabular}[t]{l}$l_i$\end{tabular}}}}%
    \put(0.54471715,0.68944741){\color[rgb]{0,0,0}\makebox(0,0)[lt]{\lineheight{1.25}\smash{\begin{tabular}[t]{l}$k_i$\end{tabular}}}}%
    \put(0.88622278,0.77395953){\color[rgb]{0,0,0}\makebox(0,0)[lt]{\lineheight{1.25}\smash{\begin{tabular}[t]{l}$r_i$\end{tabular}}}}%
  \end{picture}%
\endgroup%

%% file: image/blending.tex
\begin{tikzpicture}
	\begin{groupplot}[
			group style={
				group size=2 by 1,
				xticklabels at=edge bottom,
				horizontal sep=0pt
			},
			scale only axis,
			hide y axis,
			axis lines=middle,
			enlargelimits=false,
			restrict x to domain=-3:8, xmin=-2.5, xmax=7.5,
			restrict y to domain=0:2.5, ymin=0, ymax=2.5,
			x=1.25cm,
			y=1.25cm
		]
		\nextgroupplot[
			restrict x to domain=-3:3,  xmin=-2.5,xmax=3.0,
			restrict y to domain=0:2.5,  ymin=0, ymax=2.5,
%			x=2cm,
%			y=2cm,
			xtick={-2,-1,0.001,1,2},
			tick style= thick,
			major tick length=0.625cm,
			xticklabels={\(t_{n-2}\),\(t_{n-1}\),\(t_{n}\),\(t_{n+1}\),\(t_{n+2}\)},
			extra x ticks={-2,-1,0,1,2},
			extra x tick labels={\(\psi^1_{h}\),\(\psi^1_{h}\),\(\psi^{\alpha_1=0}_{h}\),
						\(\psi_h^{\alpha_2}\),
						\dots},
			extra x tick style={
				major tick length=0cm,
				xticklabel style={xshift=0.625cm, anchor=south}
			},
			x axis line style={-}
		]
		{
			\addplot[thick, samples=50, smooth,domain=-1:1,red] coordinates {(0,0)(0,1.5)};
			%\addplot[thin,smooth,domain=-3:3] {1.0+0.5*sin((x+4)*16*pi)};
			\node[] at (axis cs:-2,1.5){forecast};
			\draw[-{Latex[length=3mm,width=2mm]}] (axis cs:-2.5, 1.25) -- (axis cs: 0, 1.25);
			\node[] at (axis cs:0,1.75){analysis};
			\node[] at (axis cs:2,1.5){forecast};
			\draw[-] (axis cs:0, 1.25) -- (axis cs: 3, 1.25);
			\node[] at (axis cs:2,0.9){blending};
			\draw[{Latex[length=3mm,width=2mm]}-] (axis cs:0, 0.65) -- (axis cs: 3, 0.65);

		}
		\nextgroupplot[
			restrict x to domain=5:8, xmax=7.5, xmin=5,axis x discontinuity=crunch,
			restrict y to domain=0:2.5, ymax=2.5, ymin=0,
%			x=2cm,
%			y=2cm,
			xtick={5,6,7},
			tick style= thick,
			major tick length=0.625cm,
			xticklabels={,\(t_{n+k}\),\(t_{n+k+1}\)},
			extra x ticks={5,6,7},
			extra x tick labels={ \dots,
						\(\psi_h^{\alpha_k=1}\),
						\(\psi_h^{1}\)},
			extra x tick style={
				major tick length=0cm,
				xticklabel style={xshift=0.625cm, anchor=south}
			},
			x axis line style={-{Latex[length=3mm,width=2mm]}}
		]{
			\addplot[thick, samples=50, smooth,domain=-1:1,black] coordinates {(7,0)(7,0.9)};

			\draw[-{Latex[length=3mm,width=2mm]}] (axis cs:5, 0.65) -- (axis cs: 7, 0.65);
			\draw[-{Latex[length=3mm,width=2mm]}] (axis cs:3, 1.25) -- (axis cs: 7.5, 1.25);

			%\addplot[thin,smooth,domain=5:7.5] {1.0+0.5*sin((x+2)*16*pi)};
		}
	\end{groupplot}
\end{tikzpicture}

%% file: image/results_blending.pgf
%% Creator: Matplotlib, PGF backend
%%
%% To include the figure in your LaTeX document, write
%%   \input{<filename>.pgf}
%%
%% Make sure the required packages are loaded in your preamble
%%   \usepackage{pgf}
%%
%% and, on pdftex
%%   \usepackage[utf8]{inputenc}\DeclareUnicodeCharacter{2212}{-}
%%
%% or, on luatex and xetex
%%   \usepackage{unicode-math}
%%
%% Figures using additional raster images can only be included by \input if
%% they are in the same directory as the main LaTeX file. For loading figures
%% from other directories you can use the `import` package
%% and then include the figures with
%%   \import{<path to file>}{<filename>.pgf}
%%
%% Matplotlib used the following preamble
%%   \usepackage{fontspec}
%%   \setmainfont{DejaVuSerif.ttf}[Path=/home/gottfried/.local/lib/python3.8/site-packages/matplotlib/mpl-data/fonts/ttf/]
%%   \setsansfont{DejaVuSans.ttf}[Path=/home/gottfried/.local/lib/python3.8/site-packages/matplotlib/mpl-data/fonts/ttf/]
%%   \setmonofont{DejaVuSansMono.ttf}[Path=/home/gottfried/.local/lib/python3.8/site-packages/matplotlib/mpl-data/fonts/ttf/]
%%
\begingroup%
\makeatletter%
\begin{pgfpicture}%
\pgfpathrectangle{\pgfpointorigin}{\pgfqpoint{3.500000in}{3.000000in}}%
\pgfusepath{use as bounding box, clip}%
\begin{pgfscope}%
\pgfsetbuttcap%
\pgfsetmiterjoin%
\definecolor{currentfill}{rgb}{1.000000,1.000000,1.000000}%
\pgfsetfillcolor{currentfill}%
\pgfsetlinewidth{0.000000pt}%
\definecolor{currentstroke}{rgb}{1.000000,1.000000,1.000000}%
\pgfsetstrokecolor{currentstroke}%
\pgfsetdash{}{0pt}%
\pgfpathmoveto{\pgfqpoint{0.000000in}{0.000000in}}%
\pgfpathlineto{\pgfqpoint{3.500000in}{0.000000in}}%
\pgfpathlineto{\pgfqpoint{3.500000in}{3.000000in}}%
\pgfpathlineto{\pgfqpoint{0.000000in}{3.000000in}}%
\pgfpathclose%
\pgfusepath{fill}%
\end{pgfscope}%
\begin{pgfscope}%
\pgfsetbuttcap%
\pgfsetmiterjoin%
\definecolor{currentfill}{rgb}{1.000000,1.000000,1.000000}%
\pgfsetfillcolor{currentfill}%
\pgfsetlinewidth{0.000000pt}%
\definecolor{currentstroke}{rgb}{0.000000,0.000000,0.000000}%
\pgfsetstrokecolor{currentstroke}%
\pgfsetstrokeopacity{0.000000}%
\pgfsetdash{}{0pt}%
\pgfpathmoveto{\pgfqpoint{0.525000in}{0.450000in}}%
\pgfpathlineto{\pgfqpoint{3.412500in}{0.450000in}}%
\pgfpathlineto{\pgfqpoint{3.412500in}{2.850000in}}%
\pgfpathlineto{\pgfqpoint{0.525000in}{2.850000in}}%
\pgfpathclose%
\pgfusepath{fill}%
\end{pgfscope}%
\begin{pgfscope}%
\pgfsetbuttcap%
\pgfsetroundjoin%
\definecolor{currentfill}{rgb}{0.000000,0.000000,0.000000}%
\pgfsetfillcolor{currentfill}%
\pgfsetlinewidth{0.803000pt}%
\definecolor{currentstroke}{rgb}{0.000000,0.000000,0.000000}%
\pgfsetstrokecolor{currentstroke}%
\pgfsetdash{}{0pt}%
\pgfsys@defobject{currentmarker}{\pgfqpoint{0.000000in}{-0.048611in}}{\pgfqpoint{0.000000in}{0.000000in}}{%
\pgfpathmoveto{\pgfqpoint{0.000000in}{0.000000in}}%
\pgfpathlineto{\pgfqpoint{0.000000in}{-0.048611in}}%
\pgfusepath{stroke,fill}%
}%
\begin{pgfscope}%
\pgfsys@transformshift{0.656250in}{0.450000in}%
\pgfsys@useobject{currentmarker}{}%
\end{pgfscope}%
\end{pgfscope}%
\begin{pgfscope}%
\definecolor{textcolor}{rgb}{0.000000,0.000000,0.000000}%
\pgfsetstrokecolor{textcolor}%
\pgfsetfillcolor{textcolor}%
\pgftext[x=0.656250in,y=0.352778in,,top]{\color{textcolor}\rmfamily\fontsize{10.000000}{12.000000}\selectfont \(\displaystyle {0.0}\)}%
\end{pgfscope}%
\begin{pgfscope}%
\pgfsetbuttcap%
\pgfsetroundjoin%
\definecolor{currentfill}{rgb}{0.000000,0.000000,0.000000}%
\pgfsetfillcolor{currentfill}%
\pgfsetlinewidth{0.803000pt}%
\definecolor{currentstroke}{rgb}{0.000000,0.000000,0.000000}%
\pgfsetstrokecolor{currentstroke}%
\pgfsetdash{}{0pt}%
\pgfsys@defobject{currentmarker}{\pgfqpoint{0.000000in}{-0.048611in}}{\pgfqpoint{0.000000in}{0.000000in}}{%
\pgfpathmoveto{\pgfqpoint{0.000000in}{0.000000in}}%
\pgfpathlineto{\pgfqpoint{0.000000in}{-0.048611in}}%
\pgfusepath{stroke,fill}%
}%
\begin{pgfscope}%
\pgfsys@transformshift{1.181250in}{0.450000in}%
\pgfsys@useobject{currentmarker}{}%
\end{pgfscope}%
\end{pgfscope}%
\begin{pgfscope}%
\definecolor{textcolor}{rgb}{0.000000,0.000000,0.000000}%
\pgfsetstrokecolor{textcolor}%
\pgfsetfillcolor{textcolor}%
\pgftext[x=1.181250in,y=0.352778in,,top]{\color{textcolor}\rmfamily\fontsize{10.000000}{12.000000}\selectfont \(\displaystyle {0.1}\)}%
\end{pgfscope}%
\begin{pgfscope}%
\pgfsetbuttcap%
\pgfsetroundjoin%
\definecolor{currentfill}{rgb}{0.000000,0.000000,0.000000}%
\pgfsetfillcolor{currentfill}%
\pgfsetlinewidth{0.803000pt}%
\definecolor{currentstroke}{rgb}{0.000000,0.000000,0.000000}%
\pgfsetstrokecolor{currentstroke}%
\pgfsetdash{}{0pt}%
\pgfsys@defobject{currentmarker}{\pgfqpoint{0.000000in}{-0.048611in}}{\pgfqpoint{0.000000in}{0.000000in}}{%
\pgfpathmoveto{\pgfqpoint{0.000000in}{0.000000in}}%
\pgfpathlineto{\pgfqpoint{0.000000in}{-0.048611in}}%
\pgfusepath{stroke,fill}%
}%
\begin{pgfscope}%
\pgfsys@transformshift{1.706250in}{0.450000in}%
\pgfsys@useobject{currentmarker}{}%
\end{pgfscope}%
\end{pgfscope}%
\begin{pgfscope}%
\definecolor{textcolor}{rgb}{0.000000,0.000000,0.000000}%
\pgfsetstrokecolor{textcolor}%
\pgfsetfillcolor{textcolor}%
\pgftext[x=1.706250in,y=0.352778in,,top]{\color{textcolor}\rmfamily\fontsize{10.000000}{12.000000}\selectfont \(\displaystyle {0.2}\)}%
\end{pgfscope}%
\begin{pgfscope}%
\pgfsetbuttcap%
\pgfsetroundjoin%
\definecolor{currentfill}{rgb}{0.000000,0.000000,0.000000}%
\pgfsetfillcolor{currentfill}%
\pgfsetlinewidth{0.803000pt}%
\definecolor{currentstroke}{rgb}{0.000000,0.000000,0.000000}%
\pgfsetstrokecolor{currentstroke}%
\pgfsetdash{}{0pt}%
\pgfsys@defobject{currentmarker}{\pgfqpoint{0.000000in}{-0.048611in}}{\pgfqpoint{0.000000in}{0.000000in}}{%
\pgfpathmoveto{\pgfqpoint{0.000000in}{0.000000in}}%
\pgfpathlineto{\pgfqpoint{0.000000in}{-0.048611in}}%
\pgfusepath{stroke,fill}%
}%
\begin{pgfscope}%
\pgfsys@transformshift{2.231250in}{0.450000in}%
\pgfsys@useobject{currentmarker}{}%
\end{pgfscope}%
\end{pgfscope}%
\begin{pgfscope}%
\definecolor{textcolor}{rgb}{0.000000,0.000000,0.000000}%
\pgfsetstrokecolor{textcolor}%
\pgfsetfillcolor{textcolor}%
\pgftext[x=2.231250in,y=0.352778in,,top]{\color{textcolor}\rmfamily\fontsize{10.000000}{12.000000}\selectfont \(\displaystyle {0.3}\)}%
\end{pgfscope}%
\begin{pgfscope}%
\pgfsetbuttcap%
\pgfsetroundjoin%
\definecolor{currentfill}{rgb}{0.000000,0.000000,0.000000}%
\pgfsetfillcolor{currentfill}%
\pgfsetlinewidth{0.803000pt}%
\definecolor{currentstroke}{rgb}{0.000000,0.000000,0.000000}%
\pgfsetstrokecolor{currentstroke}%
\pgfsetdash{}{0pt}%
\pgfsys@defobject{currentmarker}{\pgfqpoint{0.000000in}{-0.048611in}}{\pgfqpoint{0.000000in}{0.000000in}}{%
\pgfpathmoveto{\pgfqpoint{0.000000in}{0.000000in}}%
\pgfpathlineto{\pgfqpoint{0.000000in}{-0.048611in}}%
\pgfusepath{stroke,fill}%
}%
\begin{pgfscope}%
\pgfsys@transformshift{2.756250in}{0.450000in}%
\pgfsys@useobject{currentmarker}{}%
\end{pgfscope}%
\end{pgfscope}%
\begin{pgfscope}%
\definecolor{textcolor}{rgb}{0.000000,0.000000,0.000000}%
\pgfsetstrokecolor{textcolor}%
\pgfsetfillcolor{textcolor}%
\pgftext[x=2.756250in,y=0.352778in,,top]{\color{textcolor}\rmfamily\fontsize{10.000000}{12.000000}\selectfont \(\displaystyle {0.4}\)}%
\end{pgfscope}%
\begin{pgfscope}%
\pgfsetbuttcap%
\pgfsetroundjoin%
\definecolor{currentfill}{rgb}{0.000000,0.000000,0.000000}%
\pgfsetfillcolor{currentfill}%
\pgfsetlinewidth{0.803000pt}%
\definecolor{currentstroke}{rgb}{0.000000,0.000000,0.000000}%
\pgfsetstrokecolor{currentstroke}%
\pgfsetdash{}{0pt}%
\pgfsys@defobject{currentmarker}{\pgfqpoint{0.000000in}{-0.048611in}}{\pgfqpoint{0.000000in}{0.000000in}}{%
\pgfpathmoveto{\pgfqpoint{0.000000in}{0.000000in}}%
\pgfpathlineto{\pgfqpoint{0.000000in}{-0.048611in}}%
\pgfusepath{stroke,fill}%
}%
\begin{pgfscope}%
\pgfsys@transformshift{3.281250in}{0.450000in}%
\pgfsys@useobject{currentmarker}{}%
\end{pgfscope}%
\end{pgfscope}%
\begin{pgfscope}%
\definecolor{textcolor}{rgb}{0.000000,0.000000,0.000000}%
\pgfsetstrokecolor{textcolor}%
\pgfsetfillcolor{textcolor}%
\pgftext[x=3.281250in,y=0.352778in,,top]{\color{textcolor}\rmfamily\fontsize{10.000000}{12.000000}\selectfont \(\displaystyle {0.5}\)}%
\end{pgfscope}%
\begin{pgfscope}%
\definecolor{textcolor}{rgb}{0.000000,0.000000,0.000000}%
\pgfsetstrokecolor{textcolor}%
\pgfsetfillcolor{textcolor}%
\pgftext[x=1.968750in,y=0.162809in,,top]{\color{textcolor}\rmfamily\fontsize{10.000000}{12.000000}\selectfont \(\displaystyle t_{\mathrm{blend}}/\Delta t_{\mathrm{obs}}\)}%
\end{pgfscope}%
\begin{pgfscope}%
\pgfsetbuttcap%
\pgfsetroundjoin%
\definecolor{currentfill}{rgb}{0.000000,0.000000,0.000000}%
\pgfsetfillcolor{currentfill}%
\pgfsetlinewidth{0.803000pt}%
\definecolor{currentstroke}{rgb}{0.000000,0.000000,0.000000}%
\pgfsetstrokecolor{currentstroke}%
\pgfsetdash{}{0pt}%
\pgfsys@defobject{currentmarker}{\pgfqpoint{-0.048611in}{0.000000in}}{\pgfqpoint{0.000000in}{0.000000in}}{%
\pgfpathmoveto{\pgfqpoint{0.000000in}{0.000000in}}%
\pgfpathlineto{\pgfqpoint{-0.048611in}{0.000000in}}%
\pgfusepath{stroke,fill}%
}%
\begin{pgfscope}%
\pgfsys@transformshift{0.525000in}{0.450000in}%
\pgfsys@useobject{currentmarker}{}%
\end{pgfscope}%
\end{pgfscope}%
\begin{pgfscope}%
\definecolor{textcolor}{rgb}{0.000000,0.000000,0.000000}%
\pgfsetstrokecolor{textcolor}%
\pgfsetfillcolor{textcolor}%
\pgftext[x=0.180863in, y=0.397238in, left, base]{\color{textcolor}\rmfamily\fontsize{10.000000}{12.000000}\selectfont \(\displaystyle {0.00}\)}%
\end{pgfscope}%
\begin{pgfscope}%
\pgfsetbuttcap%
\pgfsetroundjoin%
\definecolor{currentfill}{rgb}{0.000000,0.000000,0.000000}%
\pgfsetfillcolor{currentfill}%
\pgfsetlinewidth{0.803000pt}%
\definecolor{currentstroke}{rgb}{0.000000,0.000000,0.000000}%
\pgfsetstrokecolor{currentstroke}%
\pgfsetdash{}{0pt}%
\pgfsys@defobject{currentmarker}{\pgfqpoint{-0.048611in}{0.000000in}}{\pgfqpoint{0.000000in}{0.000000in}}{%
\pgfpathmoveto{\pgfqpoint{0.000000in}{0.000000in}}%
\pgfpathlineto{\pgfqpoint{-0.048611in}{0.000000in}}%
\pgfusepath{stroke,fill}%
}%
\begin{pgfscope}%
\pgfsys@transformshift{0.525000in}{0.750000in}%
\pgfsys@useobject{currentmarker}{}%
\end{pgfscope}%
\end{pgfscope}%
\begin{pgfscope}%
\definecolor{textcolor}{rgb}{0.000000,0.000000,0.000000}%
\pgfsetstrokecolor{textcolor}%
\pgfsetfillcolor{textcolor}%
\pgftext[x=0.180863in, y=0.697238in, left, base]{\color{textcolor}\rmfamily\fontsize{10.000000}{12.000000}\selectfont \(\displaystyle {0.25}\)}%
\end{pgfscope}%
\begin{pgfscope}%
\pgfsetbuttcap%
\pgfsetroundjoin%
\definecolor{currentfill}{rgb}{0.000000,0.000000,0.000000}%
\pgfsetfillcolor{currentfill}%
\pgfsetlinewidth{0.803000pt}%
\definecolor{currentstroke}{rgb}{0.000000,0.000000,0.000000}%
\pgfsetstrokecolor{currentstroke}%
\pgfsetdash{}{0pt}%
\pgfsys@defobject{currentmarker}{\pgfqpoint{-0.048611in}{0.000000in}}{\pgfqpoint{0.000000in}{0.000000in}}{%
\pgfpathmoveto{\pgfqpoint{0.000000in}{0.000000in}}%
\pgfpathlineto{\pgfqpoint{-0.048611in}{0.000000in}}%
\pgfusepath{stroke,fill}%
}%
\begin{pgfscope}%
\pgfsys@transformshift{0.525000in}{1.050000in}%
\pgfsys@useobject{currentmarker}{}%
\end{pgfscope}%
\end{pgfscope}%
\begin{pgfscope}%
\definecolor{textcolor}{rgb}{0.000000,0.000000,0.000000}%
\pgfsetstrokecolor{textcolor}%
\pgfsetfillcolor{textcolor}%
\pgftext[x=0.180863in, y=0.997238in, left, base]{\color{textcolor}\rmfamily\fontsize{10.000000}{12.000000}\selectfont \(\displaystyle {0.50}\)}%
\end{pgfscope}%
\begin{pgfscope}%
\pgfsetbuttcap%
\pgfsetroundjoin%
\definecolor{currentfill}{rgb}{0.000000,0.000000,0.000000}%
\pgfsetfillcolor{currentfill}%
\pgfsetlinewidth{0.803000pt}%
\definecolor{currentstroke}{rgb}{0.000000,0.000000,0.000000}%
\pgfsetstrokecolor{currentstroke}%
\pgfsetdash{}{0pt}%
\pgfsys@defobject{currentmarker}{\pgfqpoint{-0.048611in}{0.000000in}}{\pgfqpoint{0.000000in}{0.000000in}}{%
\pgfpathmoveto{\pgfqpoint{0.000000in}{0.000000in}}%
\pgfpathlineto{\pgfqpoint{-0.048611in}{0.000000in}}%
\pgfusepath{stroke,fill}%
}%
\begin{pgfscope}%
\pgfsys@transformshift{0.525000in}{1.350000in}%
\pgfsys@useobject{currentmarker}{}%
\end{pgfscope}%
\end{pgfscope}%
\begin{pgfscope}%
\definecolor{textcolor}{rgb}{0.000000,0.000000,0.000000}%
\pgfsetstrokecolor{textcolor}%
\pgfsetfillcolor{textcolor}%
\pgftext[x=0.180863in, y=1.297238in, left, base]{\color{textcolor}\rmfamily\fontsize{10.000000}{12.000000}\selectfont \(\displaystyle {0.75}\)}%
\end{pgfscope}%
\begin{pgfscope}%
\pgfsetbuttcap%
\pgfsetroundjoin%
\definecolor{currentfill}{rgb}{0.000000,0.000000,0.000000}%
\pgfsetfillcolor{currentfill}%
\pgfsetlinewidth{0.803000pt}%
\definecolor{currentstroke}{rgb}{0.000000,0.000000,0.000000}%
\pgfsetstrokecolor{currentstroke}%
\pgfsetdash{}{0pt}%
\pgfsys@defobject{currentmarker}{\pgfqpoint{-0.048611in}{0.000000in}}{\pgfqpoint{0.000000in}{0.000000in}}{%
\pgfpathmoveto{\pgfqpoint{0.000000in}{0.000000in}}%
\pgfpathlineto{\pgfqpoint{-0.048611in}{0.000000in}}%
\pgfusepath{stroke,fill}%
}%
\begin{pgfscope}%
\pgfsys@transformshift{0.525000in}{1.650000in}%
\pgfsys@useobject{currentmarker}{}%
\end{pgfscope}%
\end{pgfscope}%
\begin{pgfscope}%
\definecolor{textcolor}{rgb}{0.000000,0.000000,0.000000}%
\pgfsetstrokecolor{textcolor}%
\pgfsetfillcolor{textcolor}%
\pgftext[x=0.180863in, y=1.597238in, left, base]{\color{textcolor}\rmfamily\fontsize{10.000000}{12.000000}\selectfont \(\displaystyle {1.00}\)}%
\end{pgfscope}%
\begin{pgfscope}%
\pgfsetbuttcap%
\pgfsetroundjoin%
\definecolor{currentfill}{rgb}{0.000000,0.000000,0.000000}%
\pgfsetfillcolor{currentfill}%
\pgfsetlinewidth{0.803000pt}%
\definecolor{currentstroke}{rgb}{0.000000,0.000000,0.000000}%
\pgfsetstrokecolor{currentstroke}%
\pgfsetdash{}{0pt}%
\pgfsys@defobject{currentmarker}{\pgfqpoint{-0.048611in}{0.000000in}}{\pgfqpoint{0.000000in}{0.000000in}}{%
\pgfpathmoveto{\pgfqpoint{0.000000in}{0.000000in}}%
\pgfpathlineto{\pgfqpoint{-0.048611in}{0.000000in}}%
\pgfusepath{stroke,fill}%
}%
\begin{pgfscope}%
\pgfsys@transformshift{0.525000in}{1.950000in}%
\pgfsys@useobject{currentmarker}{}%
\end{pgfscope}%
\end{pgfscope}%
\begin{pgfscope}%
\definecolor{textcolor}{rgb}{0.000000,0.000000,0.000000}%
\pgfsetstrokecolor{textcolor}%
\pgfsetfillcolor{textcolor}%
\pgftext[x=0.180863in, y=1.897238in, left, base]{\color{textcolor}\rmfamily\fontsize{10.000000}{12.000000}\selectfont \(\displaystyle {1.25}\)}%
\end{pgfscope}%
\begin{pgfscope}%
\pgfsetbuttcap%
\pgfsetroundjoin%
\definecolor{currentfill}{rgb}{0.000000,0.000000,0.000000}%
\pgfsetfillcolor{currentfill}%
\pgfsetlinewidth{0.803000pt}%
\definecolor{currentstroke}{rgb}{0.000000,0.000000,0.000000}%
\pgfsetstrokecolor{currentstroke}%
\pgfsetdash{}{0pt}%
\pgfsys@defobject{currentmarker}{\pgfqpoint{-0.048611in}{0.000000in}}{\pgfqpoint{0.000000in}{0.000000in}}{%
\pgfpathmoveto{\pgfqpoint{0.000000in}{0.000000in}}%
\pgfpathlineto{\pgfqpoint{-0.048611in}{0.000000in}}%
\pgfusepath{stroke,fill}%
}%
\begin{pgfscope}%
\pgfsys@transformshift{0.525000in}{2.250000in}%
\pgfsys@useobject{currentmarker}{}%
\end{pgfscope}%
\end{pgfscope}%
\begin{pgfscope}%
\definecolor{textcolor}{rgb}{0.000000,0.000000,0.000000}%
\pgfsetstrokecolor{textcolor}%
\pgfsetfillcolor{textcolor}%
\pgftext[x=0.180863in, y=2.197238in, left, base]{\color{textcolor}\rmfamily\fontsize{10.000000}{12.000000}\selectfont \(\displaystyle {1.50}\)}%
\end{pgfscope}%
\begin{pgfscope}%
\pgfsetbuttcap%
\pgfsetroundjoin%
\definecolor{currentfill}{rgb}{0.000000,0.000000,0.000000}%
\pgfsetfillcolor{currentfill}%
\pgfsetlinewidth{0.803000pt}%
\definecolor{currentstroke}{rgb}{0.000000,0.000000,0.000000}%
\pgfsetstrokecolor{currentstroke}%
\pgfsetdash{}{0pt}%
\pgfsys@defobject{currentmarker}{\pgfqpoint{-0.048611in}{0.000000in}}{\pgfqpoint{0.000000in}{0.000000in}}{%
\pgfpathmoveto{\pgfqpoint{0.000000in}{0.000000in}}%
\pgfpathlineto{\pgfqpoint{-0.048611in}{0.000000in}}%
\pgfusepath{stroke,fill}%
}%
\begin{pgfscope}%
\pgfsys@transformshift{0.525000in}{2.550000in}%
\pgfsys@useobject{currentmarker}{}%
\end{pgfscope}%
\end{pgfscope}%
\begin{pgfscope}%
\definecolor{textcolor}{rgb}{0.000000,0.000000,0.000000}%
\pgfsetstrokecolor{textcolor}%
\pgfsetfillcolor{textcolor}%
\pgftext[x=0.180863in, y=2.497238in, left, base]{\color{textcolor}\rmfamily\fontsize{10.000000}{12.000000}\selectfont \(\displaystyle {1.75}\)}%
\end{pgfscope}%
\begin{pgfscope}%
\pgfsetbuttcap%
\pgfsetroundjoin%
\definecolor{currentfill}{rgb}{0.000000,0.000000,0.000000}%
\pgfsetfillcolor{currentfill}%
\pgfsetlinewidth{0.803000pt}%
\definecolor{currentstroke}{rgb}{0.000000,0.000000,0.000000}%
\pgfsetstrokecolor{currentstroke}%
\pgfsetdash{}{0pt}%
\pgfsys@defobject{currentmarker}{\pgfqpoint{-0.048611in}{0.000000in}}{\pgfqpoint{0.000000in}{0.000000in}}{%
\pgfpathmoveto{\pgfqpoint{0.000000in}{0.000000in}}%
\pgfpathlineto{\pgfqpoint{-0.048611in}{0.000000in}}%
\pgfusepath{stroke,fill}%
}%
\begin{pgfscope}%
\pgfsys@transformshift{0.525000in}{2.850000in}%
\pgfsys@useobject{currentmarker}{}%
\end{pgfscope}%
\end{pgfscope}%
\begin{pgfscope}%
\definecolor{textcolor}{rgb}{0.000000,0.000000,0.000000}%
\pgfsetstrokecolor{textcolor}%
\pgfsetfillcolor{textcolor}%
\pgftext[x=0.180863in, y=2.797238in, left, base]{\color{textcolor}\rmfamily\fontsize{10.000000}{12.000000}\selectfont \(\displaystyle {2.00}\)}%
\end{pgfscope}%
\begin{pgfscope}%
\definecolor{textcolor}{rgb}{0.000000,0.000000,0.000000}%
\pgfsetstrokecolor{textcolor}%
\pgfsetfillcolor{textcolor}%
\pgftext[x=0.125308in,y=1.650000in,,bottom,rotate=90.000000]{\color{textcolor}\rmfamily\fontsize{10.000000}{12.000000}\selectfont \(\displaystyle \mathrm{TRMSE}(q)\)}%
\end{pgfscope}%
\begin{pgfscope}%
\pgfpathrectangle{\pgfqpoint{0.525000in}{0.450000in}}{\pgfqpoint{2.887500in}{2.400000in}}%
\pgfusepath{clip}%
\pgfsetrectcap%
\pgfsetroundjoin%
\pgfsetlinewidth{1.003750pt}%
\definecolor{currentstroke}{rgb}{1.000000,0.600000,0.000000}%
\pgfsetstrokecolor{currentstroke}%
\pgfsetdash{}{0pt}%
\pgfpathmoveto{\pgfqpoint{0.656250in}{2.694752in}}%
\pgfpathlineto{\pgfqpoint{3.281250in}{2.694752in}}%
\pgfusepath{stroke}%
\end{pgfscope}%
\begin{pgfscope}%
\pgfpathrectangle{\pgfqpoint{0.525000in}{0.450000in}}{\pgfqpoint{2.887500in}{2.400000in}}%
\pgfusepath{clip}%
\pgfsetrectcap%
\pgfsetroundjoin%
\pgfsetlinewidth{1.003750pt}%
\definecolor{currentstroke}{rgb}{0.000000,0.200000,0.400000}%
\pgfsetstrokecolor{currentstroke}%
\pgfsetdash{}{0pt}%
\pgfpathmoveto{\pgfqpoint{0.656250in}{2.690759in}}%
\pgfpathlineto{\pgfqpoint{0.984375in}{1.983325in}}%
\pgfpathlineto{\pgfqpoint{1.312500in}{1.519134in}}%
\pgfpathlineto{\pgfqpoint{1.640625in}{1.227934in}}%
\pgfpathlineto{\pgfqpoint{1.968750in}{0.805395in}}%
\pgfpathlineto{\pgfqpoint{2.296875in}{0.916976in}}%
\pgfpathlineto{\pgfqpoint{2.625000in}{0.922930in}}%
\pgfpathlineto{\pgfqpoint{2.953125in}{0.842513in}}%
\pgfpathlineto{\pgfqpoint{3.281250in}{0.825388in}}%
\pgfusepath{stroke}%
\end{pgfscope}%
\begin{pgfscope}%
\pgfpathrectangle{\pgfqpoint{0.525000in}{0.450000in}}{\pgfqpoint{2.887500in}{2.400000in}}%
\pgfusepath{clip}%
\pgfsetbuttcap%
\pgfsetroundjoin%
\definecolor{currentfill}{rgb}{0.000000,0.200000,0.400000}%
\pgfsetfillcolor{currentfill}%
\pgfsetlinewidth{1.003750pt}%
\definecolor{currentstroke}{rgb}{0.000000,0.200000,0.400000}%
\pgfsetstrokecolor{currentstroke}%
\pgfsetdash{}{0pt}%
\pgfsys@defobject{currentmarker}{\pgfqpoint{-0.027778in}{-0.027778in}}{\pgfqpoint{0.027778in}{0.027778in}}{%
\pgfpathmoveto{\pgfqpoint{-0.027778in}{-0.027778in}}%
\pgfpathlineto{\pgfqpoint{0.027778in}{0.027778in}}%
\pgfpathmoveto{\pgfqpoint{-0.027778in}{0.027778in}}%
\pgfpathlineto{\pgfqpoint{0.027778in}{-0.027778in}}%
\pgfusepath{stroke,fill}%
}%
\begin{pgfscope}%
\pgfsys@transformshift{0.656250in}{2.690759in}%
\pgfsys@useobject{currentmarker}{}%
\end{pgfscope}%
\begin{pgfscope}%
\pgfsys@transformshift{0.984375in}{1.983325in}%
\pgfsys@useobject{currentmarker}{}%
\end{pgfscope}%
\begin{pgfscope}%
\pgfsys@transformshift{1.312500in}{1.519134in}%
\pgfsys@useobject{currentmarker}{}%
\end{pgfscope}%
\begin{pgfscope}%
\pgfsys@transformshift{1.640625in}{1.227934in}%
\pgfsys@useobject{currentmarker}{}%
\end{pgfscope}%
\begin{pgfscope}%
\pgfsys@transformshift{1.968750in}{0.805395in}%
\pgfsys@useobject{currentmarker}{}%
\end{pgfscope}%
\begin{pgfscope}%
\pgfsys@transformshift{2.296875in}{0.916976in}%
\pgfsys@useobject{currentmarker}{}%
\end{pgfscope}%
\begin{pgfscope}%
\pgfsys@transformshift{2.625000in}{0.922930in}%
\pgfsys@useobject{currentmarker}{}%
\end{pgfscope}%
\begin{pgfscope}%
\pgfsys@transformshift{2.953125in}{0.842513in}%
\pgfsys@useobject{currentmarker}{}%
\end{pgfscope}%
\begin{pgfscope}%
\pgfsys@transformshift{3.281250in}{0.825388in}%
\pgfsys@useobject{currentmarker}{}%
\end{pgfscope}%
\end{pgfscope}%
\begin{pgfscope}%
\pgfsetrectcap%
\pgfsetmiterjoin%
\pgfsetlinewidth{0.803000pt}%
\definecolor{currentstroke}{rgb}{0.000000,0.000000,0.000000}%
\pgfsetstrokecolor{currentstroke}%
\pgfsetdash{}{0pt}%
\pgfpathmoveto{\pgfqpoint{0.525000in}{0.450000in}}%
\pgfpathlineto{\pgfqpoint{0.525000in}{2.850000in}}%
\pgfusepath{stroke}%
\end{pgfscope}%
\begin{pgfscope}%
\pgfsetrectcap%
\pgfsetmiterjoin%
\pgfsetlinewidth{0.803000pt}%
\definecolor{currentstroke}{rgb}{0.000000,0.000000,0.000000}%
\pgfsetstrokecolor{currentstroke}%
\pgfsetdash{}{0pt}%
\pgfpathmoveto{\pgfqpoint{3.412500in}{0.450000in}}%
\pgfpathlineto{\pgfqpoint{3.412500in}{2.850000in}}%
\pgfusepath{stroke}%
\end{pgfscope}%
\begin{pgfscope}%
\pgfsetrectcap%
\pgfsetmiterjoin%
\pgfsetlinewidth{0.803000pt}%
\definecolor{currentstroke}{rgb}{0.000000,0.000000,0.000000}%
\pgfsetstrokecolor{currentstroke}%
\pgfsetdash{}{0pt}%
\pgfpathmoveto{\pgfqpoint{0.525000in}{0.450000in}}%
\pgfpathlineto{\pgfqpoint{3.412500in}{0.450000in}}%
\pgfusepath{stroke}%
\end{pgfscope}%
\begin{pgfscope}%
\pgfsetrectcap%
\pgfsetmiterjoin%
\pgfsetlinewidth{0.803000pt}%
\definecolor{currentstroke}{rgb}{0.000000,0.000000,0.000000}%
\pgfsetstrokecolor{currentstroke}%
\pgfsetdash{}{0pt}%
\pgfpathmoveto{\pgfqpoint{0.525000in}{2.850000in}}%
\pgfpathlineto{\pgfqpoint{3.412500in}{2.850000in}}%
\pgfusepath{stroke}%
\end{pgfscope}%
\begin{pgfscope}%
\pgfsetbuttcap%
\pgfsetmiterjoin%
\definecolor{currentfill}{rgb}{1.000000,1.000000,1.000000}%
\pgfsetfillcolor{currentfill}%
\pgfsetfillopacity{0.800000}%
\pgfsetlinewidth{1.003750pt}%
\definecolor{currentstroke}{rgb}{0.800000,0.800000,0.800000}%
\pgfsetstrokecolor{currentstroke}%
\pgfsetstrokeopacity{0.800000}%
\pgfsetdash{}{0pt}%
\pgfpathmoveto{\pgfqpoint{2.613704in}{1.514596in}}%
\pgfpathlineto{\pgfqpoint{3.354167in}{1.514596in}}%
\pgfpathquadraticcurveto{\pgfqpoint{3.370833in}{1.514596in}}{\pgfqpoint{3.370833in}{1.531262in}}%
\pgfpathlineto{\pgfqpoint{3.370833in}{1.768738in}}%
\pgfpathquadraticcurveto{\pgfqpoint{3.370833in}{1.785404in}}{\pgfqpoint{3.354167in}{1.785404in}}%
\pgfpathlineto{\pgfqpoint{2.613704in}{1.785404in}}%
\pgfpathquadraticcurveto{\pgfqpoint{2.597038in}{1.785404in}}{\pgfqpoint{2.597038in}{1.768738in}}%
\pgfpathlineto{\pgfqpoint{2.597038in}{1.531262in}}%
\pgfpathquadraticcurveto{\pgfqpoint{2.597038in}{1.514596in}}{\pgfqpoint{2.613704in}{1.514596in}}%
\pgfpathclose%
\pgfusepath{stroke,fill}%
\end{pgfscope}%
\begin{pgfscope}%
\pgfsetrectcap%
\pgfsetroundjoin%
\pgfsetlinewidth{1.003750pt}%
\definecolor{currentstroke}{rgb}{1.000000,0.600000,0.000000}%
\pgfsetstrokecolor{currentstroke}%
\pgfsetdash{}{0pt}%
\pgfpathmoveto{\pgfqpoint{2.630371in}{1.717924in}}%
\pgfpathlineto{\pgfqpoint{2.797038in}{1.717924in}}%
\pgfusepath{stroke}%
\end{pgfscope}%
\begin{pgfscope}%
\definecolor{textcolor}{rgb}{0.000000,0.000000,0.000000}%
\pgfsetstrokecolor{textcolor}%
\pgfsetfillcolor{textcolor}%
\pgftext[x=2.863704in,y=1.688757in,left,base]{\color{textcolor}\rmfamily\fontsize{6.000000}{7.200000}\selectfont Plain EnKF}%
\end{pgfscope}%
\begin{pgfscope}%
\pgfsetrectcap%
\pgfsetroundjoin%
\pgfsetlinewidth{1.003750pt}%
\definecolor{currentstroke}{rgb}{0.000000,0.200000,0.400000}%
\pgfsetstrokecolor{currentstroke}%
\pgfsetdash{}{0pt}%
\pgfpathmoveto{\pgfqpoint{2.630371in}{1.595610in}}%
\pgfpathlineto{\pgfqpoint{2.797038in}{1.595610in}}%
\pgfusepath{stroke}%
\end{pgfscope}%
\begin{pgfscope}%
\pgfsetbuttcap%
\pgfsetroundjoin%
\definecolor{currentfill}{rgb}{0.000000,0.200000,0.400000}%
\pgfsetfillcolor{currentfill}%
\pgfsetlinewidth{1.003750pt}%
\definecolor{currentstroke}{rgb}{0.000000,0.200000,0.400000}%
\pgfsetstrokecolor{currentstroke}%
\pgfsetdash{}{0pt}%
\pgfsys@defobject{currentmarker}{\pgfqpoint{-0.027778in}{-0.027778in}}{\pgfqpoint{0.027778in}{0.027778in}}{%
\pgfpathmoveto{\pgfqpoint{-0.027778in}{-0.027778in}}%
\pgfpathlineto{\pgfqpoint{0.027778in}{0.027778in}}%
\pgfpathmoveto{\pgfqpoint{-0.027778in}{0.027778in}}%
\pgfpathlineto{\pgfqpoint{0.027778in}{-0.027778in}}%
\pgfusepath{stroke,fill}%
}%
\begin{pgfscope}%
\pgfsys@transformshift{2.713704in}{1.595610in}%
\pgfsys@useobject{currentmarker}{}%
\end{pgfscope}%
\end{pgfscope}%
\begin{pgfscope}%
\definecolor{textcolor}{rgb}{0.000000,0.000000,0.000000}%
\pgfsetstrokecolor{textcolor}%
\pgfsetfillcolor{textcolor}%
\pgftext[x=2.863704in,y=1.566443in,left,base]{\color{textcolor}\rmfamily\fontsize{6.000000}{7.200000}\selectfont Blending}%
\end{pgfscope}%
\end{pgfpicture}%
\makeatother%
\endgroup%

%% file: image/results_vel_blending.pgf
%% Creator: Matplotlib, PGF backend
%%
%% To include the figure in your LaTeX document, write
%%   \input{<filename>.pgf}
%%
%% Make sure the required packages are loaded in your preamble
%%   \usepackage{pgf}
%%
%% and, on pdftex
%%   \usepackage[utf8]{inputenc}\DeclareUnicodeCharacter{2212}{-}
%%
%% or, on luatex and xetex
%%   \usepackage{unicode-math}
%%
%% Figures using additional raster images can only be included by \input if
%% they are in the same directory as the main LaTeX file. For loading figures
%% from other directories you can use the `import` package
%% and then include the figures with
%%   \import{<path to file>}{<filename>.pgf}
%%
%% Matplotlib used the following preamble
%%   \usepackage{fontspec}
%%   \setmainfont{DejaVuSerif.ttf}[Path=/home/gottfried/.local/lib/python3.8/site-packages/matplotlib/mpl-data/fonts/ttf/]
%%   \setsansfont{DejaVuSans.ttf}[Path=/home/gottfried/.local/lib/python3.8/site-packages/matplotlib/mpl-data/fonts/ttf/]
%%   \setmonofont{DejaVuSansMono.ttf}[Path=/home/gottfried/.local/lib/python3.8/site-packages/matplotlib/mpl-data/fonts/ttf/]
%%
\begingroup%
\makeatletter%
\begin{pgfpicture}%
\pgfpathrectangle{\pgfpointorigin}{\pgfqpoint{3.500000in}{3.000000in}}%
\pgfusepath{use as bounding box, clip}%
\begin{pgfscope}%
\pgfsetbuttcap%
\pgfsetmiterjoin%
\definecolor{currentfill}{rgb}{1.000000,1.000000,1.000000}%
\pgfsetfillcolor{currentfill}%
\pgfsetlinewidth{0.000000pt}%
\definecolor{currentstroke}{rgb}{1.000000,1.000000,1.000000}%
\pgfsetstrokecolor{currentstroke}%
\pgfsetdash{}{0pt}%
\pgfpathmoveto{\pgfqpoint{0.000000in}{0.000000in}}%
\pgfpathlineto{\pgfqpoint{3.500000in}{0.000000in}}%
\pgfpathlineto{\pgfqpoint{3.500000in}{3.000000in}}%
\pgfpathlineto{\pgfqpoint{0.000000in}{3.000000in}}%
\pgfpathclose%
\pgfusepath{fill}%
\end{pgfscope}%
\begin{pgfscope}%
\pgfsetbuttcap%
\pgfsetmiterjoin%
\definecolor{currentfill}{rgb}{1.000000,1.000000,1.000000}%
\pgfsetfillcolor{currentfill}%
\pgfsetlinewidth{0.000000pt}%
\definecolor{currentstroke}{rgb}{0.000000,0.000000,0.000000}%
\pgfsetstrokecolor{currentstroke}%
\pgfsetstrokeopacity{0.000000}%
\pgfsetdash{}{0pt}%
\pgfpathmoveto{\pgfqpoint{0.525000in}{0.450000in}}%
\pgfpathlineto{\pgfqpoint{3.412500in}{0.450000in}}%
\pgfpathlineto{\pgfqpoint{3.412500in}{2.850000in}}%
\pgfpathlineto{\pgfqpoint{0.525000in}{2.850000in}}%
\pgfpathclose%
\pgfusepath{fill}%
\end{pgfscope}%
\begin{pgfscope}%
\pgfsetbuttcap%
\pgfsetroundjoin%
\definecolor{currentfill}{rgb}{0.000000,0.000000,0.000000}%
\pgfsetfillcolor{currentfill}%
\pgfsetlinewidth{0.803000pt}%
\definecolor{currentstroke}{rgb}{0.000000,0.000000,0.000000}%
\pgfsetstrokecolor{currentstroke}%
\pgfsetdash{}{0pt}%
\pgfsys@defobject{currentmarker}{\pgfqpoint{0.000000in}{-0.048611in}}{\pgfqpoint{0.000000in}{0.000000in}}{%
\pgfpathmoveto{\pgfqpoint{0.000000in}{0.000000in}}%
\pgfpathlineto{\pgfqpoint{0.000000in}{-0.048611in}}%
\pgfusepath{stroke,fill}%
}%
\begin{pgfscope}%
\pgfsys@transformshift{0.656250in}{0.450000in}%
\pgfsys@useobject{currentmarker}{}%
\end{pgfscope}%
\end{pgfscope}%
\begin{pgfscope}%
\definecolor{textcolor}{rgb}{0.000000,0.000000,0.000000}%
\pgfsetstrokecolor{textcolor}%
\pgfsetfillcolor{textcolor}%
\pgftext[x=0.656250in,y=0.352778in,,top]{\color{textcolor}\rmfamily\fontsize{10.000000}{12.000000}\selectfont \(\displaystyle {0.0}\)}%
\end{pgfscope}%
\begin{pgfscope}%
\pgfsetbuttcap%
\pgfsetroundjoin%
\definecolor{currentfill}{rgb}{0.000000,0.000000,0.000000}%
\pgfsetfillcolor{currentfill}%
\pgfsetlinewidth{0.803000pt}%
\definecolor{currentstroke}{rgb}{0.000000,0.000000,0.000000}%
\pgfsetstrokecolor{currentstroke}%
\pgfsetdash{}{0pt}%
\pgfsys@defobject{currentmarker}{\pgfqpoint{0.000000in}{-0.048611in}}{\pgfqpoint{0.000000in}{0.000000in}}{%
\pgfpathmoveto{\pgfqpoint{0.000000in}{0.000000in}}%
\pgfpathlineto{\pgfqpoint{0.000000in}{-0.048611in}}%
\pgfusepath{stroke,fill}%
}%
\begin{pgfscope}%
\pgfsys@transformshift{1.181250in}{0.450000in}%
\pgfsys@useobject{currentmarker}{}%
\end{pgfscope}%
\end{pgfscope}%
\begin{pgfscope}%
\definecolor{textcolor}{rgb}{0.000000,0.000000,0.000000}%
\pgfsetstrokecolor{textcolor}%
\pgfsetfillcolor{textcolor}%
\pgftext[x=1.181250in,y=0.352778in,,top]{\color{textcolor}\rmfamily\fontsize{10.000000}{12.000000}\selectfont \(\displaystyle {0.1}\)}%
\end{pgfscope}%
\begin{pgfscope}%
\pgfsetbuttcap%
\pgfsetroundjoin%
\definecolor{currentfill}{rgb}{0.000000,0.000000,0.000000}%
\pgfsetfillcolor{currentfill}%
\pgfsetlinewidth{0.803000pt}%
\definecolor{currentstroke}{rgb}{0.000000,0.000000,0.000000}%
\pgfsetstrokecolor{currentstroke}%
\pgfsetdash{}{0pt}%
\pgfsys@defobject{currentmarker}{\pgfqpoint{0.000000in}{-0.048611in}}{\pgfqpoint{0.000000in}{0.000000in}}{%
\pgfpathmoveto{\pgfqpoint{0.000000in}{0.000000in}}%
\pgfpathlineto{\pgfqpoint{0.000000in}{-0.048611in}}%
\pgfusepath{stroke,fill}%
}%
\begin{pgfscope}%
\pgfsys@transformshift{1.706250in}{0.450000in}%
\pgfsys@useobject{currentmarker}{}%
\end{pgfscope}%
\end{pgfscope}%
\begin{pgfscope}%
\definecolor{textcolor}{rgb}{0.000000,0.000000,0.000000}%
\pgfsetstrokecolor{textcolor}%
\pgfsetfillcolor{textcolor}%
\pgftext[x=1.706250in,y=0.352778in,,top]{\color{textcolor}\rmfamily\fontsize{10.000000}{12.000000}\selectfont \(\displaystyle {0.2}\)}%
\end{pgfscope}%
\begin{pgfscope}%
\pgfsetbuttcap%
\pgfsetroundjoin%
\definecolor{currentfill}{rgb}{0.000000,0.000000,0.000000}%
\pgfsetfillcolor{currentfill}%
\pgfsetlinewidth{0.803000pt}%
\definecolor{currentstroke}{rgb}{0.000000,0.000000,0.000000}%
\pgfsetstrokecolor{currentstroke}%
\pgfsetdash{}{0pt}%
\pgfsys@defobject{currentmarker}{\pgfqpoint{0.000000in}{-0.048611in}}{\pgfqpoint{0.000000in}{0.000000in}}{%
\pgfpathmoveto{\pgfqpoint{0.000000in}{0.000000in}}%
\pgfpathlineto{\pgfqpoint{0.000000in}{-0.048611in}}%
\pgfusepath{stroke,fill}%
}%
\begin{pgfscope}%
\pgfsys@transformshift{2.231250in}{0.450000in}%
\pgfsys@useobject{currentmarker}{}%
\end{pgfscope}%
\end{pgfscope}%
\begin{pgfscope}%
\definecolor{textcolor}{rgb}{0.000000,0.000000,0.000000}%
\pgfsetstrokecolor{textcolor}%
\pgfsetfillcolor{textcolor}%
\pgftext[x=2.231250in,y=0.352778in,,top]{\color{textcolor}\rmfamily\fontsize{10.000000}{12.000000}\selectfont \(\displaystyle {0.3}\)}%
\end{pgfscope}%
\begin{pgfscope}%
\pgfsetbuttcap%
\pgfsetroundjoin%
\definecolor{currentfill}{rgb}{0.000000,0.000000,0.000000}%
\pgfsetfillcolor{currentfill}%
\pgfsetlinewidth{0.803000pt}%
\definecolor{currentstroke}{rgb}{0.000000,0.000000,0.000000}%
\pgfsetstrokecolor{currentstroke}%
\pgfsetdash{}{0pt}%
\pgfsys@defobject{currentmarker}{\pgfqpoint{0.000000in}{-0.048611in}}{\pgfqpoint{0.000000in}{0.000000in}}{%
\pgfpathmoveto{\pgfqpoint{0.000000in}{0.000000in}}%
\pgfpathlineto{\pgfqpoint{0.000000in}{-0.048611in}}%
\pgfusepath{stroke,fill}%
}%
\begin{pgfscope}%
\pgfsys@transformshift{2.756250in}{0.450000in}%
\pgfsys@useobject{currentmarker}{}%
\end{pgfscope}%
\end{pgfscope}%
\begin{pgfscope}%
\definecolor{textcolor}{rgb}{0.000000,0.000000,0.000000}%
\pgfsetstrokecolor{textcolor}%
\pgfsetfillcolor{textcolor}%
\pgftext[x=2.756250in,y=0.352778in,,top]{\color{textcolor}\rmfamily\fontsize{10.000000}{12.000000}\selectfont \(\displaystyle {0.4}\)}%
\end{pgfscope}%
\begin{pgfscope}%
\pgfsetbuttcap%
\pgfsetroundjoin%
\definecolor{currentfill}{rgb}{0.000000,0.000000,0.000000}%
\pgfsetfillcolor{currentfill}%
\pgfsetlinewidth{0.803000pt}%
\definecolor{currentstroke}{rgb}{0.000000,0.000000,0.000000}%
\pgfsetstrokecolor{currentstroke}%
\pgfsetdash{}{0pt}%
\pgfsys@defobject{currentmarker}{\pgfqpoint{0.000000in}{-0.048611in}}{\pgfqpoint{0.000000in}{0.000000in}}{%
\pgfpathmoveto{\pgfqpoint{0.000000in}{0.000000in}}%
\pgfpathlineto{\pgfqpoint{0.000000in}{-0.048611in}}%
\pgfusepath{stroke,fill}%
}%
\begin{pgfscope}%
\pgfsys@transformshift{3.281250in}{0.450000in}%
\pgfsys@useobject{currentmarker}{}%
\end{pgfscope}%
\end{pgfscope}%
\begin{pgfscope}%
\definecolor{textcolor}{rgb}{0.000000,0.000000,0.000000}%
\pgfsetstrokecolor{textcolor}%
\pgfsetfillcolor{textcolor}%
\pgftext[x=3.281250in,y=0.352778in,,top]{\color{textcolor}\rmfamily\fontsize{10.000000}{12.000000}\selectfont \(\displaystyle {0.5}\)}%
\end{pgfscope}%
\begin{pgfscope}%
\definecolor{textcolor}{rgb}{0.000000,0.000000,0.000000}%
\pgfsetstrokecolor{textcolor}%
\pgfsetfillcolor{textcolor}%
\pgftext[x=1.968750in,y=0.162809in,,top]{\color{textcolor}\rmfamily\fontsize{10.000000}{12.000000}\selectfont \(\displaystyle t_{\mathrm{blend}}/\Delta t_{\mathrm{obs}}\)}%
\end{pgfscope}%
\begin{pgfscope}%
\pgfsetbuttcap%
\pgfsetroundjoin%
\definecolor{currentfill}{rgb}{0.000000,0.000000,0.000000}%
\pgfsetfillcolor{currentfill}%
\pgfsetlinewidth{0.803000pt}%
\definecolor{currentstroke}{rgb}{0.000000,0.000000,0.000000}%
\pgfsetstrokecolor{currentstroke}%
\pgfsetdash{}{0pt}%
\pgfsys@defobject{currentmarker}{\pgfqpoint{-0.048611in}{0.000000in}}{\pgfqpoint{0.000000in}{0.000000in}}{%
\pgfpathmoveto{\pgfqpoint{0.000000in}{0.000000in}}%
\pgfpathlineto{\pgfqpoint{-0.048611in}{0.000000in}}%
\pgfusepath{stroke,fill}%
}%
\begin{pgfscope}%
\pgfsys@transformshift{0.525000in}{0.488244in}%
\pgfsys@useobject{currentmarker}{}%
\end{pgfscope}%
\end{pgfscope}%
\begin{pgfscope}%
\definecolor{textcolor}{rgb}{0.000000,0.000000,0.000000}%
\pgfsetstrokecolor{textcolor}%
\pgfsetfillcolor{textcolor}%
\pgftext[x=0.250308in, y=0.435483in, left, base]{\color{textcolor}\rmfamily\fontsize{10.000000}{12.000000}\selectfont \(\displaystyle {0.0}\)}%
\end{pgfscope}%
\begin{pgfscope}%
\pgfsetbuttcap%
\pgfsetroundjoin%
\definecolor{currentfill}{rgb}{0.000000,0.000000,0.000000}%
\pgfsetfillcolor{currentfill}%
\pgfsetlinewidth{0.803000pt}%
\definecolor{currentstroke}{rgb}{0.000000,0.000000,0.000000}%
\pgfsetstrokecolor{currentstroke}%
\pgfsetdash{}{0pt}%
\pgfsys@defobject{currentmarker}{\pgfqpoint{-0.048611in}{0.000000in}}{\pgfqpoint{0.000000in}{0.000000in}}{%
\pgfpathmoveto{\pgfqpoint{0.000000in}{0.000000in}}%
\pgfpathlineto{\pgfqpoint{-0.048611in}{0.000000in}}%
\pgfusepath{stroke,fill}%
}%
\begin{pgfscope}%
\pgfsys@transformshift{0.525000in}{0.969635in}%
\pgfsys@useobject{currentmarker}{}%
\end{pgfscope}%
\end{pgfscope}%
\begin{pgfscope}%
\definecolor{textcolor}{rgb}{0.000000,0.000000,0.000000}%
\pgfsetstrokecolor{textcolor}%
\pgfsetfillcolor{textcolor}%
\pgftext[x=0.250308in, y=0.916873in, left, base]{\color{textcolor}\rmfamily\fontsize{10.000000}{12.000000}\selectfont \(\displaystyle {0.1}\)}%
\end{pgfscope}%
\begin{pgfscope}%
\pgfsetbuttcap%
\pgfsetroundjoin%
\definecolor{currentfill}{rgb}{0.000000,0.000000,0.000000}%
\pgfsetfillcolor{currentfill}%
\pgfsetlinewidth{0.803000pt}%
\definecolor{currentstroke}{rgb}{0.000000,0.000000,0.000000}%
\pgfsetstrokecolor{currentstroke}%
\pgfsetdash{}{0pt}%
\pgfsys@defobject{currentmarker}{\pgfqpoint{-0.048611in}{0.000000in}}{\pgfqpoint{0.000000in}{0.000000in}}{%
\pgfpathmoveto{\pgfqpoint{0.000000in}{0.000000in}}%
\pgfpathlineto{\pgfqpoint{-0.048611in}{0.000000in}}%
\pgfusepath{stroke,fill}%
}%
\begin{pgfscope}%
\pgfsys@transformshift{0.525000in}{1.451026in}%
\pgfsys@useobject{currentmarker}{}%
\end{pgfscope}%
\end{pgfscope}%
\begin{pgfscope}%
\definecolor{textcolor}{rgb}{0.000000,0.000000,0.000000}%
\pgfsetstrokecolor{textcolor}%
\pgfsetfillcolor{textcolor}%
\pgftext[x=0.250308in, y=1.398264in, left, base]{\color{textcolor}\rmfamily\fontsize{10.000000}{12.000000}\selectfont \(\displaystyle {0.2}\)}%
\end{pgfscope}%
\begin{pgfscope}%
\pgfsetbuttcap%
\pgfsetroundjoin%
\definecolor{currentfill}{rgb}{0.000000,0.000000,0.000000}%
\pgfsetfillcolor{currentfill}%
\pgfsetlinewidth{0.803000pt}%
\definecolor{currentstroke}{rgb}{0.000000,0.000000,0.000000}%
\pgfsetstrokecolor{currentstroke}%
\pgfsetdash{}{0pt}%
\pgfsys@defobject{currentmarker}{\pgfqpoint{-0.048611in}{0.000000in}}{\pgfqpoint{0.000000in}{0.000000in}}{%
\pgfpathmoveto{\pgfqpoint{0.000000in}{0.000000in}}%
\pgfpathlineto{\pgfqpoint{-0.048611in}{0.000000in}}%
\pgfusepath{stroke,fill}%
}%
\begin{pgfscope}%
\pgfsys@transformshift{0.525000in}{1.932416in}%
\pgfsys@useobject{currentmarker}{}%
\end{pgfscope}%
\end{pgfscope}%
\begin{pgfscope}%
\definecolor{textcolor}{rgb}{0.000000,0.000000,0.000000}%
\pgfsetstrokecolor{textcolor}%
\pgfsetfillcolor{textcolor}%
\pgftext[x=0.250308in, y=1.879655in, left, base]{\color{textcolor}\rmfamily\fontsize{10.000000}{12.000000}\selectfont \(\displaystyle {0.3}\)}%
\end{pgfscope}%
\begin{pgfscope}%
\pgfsetbuttcap%
\pgfsetroundjoin%
\definecolor{currentfill}{rgb}{0.000000,0.000000,0.000000}%
\pgfsetfillcolor{currentfill}%
\pgfsetlinewidth{0.803000pt}%
\definecolor{currentstroke}{rgb}{0.000000,0.000000,0.000000}%
\pgfsetstrokecolor{currentstroke}%
\pgfsetdash{}{0pt}%
\pgfsys@defobject{currentmarker}{\pgfqpoint{-0.048611in}{0.000000in}}{\pgfqpoint{0.000000in}{0.000000in}}{%
\pgfpathmoveto{\pgfqpoint{0.000000in}{0.000000in}}%
\pgfpathlineto{\pgfqpoint{-0.048611in}{0.000000in}}%
\pgfusepath{stroke,fill}%
}%
\begin{pgfscope}%
\pgfsys@transformshift{0.525000in}{2.413807in}%
\pgfsys@useobject{currentmarker}{}%
\end{pgfscope}%
\end{pgfscope}%
\begin{pgfscope}%
\definecolor{textcolor}{rgb}{0.000000,0.000000,0.000000}%
\pgfsetstrokecolor{textcolor}%
\pgfsetfillcolor{textcolor}%
\pgftext[x=0.250308in, y=2.361046in, left, base]{\color{textcolor}\rmfamily\fontsize{10.000000}{12.000000}\selectfont \(\displaystyle {0.4}\)}%
\end{pgfscope}%
\begin{pgfscope}%
\definecolor{textcolor}{rgb}{0.000000,0.000000,0.000000}%
\pgfsetstrokecolor{textcolor}%
\pgfsetfillcolor{textcolor}%
\pgftext[x=0.194753in,y=1.650000in,,bottom,rotate=90.000000]{\color{textcolor}\rmfamily\fontsize{10.000000}{12.000000}\selectfont \(\displaystyle \mathrm{TRMSE}(p_\parallel)\)}%
\end{pgfscope}%
\begin{pgfscope}%
\pgfpathrectangle{\pgfqpoint{0.525000in}{0.450000in}}{\pgfqpoint{2.887500in}{2.400000in}}%
\pgfusepath{clip}%
\pgfsetrectcap%
\pgfsetroundjoin%
\pgfsetlinewidth{1.003750pt}%
\definecolor{currentstroke}{rgb}{1.000000,0.600000,0.000000}%
\pgfsetstrokecolor{currentstroke}%
\pgfsetdash{}{0pt}%
\pgfpathmoveto{\pgfqpoint{0.656250in}{2.740909in}}%
\pgfpathlineto{\pgfqpoint{3.281250in}{2.740909in}}%
\pgfusepath{stroke}%
\end{pgfscope}%
\begin{pgfscope}%
\pgfpathrectangle{\pgfqpoint{0.525000in}{0.450000in}}{\pgfqpoint{2.887500in}{2.400000in}}%
\pgfusepath{clip}%
\pgfsetrectcap%
\pgfsetroundjoin%
\pgfsetlinewidth{1.003750pt}%
\definecolor{currentstroke}{rgb}{0.000000,0.200000,0.400000}%
\pgfsetstrokecolor{currentstroke}%
\pgfsetdash{}{0pt}%
\pgfpathmoveto{\pgfqpoint{0.656250in}{2.733229in}}%
\pgfpathlineto{\pgfqpoint{0.984375in}{1.084630in}}%
\pgfpathlineto{\pgfqpoint{1.312500in}{0.900784in}}%
\pgfpathlineto{\pgfqpoint{1.640625in}{0.785834in}}%
\pgfpathlineto{\pgfqpoint{1.968750in}{0.565358in}}%
\pgfpathlineto{\pgfqpoint{2.296875in}{0.563550in}}%
\pgfpathlineto{\pgfqpoint{2.625000in}{0.565488in}}%
\pgfpathlineto{\pgfqpoint{2.953125in}{0.561530in}}%
\pgfpathlineto{\pgfqpoint{3.281250in}{0.559091in}}%
\pgfusepath{stroke}%
\end{pgfscope}%
\begin{pgfscope}%
\pgfpathrectangle{\pgfqpoint{0.525000in}{0.450000in}}{\pgfqpoint{2.887500in}{2.400000in}}%
\pgfusepath{clip}%
\pgfsetbuttcap%
\pgfsetroundjoin%
\definecolor{currentfill}{rgb}{0.000000,0.200000,0.400000}%
\pgfsetfillcolor{currentfill}%
\pgfsetlinewidth{1.003750pt}%
\definecolor{currentstroke}{rgb}{0.000000,0.200000,0.400000}%
\pgfsetstrokecolor{currentstroke}%
\pgfsetdash{}{0pt}%
\pgfsys@defobject{currentmarker}{\pgfqpoint{-0.027778in}{-0.027778in}}{\pgfqpoint{0.027778in}{0.027778in}}{%
\pgfpathmoveto{\pgfqpoint{-0.027778in}{-0.027778in}}%
\pgfpathlineto{\pgfqpoint{0.027778in}{0.027778in}}%
\pgfpathmoveto{\pgfqpoint{-0.027778in}{0.027778in}}%
\pgfpathlineto{\pgfqpoint{0.027778in}{-0.027778in}}%
\pgfusepath{stroke,fill}%
}%
\begin{pgfscope}%
\pgfsys@transformshift{0.656250in}{2.733229in}%
\pgfsys@useobject{currentmarker}{}%
\end{pgfscope}%
\begin{pgfscope}%
\pgfsys@transformshift{0.984375in}{1.084630in}%
\pgfsys@useobject{currentmarker}{}%
\end{pgfscope}%
\begin{pgfscope}%
\pgfsys@transformshift{1.312500in}{0.900784in}%
\pgfsys@useobject{currentmarker}{}%
\end{pgfscope}%
\begin{pgfscope}%
\pgfsys@transformshift{1.640625in}{0.785834in}%
\pgfsys@useobject{currentmarker}{}%
\end{pgfscope}%
\begin{pgfscope}%
\pgfsys@transformshift{1.968750in}{0.565358in}%
\pgfsys@useobject{currentmarker}{}%
\end{pgfscope}%
\begin{pgfscope}%
\pgfsys@transformshift{2.296875in}{0.563550in}%
\pgfsys@useobject{currentmarker}{}%
\end{pgfscope}%
\begin{pgfscope}%
\pgfsys@transformshift{2.625000in}{0.565488in}%
\pgfsys@useobject{currentmarker}{}%
\end{pgfscope}%
\begin{pgfscope}%
\pgfsys@transformshift{2.953125in}{0.561530in}%
\pgfsys@useobject{currentmarker}{}%
\end{pgfscope}%
\begin{pgfscope}%
\pgfsys@transformshift{3.281250in}{0.559091in}%
\pgfsys@useobject{currentmarker}{}%
\end{pgfscope}%
\end{pgfscope}%
\begin{pgfscope}%
\pgfsetrectcap%
\pgfsetmiterjoin%
\pgfsetlinewidth{0.803000pt}%
\definecolor{currentstroke}{rgb}{0.000000,0.000000,0.000000}%
\pgfsetstrokecolor{currentstroke}%
\pgfsetdash{}{0pt}%
\pgfpathmoveto{\pgfqpoint{0.525000in}{0.450000in}}%
\pgfpathlineto{\pgfqpoint{0.525000in}{2.850000in}}%
\pgfusepath{stroke}%
\end{pgfscope}%
\begin{pgfscope}%
\pgfsetrectcap%
\pgfsetmiterjoin%
\pgfsetlinewidth{0.803000pt}%
\definecolor{currentstroke}{rgb}{0.000000,0.000000,0.000000}%
\pgfsetstrokecolor{currentstroke}%
\pgfsetdash{}{0pt}%
\pgfpathmoveto{\pgfqpoint{3.412500in}{0.450000in}}%
\pgfpathlineto{\pgfqpoint{3.412500in}{2.850000in}}%
\pgfusepath{stroke}%
\end{pgfscope}%
\begin{pgfscope}%
\pgfsetrectcap%
\pgfsetmiterjoin%
\pgfsetlinewidth{0.803000pt}%
\definecolor{currentstroke}{rgb}{0.000000,0.000000,0.000000}%
\pgfsetstrokecolor{currentstroke}%
\pgfsetdash{}{0pt}%
\pgfpathmoveto{\pgfqpoint{0.525000in}{0.450000in}}%
\pgfpathlineto{\pgfqpoint{3.412500in}{0.450000in}}%
\pgfusepath{stroke}%
\end{pgfscope}%
\begin{pgfscope}%
\pgfsetrectcap%
\pgfsetmiterjoin%
\pgfsetlinewidth{0.803000pt}%
\definecolor{currentstroke}{rgb}{0.000000,0.000000,0.000000}%
\pgfsetstrokecolor{currentstroke}%
\pgfsetdash{}{0pt}%
\pgfpathmoveto{\pgfqpoint{0.525000in}{2.850000in}}%
\pgfpathlineto{\pgfqpoint{3.412500in}{2.850000in}}%
\pgfusepath{stroke}%
\end{pgfscope}%
\begin{pgfscope}%
\pgfsetbuttcap%
\pgfsetmiterjoin%
\definecolor{currentfill}{rgb}{1.000000,1.000000,1.000000}%
\pgfsetfillcolor{currentfill}%
\pgfsetfillopacity{0.800000}%
\pgfsetlinewidth{1.003750pt}%
\definecolor{currentstroke}{rgb}{0.800000,0.800000,0.800000}%
\pgfsetstrokecolor{currentstroke}%
\pgfsetstrokeopacity{0.800000}%
\pgfsetdash{}{0pt}%
\pgfpathmoveto{\pgfqpoint{2.613704in}{1.514596in}}%
\pgfpathlineto{\pgfqpoint{3.354167in}{1.514596in}}%
\pgfpathquadraticcurveto{\pgfqpoint{3.370833in}{1.514596in}}{\pgfqpoint{3.370833in}{1.531262in}}%
\pgfpathlineto{\pgfqpoint{3.370833in}{1.768738in}}%
\pgfpathquadraticcurveto{\pgfqpoint{3.370833in}{1.785404in}}{\pgfqpoint{3.354167in}{1.785404in}}%
\pgfpathlineto{\pgfqpoint{2.613704in}{1.785404in}}%
\pgfpathquadraticcurveto{\pgfqpoint{2.597038in}{1.785404in}}{\pgfqpoint{2.597038in}{1.768738in}}%
\pgfpathlineto{\pgfqpoint{2.597038in}{1.531262in}}%
\pgfpathquadraticcurveto{\pgfqpoint{2.597038in}{1.514596in}}{\pgfqpoint{2.613704in}{1.514596in}}%
\pgfpathclose%
\pgfusepath{stroke,fill}%
\end{pgfscope}%
\begin{pgfscope}%
\pgfsetrectcap%
\pgfsetroundjoin%
\pgfsetlinewidth{1.003750pt}%
\definecolor{currentstroke}{rgb}{1.000000,0.600000,0.000000}%
\pgfsetstrokecolor{currentstroke}%
\pgfsetdash{}{0pt}%
\pgfpathmoveto{\pgfqpoint{2.630371in}{1.717924in}}%
\pgfpathlineto{\pgfqpoint{2.797038in}{1.717924in}}%
\pgfusepath{stroke}%
\end{pgfscope}%
\begin{pgfscope}%
\definecolor{textcolor}{rgb}{0.000000,0.000000,0.000000}%
\pgfsetstrokecolor{textcolor}%
\pgfsetfillcolor{textcolor}%
\pgftext[x=2.863704in,y=1.688757in,left,base]{\color{textcolor}\rmfamily\fontsize{6.000000}{7.200000}\selectfont Plain EnKF}%
\end{pgfscope}%
\begin{pgfscope}%
\pgfsetrectcap%
\pgfsetroundjoin%
\pgfsetlinewidth{1.003750pt}%
\definecolor{currentstroke}{rgb}{0.000000,0.200000,0.400000}%
\pgfsetstrokecolor{currentstroke}%
\pgfsetdash{}{0pt}%
\pgfpathmoveto{\pgfqpoint{2.630371in}{1.595610in}}%
\pgfpathlineto{\pgfqpoint{2.797038in}{1.595610in}}%
\pgfusepath{stroke}%
\end{pgfscope}%
\begin{pgfscope}%
\pgfsetbuttcap%
\pgfsetroundjoin%
\definecolor{currentfill}{rgb}{0.000000,0.200000,0.400000}%
\pgfsetfillcolor{currentfill}%
\pgfsetlinewidth{1.003750pt}%
\definecolor{currentstroke}{rgb}{0.000000,0.200000,0.400000}%
\pgfsetstrokecolor{currentstroke}%
\pgfsetdash{}{0pt}%
\pgfsys@defobject{currentmarker}{\pgfqpoint{-0.027778in}{-0.027778in}}{\pgfqpoint{0.027778in}{0.027778in}}{%
\pgfpathmoveto{\pgfqpoint{-0.027778in}{-0.027778in}}%
\pgfpathlineto{\pgfqpoint{0.027778in}{0.027778in}}%
\pgfpathmoveto{\pgfqpoint{-0.027778in}{0.027778in}}%
\pgfpathlineto{\pgfqpoint{0.027778in}{-0.027778in}}%
\pgfusepath{stroke,fill}%
}%
\begin{pgfscope}%
\pgfsys@transformshift{2.713704in}{1.595610in}%
\pgfsys@useobject{currentmarker}{}%
\end{pgfscope}%
\end{pgfscope}%
\begin{pgfscope}%
\definecolor{textcolor}{rgb}{0.000000,0.000000,0.000000}%
\pgfsetstrokecolor{textcolor}%
\pgfsetfillcolor{textcolor}%
\pgftext[x=2.863704in,y=1.566443in,left,base]{\color{textcolor}\rmfamily\fontsize{6.000000}{7.200000}\selectfont Blending}%
\end{pgfscope}%
\end{pgfpicture}%
\makeatother%
\endgroup%

%% file: publication.bbl
\begin{thebibliography}{52}
\providecommand{\natexlab}[1]{#1}
\providecommand{\url}[1]{\texttt{#1}}
\expandafter\ifx\csname urlstyle\endcsname\relax
  \providecommand{\doi}[1]{doi: #1}\else
  \providecommand{\doi}{doi: \begingroup \urlstyle{rm}\Url}\fi

\bibitem[Ascher and Petzold(1998)]{Ascher1998}
U.~M. Ascher and L.~R. Petzold.
\newblock \emph{Computer Methods for Ordinary Differential Equations and
  Differential-Algebraic Equations}.
\newblock {Society for Industrial and Applied Mathematics}, {Philadelphia},
  1998.

\bibitem[Ascher et~al.(1994)Ascher, Chin, and Reich]{ascher1994}
U.~M. Ascher, H.~Chin, and S.~Reich.
\newblock Stabilization of {{DAEs}} and invariant manifolds.
\newblock \emph{Numerische Mathematik}, 67\penalty0 (2):\penalty0 131--149,
  1994.

\bibitem[Benacchio et~al.(2014)Benacchio, O'Neill, and
  Klein]{BenacchioEtAl2014}
T.~Benacchio, W.~P. O'Neill, and R.~Klein.
\newblock {A} blended soundproof-to-compressible numerical model for small- to
  mesoscale atmospheric dynamics.
\newblock \emph{{M}onthly {W}eather {R}eview}, 142\penalty0 (12):\penalty0
  4416--4438, 2014.

\bibitem[Benettin et~al.(1987)Benettin, Galgani, and
  Giorgilli]{BenettinEtAl1987}
G.~Benettin, L.~Galgani, and A.~Giorgilli.
\newblock {R}ealization of holonomic constraints and freezing of high frequency
  degrees of freedom in the light of classical perturbation theory .1.
\newblock \emph{{C}ommunications in {M}athematical {P}hysics}, 113\penalty0
  (1):\penalty0 87--103, 1987.

\bibitem[Benettin et~al.(1989)Benettin, Galgani, and
  Giorgilli]{BenettinEtAl1989}
G.~Benettin, L.~Galgani, and A.~Giorgilli.
\newblock {R}ealization of holonomic constraints and freezing of high frequency
  degrees of freedom in the light of classical perturbation theory .2.
\newblock \emph{{C}ommunications in {M}athematical {P}hysics}, 121\penalty0
  (4):\penalty0 557--601, 1989.

\bibitem[Bergemann and Reich(2010)]{Kay2010}
K.~Bergemann and S.~Reich.
\newblock A mollified ensemble {K}alman filter.
\newblock \emph{Quarterly Journal of the Royal Meteorological Society},
  136:\penalty0 1636--1643, 2010.

\bibitem[Bloom et~al.(1996)Bloom, Takacs, Da~Silva, and Ledvina]{Bloom1996}
S.~Bloom, L.~L. Takacs, A.~Da~Silva, and D.~Ledvina.
\newblock Data assimilation using incremental analysis updates.
\newblock \emph{Monthly Weather Review}, 124:\penalty0 1256--1271, 1996.

\bibitem[Bokhove and Shepherd(1996)]{Bokhove1996}
O.~Bokhove and T.~G. Shepherd.
\newblock On hamiltonian balanced dynamics and the slowest invariant manifold.
\newblock \emph{Journal of the Atmospheric Sciences}, 53\penalty0 (2):\penalty0
  276--297, 1996.

\bibitem[Bornemann and Sch\"utte(1997)]{Bornemann1997}
F.~A. Bornemann and C.~Sch\"utte.
\newblock Homogenization of {H}amiltonian systems with a strong constraining
  potential.
\newblock \emph{Physica D}, 102\penalty0 (1-2):\penalty0 57--77, 1997.

\bibitem[Broyden(1970)]{broyden_convergence_1970}
C.~G. Broyden.
\newblock The {Convergence} of a {Class} of {Double}-rank {Minimization}
  {Algorithms} 1. {General} {Considerations}.
\newblock \emph{IMA Journal of Applied Mathematics}, 6\penalty0 (1):\penalty0
  76--90, 1970.

\bibitem[Camassa(1995)]{Camassa1995}
R.~Camassa.
\newblock On the geometry of an atmospheric slow manifold.
\newblock \emph{Physica D: Nonlinear Phenomena}, 84\penalty0 (3-4):\penalty0
  357--397, 1995.

\bibitem[Chorin(1967)]{Chorin1967}
A.~J. Chorin.
\newblock The numerical solution of {N}avier-{S}tokes equations for an
  imcompressible fluid.
\newblock \emph{{B}ulletin of the {A}merican {M}athematical {S}ociety},
  73\penalty0 (6):\penalty0 928--931, 1967.

\bibitem[Chorin and Tu(2009)]{ChorinTu2009}
A.~J. Chorin and X.~Tu.
\newblock Implicit sampling for particle filters.
\newblock \emph{P. Natl. Acad. Sci. USA}, 106\penalty0 (41):\penalty0
  17249--17254, 2009.

\bibitem[Chorin et~al.(2010)Chorin, Morzfeld, and Tu]{ChorinEtAl2010}
A.~J. Chorin, M.~Morzfeld, and X.~Tu.
\newblock Implicit particle filters for data assimilation.
\newblock \emph{Comm. Appl. Math. Comput. Sci.}, 5\penalty0 (2):\penalty0
  221--240, 2010.

\bibitem[Cordier et~al.(2012)Cordier, Degond, and Kumbaro]{CordierEtAl2012}
F.~Cordier, P.~Degond, and A.~Kumbaro.
\newblock An {A}symptotic-{P}reserving all-speed scheme for the {E}uler and
  {N}avier--{S}tokes equations.
\newblock \emph{J. Comput. Phys.}, 231:\penalty0 5685---5704, 2012.

\bibitem[Cotter(2013)]{Cotter2013}
C.~Cotter.
\newblock Data assimilation on the exponentially accurate slow manifold.
\newblock \emph{Philosophical Transactions of the Royal Society A:
  Mathematical, Physical and Engineering Sciences}, 371\penalty0
  (1991):\penalty0 20120300, 2013.

\bibitem[Evensen(2003)]{Evensen2003}
G.~Evensen.
\newblock {T}he {E}nsemble {K}alman {F}ilter: theoretical formulation and
  practical implementation.
\newblock \emph{Ocean Dynamics}, 53:\penalty0 343, 2003.

\bibitem[Fenichel(1979)]{Fenichel1979}
N.~Fenichel.
\newblock Geometric singular perturbation theory for ordinary differential
  equations.
\newblock \emph{Journal of Differential Equations}, 31\penalty0 (1):\penalty0
  53--98, 1979.

\bibitem[Fletcher(1970)]{fletcher_new_1970}
R.~Fletcher.
\newblock A new approach to variable metric algorithms.
\newblock \emph{The Computer Journal}, 13\penalty0 (3):\penalty0 317--322,
  1970.

\bibitem[Gear(1986)]{Gear1986}
C.~W. Gear.
\newblock Maintaining solution invariants in the numerical solution of odes.
\newblock \emph{SIAM J. Sci. Stat. Comput.}, 7\penalty0 (3):\penalty0 734--743,
  1986.

\bibitem[Goldfarb(1970)]{goldfarb_family_1970}
D.~Goldfarb.
\newblock A family of variable-metric methods derived by variational means.
\newblock \emph{Mathematics of Computation}, 24\penalty0 (109):\penalty0
  23--23, 1970.

\bibitem[Golub and Van~Loan(1996)]{golub1996}
G.~H. Golub and C.~F. Van~Loan.
\newblock \emph{Matrix Computations}.
\newblock Johns {{Hopkins}} Studies in the Mathematical Sciences. {Johns
  Hopkins University Press}, {Baltimore}, 3rd ed edition, 1996.

\bibitem[Gottwald(2014)]{gottwald2014}
G.~A. Gottwald.
\newblock Controlling balance in an ensemble {{Kalman}} filter.
\newblock \emph{Nonlinear Processes in Geophysics}, 21\penalty0 (2):\penalty0
  417--426, 2014.

\bibitem[Gottwald et~al.(2011)Gottwald, Mitchell, and Reich]{GottwaldEtAl2011}
G.~A. Gottwald, L.~Mitchell, and S.~Reich.
\newblock Controlling overestimation of error covariance in ensemble kalman
  filters with sparse observations: A variance-limiting kalman filter.
\newblock \emph{Monthly Weather Review}, 139\penalty0 (8):\penalty0 2650--2667,
  2011.
\newblock \doi{10.1175/2011MWR3557.1}.
\newblock URL
  \url{https://journals.ametsoc.org/view/journals/mwre/139/8/2011mwr3557.1.xml}.

\bibitem[Hairer et~al.(2010)Hairer, Lubich, and Wanner]{HairerEtAl2010}
E.~Hairer, C.~Lubich, and G.~Wanner.
\newblock \emph{Geometric numerical integration: structure-preserving
  algorithms for ordinary differential equations}.
\newblock Springer series in computational mathematics. Springer, Berlin u.a.,
  2. edition, 2010.

\bibitem[Harlim and Majda(2010)]{HarlimMajda2010}
J.~Harlim and A.~J. Majda.
\newblock Catastrophic filter divergence in filtering nonlinear dissipative
  systems.
\newblock \emph{Communications in Mathematical Sciences}, 8\penalty0
  (1):\penalty0 27--43, 2010.

\bibitem[Jin(2012)]{Jin2012}
S.~Jin.
\newblock Asymptotic preserving (ap) schemes for multiscale kinetic and
  hyperbolic equations: a review.
\newblock \emph{Riv. Mat. Univ. Parma}, 3:\penalty0 177--216, 2012.

\bibitem[Kalman(1960)]{Kalman1960}
R.~E. Kalman.
\newblock A new approach to linear filtering and prediction problems.
\newblock \emph{Journal of Basic Engineering}, 82:\penalty0 35, 1960.

\bibitem[Kalnay(2002)]{kalnay_atmospheric_2002}
E.~Kalnay.
\newblock \emph{Atmospheric {Modeling}, {Data} {Assimilation} and
  {Predictability}}.
\newblock Cambridge University Press, 1 edition, 2002.

\bibitem[Kelly et~al.(2015)Kelly, Majda, and Tong]{KellyEtAl2015}
D.~Kelly, A.~J. Majda, and X.~T. Tong.
\newblock Concrete ensemble {K}alman filters with rigorous catastrophic filter
  divergence.
\newblock \emph{P. Natl. Acad. Sci. USA}, 112\penalty0 (34):\penalty0
  10589--10594, 2015.

\bibitem[Kepert(2009)]{kepert2009}
J.~D. Kepert.
\newblock Covariance localisation and balance in an {{Ensemble Kalman Filter}}.
\newblock \emph{Quarterly Journal of the Royal Meteorological Society},
  135\penalty0 (642):\penalty0 1157--1176, 2009.

\bibitem[Klein(2010)]{Klein2010}
R.~Klein.
\newblock Scale-dependent asymptotic models for atmospheric flows.
\newblock \emph{Ann. Rev. Fluid Mech.}, 42:\penalty0 249--274, 2010.

\bibitem[Klein et~al.(2001)Klein, Botta, Hofmann, Meister, Munz, Roller, and
  Sonar]{KleinEtAl2001}
R.~Klein, N.~Botta, L.~Hofmann, A.~Meister, C.~Munz, S.~Roller, and T.~Sonar.
\newblock Asymptotic adaptive methods for multiscale problems in fluid
  mechanics.
\newblock \emph{J. Engrg. Math.}, 39:\penalty0 261--343, 2001.

\bibitem[Kuehn(2015)]{Kuehn2015}
C.~Kuehn.
\newblock \emph{Geometric Singular Perturbation Theory}, pages 53--70.
\newblock Springer International Publishing, Cham, 2015.

\bibitem[Leimkuhler and Reich(2005)]{leimkuhler_reich_2005}
B.~Leimkuhler and S.~Reich.
\newblock \emph{Simulating Hamiltonian Dynamics}.
\newblock Cambridge Monographs on Applied and Computational Mathematics.
  Cambridge University Press, 2005.

\bibitem[Leimkuhler and Skeel(1994)]{leimkuhler1994}
B.~J. Leimkuhler and R.~D. Skeel.
\newblock Symplectic {{Numerical Integrators}} in {{Constrained Hamiltonian
  Systems}}.
\newblock \emph{Journal of Computational Physics}, 112\penalty0 (1):\penalty0
  117--125, 1994.

\bibitem[Lorenz(1963)]{Lorenz1963}
E.~N. Lorenz.
\newblock Deterministic nonperiodic flow.
\newblock \emph{Journal of the Atmospheric Sciences}, 20\penalty0 (2):\penalty0
  130--141, 1963.

\bibitem[Lorenz(2006)]{Lorenz2006}
E.~N. Lorenz.
\newblock Predictability -- a problem partly solved.
\newblock In T.~Palmer and R.~Hagedorn, editors, \emph{Predictability of
  Weather and Climate}, pages 40--58. Cambridge University Press, 2006.

\bibitem[Lynch(2002)]{Lynch2002}
P.~Lynch.
\newblock The swinging spring: a simple model for atmospheric balance.
\newblock In J.~Norbury and I.~Roulstone, editors, \emph{Large-Scale
  Atmosphere-Ocean Dynamics: Volume II: Geometric Methods and Models}, page~64,
  2002.

\bibitem[Lynch(2014)]{Lynch2014}
P.~Lynch.
\newblock \emph{The Emergence of Numerical Weather Prediction: Richardson's
  Dream}.
\newblock Cambridge University Press, 2014.

\bibitem[Lynch and Huang(1992)]{Lynch1992}
P.~Lynch and X.-Y. Huang.
\newblock Initialization of the {H}{I}{R}{L}{A}{M} model using a digital
  filter.
\newblock \emph{Monthly Weather Review}, 120:\penalty0 1019--1034, 1992.

\bibitem[Mor{\'e} et~al.(1980)Mor{\'e}, Garbow, and Hillstrom]{More:126569}
J.~J. Mor{\'e}, B.~S. Garbow, and K.~E. Hillstrom.
\newblock {User guide for MINPACK-1}.
\newblock Technical Report ANL-80-74, Argonne Nat. Lab., Argonne, IL, 1980.

\bibitem[Reich(1995)]{Reich1995}
S.~Reich.
\newblock Smoothed dynamics of highly oscillatory {H}amiltonian systems.
\newblock \emph{Physica D}, 89:\penalty0 28---42, 1995.

\bibitem[Reich(2000)]{Reich1999}
S.~Reich.
\newblock {S}moothed {L}angevin dynamics of highly oscillatory systems.
\newblock \emph{Physica D}, 138:\penalty0 210--224, 2000.

\bibitem[Reich and Cotter(2015)]{reichbook}
S.~Reich and C.~Cotter.
\newblock \emph{Probabilistic forecasting and {B}ayesian data assimilation: a
  tutorial}.
\newblock Cambridge University Press, 2015.

\bibitem[Rubin and Ungar(1957)]{RubinUngar1957}
H.~Rubin and P.~Ungar.
\newblock {M}otion under a {S}trong {C}onstraining {F}orce.
\newblock \emph{{C}ommunications on {P}ure and {A}pplied {M}athematics},
  10\penalty0 (1):\penalty0 65--87, 1957.

\bibitem[Shanno(1970)]{shanno_conditioning_1970}
D.~F. Shanno.
\newblock Conditioning of quasi-{Newton} methods for function minimization.
\newblock \emph{Mathematics of Computation}, 24\penalty0 (111):\penalty0
  647--647, 1970.

\bibitem[Takens(1980)]{Takens1980}
F.~Takens.
\newblock \emph{Motion under the influence of a strong constraining force}.
\newblock Springer Berlin Heidelberg, 1980.

\bibitem[{Virtanen} et~al.(2020){Virtanen}, {Gommers}, {Oliphant}, {Haberland},
  {Reddy}, {Cournapeau}, {Burovski}, {Peterson}, {Weckesser}, {Bright}, {van
  der Walt}, {Brett}, {Wilson}, {Jarrod Millman}, {Mayorov}, {Nelson}, {Jones},
  {Kern}, {Larson}, {Carey}, {Polat}, {Feng}, {Moore}, {VanderPlas}, {Laxalde},
  {Perktold}, {Cimrman}, {Henriksen}, {Quintero}, {Harris}, {Archibald},
  {Ribeiro}, {Pedregosa}, {van Mulbregt}, and {Contributors}]{2020SciPy-NMeth}
P.~{Virtanen}, R.~{Gommers}, T.~E. {Oliphant}, M.~{Haberland}, T.~{Reddy},
  D.~{Cournapeau}, E.~{Burovski}, P.~{Peterson}, W.~{Weckesser}, J.~{Bright},
  S.~J. {van der Walt}, M.~{Brett}, J.~{Wilson}, K.~{Jarrod Millman},
  N.~{Mayorov}, A.~R.~J. {Nelson}, E.~{Jones}, R.~{Kern}, E.~{Larson},
  C.~{Carey}, {\.I}.~{Polat}, Y.~{Feng}, E.~W. {Moore}, J.~{VanderPlas},
  D.~{Laxalde}, J.~{Perktold}, R.~{Cimrman}, I.~{Henriksen}, E.~A. {Quintero},
  C.~R. {Harris}, A.~M. {Archibald}, A.~H. {Ribeiro}, F.~{Pedregosa}, P.~{van
  Mulbregt}, and S.~.~. {Contributors}.
\newblock {SciPy 1.0: Fundamental Algorithms for Scientific Computing in
  Python}.
\newblock \emph{Nature Methods}, 17:\penalty0 261--272, 2020.

\bibitem[Warming and Hyett(1974)]{WarmingHyett1974}
R.~F. Warming and F.~Hyett.
\newblock The modified equation approach to the stability and accuracy analysis
  of finite difference methods.
\newblock \emph{J. Comput. Phys}, 14\penalty0 (2):\penalty0 159--179, 1974.

\bibitem[Whitaker and Hamill(2002)]{Whitaker2002}
J.~S. Whitaker and T.~M. Hamill.
\newblock Ensemble {{Data Assimilation}} without {{Perturbed Observations}}.
\newblock \emph{Monthly Weather Review}, 130\penalty0 (7):\penalty0 1913--1924,
  July 2002.
\newblock ISSN 0027-0644.

\bibitem[Zhou et~al.(2000)Zhou, Reich, and Brooks]{Zhou2000}
J.~Zhou, S.~Reich, and B.~Brooks.
\newblock Elastic molecular dynamics with self-consistent flexible constraints.
\newblock \emph{J. Chem. Phys.}, 112:\penalty0 7919---7929, 2000.

\end{thebibliography}
